\documentclass[10pt,reqno]{amsart}

\usepackage[utf8]{inputenc}
\usepackage[T1]{fontenc}
\usepackage{mlmodern}

\usepackage{bbm}            

\DeclareFontFamily{U}{wncy}{}
\DeclareFontShape{U}{wncy}{m}{n}{<->wncyr10}{}
\DeclareSymbolFont{mcy}{U}{wncy}{m}{n}
\DeclareMathSymbol{\Sh}{\mathord}{mcy}{"58}
\DeclareMathOperator \KL{D_{KL}\mkern-1mu}

\usepackage{amsmath}
\usepackage{mathtools}     
\usepackage{amsfonts}
\usepackage{amssymb}
\usepackage{amsthm}
\usepackage{mathrsfs}      

\numberwithin{equation}{section}

\DeclareFontFamily{U}{mathx}{} 
\DeclareFontShape{U}{mathx}{m}{n}{<-> mathx10}{}
\DeclareSymbolFont{mathx}{U}{mathx}{m}{n}
\DeclareMathAccent{\widehat}{0}{mathx}{"70}
\DeclareMathAccent{\widecheck}{0}{mathx}{"71}

\usepackage{setspace}
\usepackage{enumitem}
\usepackage{float}
\usepackage{makecell}
\usepackage{tablefootnote}
\usepackage{caption}
\usepackage[labelformat=simple]{subcaption}
\usepackage{csquotes}      

\restylefloat{table}
\newlist{enumlemma}{enumerate}{3} 
\setlist[enumlemma]{label*={ (\alph*)}, ref={(\alph*)}}
\newlist{enumcount}{enumerate}{3} 
\setlist[enumcount]{label*={ (\arabic*)}, ref={(\arabic*)}}

\usepackage[dvipsnames]{xcolor}
\usepackage[normalem]{ulem}

\definecolor{green}{rgb}{0.1, 0.5, 0.1}

\usepackage{marginnote}
\usepackage[colorinlistoftodos]{todonotes}

\usepackage[most]{tcolorbox}
\tcbuselibrary{theorems, skins, breakable}

\tcbset{
    sharp/.style={
        sharp corners,
        breakable
    }
}

\newtcolorbox{sumbox}[1][]{colback=blue!5!white, colframe=blue!75!black, fonttitle=\bfseries, title={Summary: #1}, breakable, sharp corners}
\newtcolorbox{chkbox}[1][]{colback=yellow!10!white, colframe=yellow!50!black, fonttitle=\bfseries, title={Check: #1}, breakable, sharp corners}

\usepackage[colorlinks=true, linkcolor=blue, citecolor=blue, linktocpage=true, urlcolor=black]{hyperref}

\usepackage[alphabetic, abbrev,nobysame,lite]{amsrefs}

\usepackage[capitalize, sort]{cleveref}

\crefformat{section}{\S#2#1#3}

\usepackage[toc, title]{appendix}

\theoremstyle{plain}
\newtheorem{theorem}{Theorem}[section]

\newtheorem{proposition}[theorem]{Proposition} 
\newtheorem{lemma}[theorem]{Lemma}
\newtheorem{corollary}[theorem]{Corollary}

\newtheorem{claim}[theorem]{Claim}
\crefname{claim}{Claim}{Claims}

\theoremstyle{definition}
\newtheorem{definition}[theorem]{Definition}
\newtheorem{remark}[theorem]{Remark}

\newtheorem{example}[theorem]{Example}

\theoremstyle{remark}

\newtheorem*{remark*}{Remark}

\newcommand{\eqnote}[1]{ \hfill {\text{\footnotesize #1}}}

\long\def\symbolfootnote[#1]#2{\begingroup\def\thefootnote{\fnsymbol{footnote}}\footnote[#1]{#2}\endgroup}

\DeclareMathOperator{\diam}{diam}

\DeclareMathOperator{\supp}{supp}

\DeclareMathOperator{\id}{Id}

\DeclareMathOperator{\Emb}{Emb}

\def\Id{\id}

\def\ptpart{\mathcal E}

\def\loc{\mathrm{loc}}

\providecommand{\restrict}{}
\renewcommand{\restrict}[1]{\raisebox{-.0ex}{$\upharpoonright$}_{#1}}

\newcommand{\hol}[1][\beta]{\, \mathrm{H\"{o}l}^{{#1}}}
\def\Hol{\hol}

\newcommand{\inv}{^{-1}}

\renewcommand{\emph}[1]{{\bf #1}} 
\renewcommand{\hat}{\widehat}
\renewcommand{\phi}{\varphi}

\renewcommand{\bar}{\overline}

\newcommand{\diff}{\mathrm{Diff}}
\newcommand{\Diff}{\diff}

\providecommand{\N}{}\renewcommand{\N}{\mathbb{N}}

\providecommand{\R}{}\renewcommand{\R}{\mathbb{R}}
\providecommand{\T}{}\renewcommand{\T}{\mathbb{T}}

\providecommand{\Z}{}\renewcommand{\Z}{\mathbb{Z}}

\providecommand{\Nbb}{}\renewcommand{\Nbb}{\mathbb{N}}
\providecommand{\Rbb}{}\renewcommand{\Rbb}{\mathbb{R}}

\providecommand{\Zbb}{}\renewcommand{\Zbb}{\mathbb{Z}}

\def\calP{\mathcal P}\def\calQ{\mathcal Q}

\newcommand{\Fcal}{\mathcal{F}}

\newcommand{\Qcal}{\mathcal{Q}}

\newcommand{\Pcal}{\mathcal{P}}\newcommand{\Rcal}{\mathcal{R}}

\def\partP{\mathcal P}
\def\partQ{\mathcal Q}
\def\partR{\mathcal R}
\def\partK{\mathcal K}

\def\sigA{\mathscr A}\def\sigB{\mathscr B}\def\sigC{\mathscr C}
\def\sigF{\mathscr F}\def\sigG{\mathscr G}
\newcommand{\Asc}{\mathscr{A}}\newcommand{\Bsc}{\mathscr{B}}
\newcommand{\Csc}{\mathscr{C}}
\newcommand{\Fsc}{\mathscr{F}}

\newcommand{\wtd}{\widetilde}
\newcommand{\what}{\widehat}
\newcommand{\wcheck}{\widecheck}
\newcommand{\td}{\tilde}

\newcommand{\refine}{\mathrm{ref}}

\def\siname{invariant}
\def\maxent{h^{\mathrm{ref}}}
\def\maxFent{\maxent_\mu(T\mid \sigF)}

\def\transent{h^{\mathrm{tr}}}

\def\ptpart{\mathcal E}
\def\trivpart{\mathcal N}
\def\furstent{h^{\mathrm{Furst}}}

\title[Exact dimensionality]{Exact dimensionality of hyperbolic stationary measures}

\author[A.~Brown]{Aaron Brown}
\author[H.~Lee]{Homin Lee}
\address{Northwestern University, Evanston, IL 60208, USA}
\email{awb@northwestern.edu}
\address{HCMC, Korea Institute for Advanced Study, Seoul, South Korea}
\email{hominlee@kias.re.kr}

\def\calP{\mathcal P}

\title[Exact dimensionality]{Exact dimensionality of hyperbolic stationary measures}

\begin{document}

\begin{abstract}
For finitely generated i.i.d.~random walks 
we show all ergodic hyperbolic stationary measures are exact dimensional with dimension the sum of certain stable and unstable dimensions.  
We derive similar results for invariant measures for non-invertible skew-products with application to systems admitting stable holonomies and derive a new invariance principle in terms of pointwise dimension of fiber measures.  
\end{abstract}
\maketitle

\let\oldtocsection=\tocsection
\let\oldtocsubsection=\tocsubsection 
\let\oldtocsubsubsection=\tocsubsubsection
 
\renewcommand{\tocsection}[2]{\hspace{0em}\oldtocsection{#1}{#2}}
\renewcommand{\tocsubsection}[2]{\hspace{1em}\oldtocsubsection{#1}{#2}}
\renewcommand{\tocsubsubsection}[2]{\hspace{2em}\oldtocsubsubsection{#1}{#2}}
\setcounter{tocdepth}{2}

\tableofcontents

\section{Introduction}\label{sec:intro} 

Given a compact manifold and a finitely supported probability measure $\nu_0$ on the group $\diff^{r}(M)$ of $C^r$ diffeomorphisms of $M$ ($r>1$), we study the ergodic $\nu_0$-stationary measures $\mu_0$ for the random walk on $M$ distributed with law $\nu_0$.  
When an ergodic stationary measure $\mu_0$ is hyperbolic (i.e., with all Lyapunov exponents non-zero), we show the stationary measure $\mu_0$ is exact dimensional; moreover the dimension coincides with the sum of certain associated conditional dimensions along stable and unstable manifolds.  This extends to the setting of hyperbolic stationary measures the exact dimensionality results for probability measures invariant under single diffeomorphisms established by  Young \cite{MR684248} on surfaces and by \cite{MR1709302} on higher-dimensional manifolds.  It also extends the results of  \cite{BRH} showing exact dimensionality of ergodic, contracting on average (i.e.~with all Lyapunov exponents negative) stationary measures  to the setting of ergodic hyperbolic stationary measures.  

As a direct consequence of our results, we have the following analogue of Young's theorem from \cite{MR684248} for i.i.d.~random walks on surfaces:  

\begin{corollary}\label{cor:Young}
Let $M$ be a compact surface, let $\nu_0$ be a finitely supported  probability measure  on the group $\diff^{1+\beta}(M)$, and let $\mu_0$ be an ergodic, hyperbolic, $\nu_0$-stationary measure with Lyapunov exponents $\lambda^s<0<\lambda^u$. Then for  $\mu_0$-a.e.~$x\in M$  \[
\lim_{r\to 0} \dfrac{ \log \mu_0(B(x,r))}{\log r} =
\frac{1}{\lambda^{u}} h_{\mu_{0}}(\nu_0)
-\frac{1}{\lambda^{s}}\Bigl(h_{\mu_{0}}(\nu_0)+\furstent_{\mu_{0}}(\nu_0)\Bigr).
\]
In particular,  $\mu_0$ is exact dimensional.  

\end{corollary}

Above $h_{\mu_{0}}(\nu_0)$ is the ``standard'' metric entropy  and ${h_{\mu_{0}}^{\textrm{Furst}}(\nu_0)}$ is the Furstenberg entropy of the stationary measure $\mu_0$ for the i.i.d.~random walk determined by $\nu_0$.  See e.g.~\cite[\S2.2.1]{BRH} and \cite[\S2.2.3]{BRH}   for definitions of the ``standard'' and Furstenberg entropies in the setting of i.i.d.~random walks.  

One might question the asymmetry in \cref{cor:Young} (and our more general \cref{thm:mainiid} below) between the stable and unstable contributions to the dimension of $\mu_0$.  Roughly, this is due to an asymmetry in entropy formulas associated with stable and unstable manifolds and the measure $\mu_0$.  When stable and unstable manifolds are 1-dimensional, the  product of the dimension of the measure $\mu_0$ along stable manifolds and $-\lambda^s$ gives a quantity that in \cite{BRH} was called ``transverse  entropy''; in this setting, the ``transverse entropy'' is given by $h_{\mu_{0}}(\nu_0)+\furstent_{\mu_{0}}(\nu_0). $  
On the other hand, the  product of the dimension of the measure $\mu_0$ along unstable manifolds and the positive Lyapunov exponent gives the standard  entropy $h_{\mu_{0}}(\nu_0)$.  The (conditional) measures of Riemannian balls decay along stable manifolds like ``transverse entropy'' whereas along unstable manifolds, the (conditional) measures of Riemannian balls decay like standard entropy; the total decay of Riemannian balls is the sum of these two decay rates.  
See \cref{sss:asymmetry} below for more detailed discussion.  

For the remainder of this introduction, we define our main objects of study and state our general theorems in the setting of countably generated, i.i.d.~random walks.  In \cref{sec:general} we will present our results in a more general framework of  skew products which includes results for fibered skew product systems admitting systems of stable holonomies.

\subsection{Stationary measures and  skew extensions for i.i.d.\ random walks}  
Let $M$ be a compact manifold (or a more general compact metric space).  A distribution $\nu_0$ on the space of continuous self-maps $C^0(M,M)$ generates an i.i.d.\ random walk on $M$.  
We say a Borel probability measure $\mu_0$ on $M$ is \emph{$\nu_0$-stationary} if $$\int f_*\mu_0 \, d \nu_0(f) = \mu_0.$$ 
By compactness and continuity, stationary measures exist.  

In our  applications of interest, we will be concerned with the case $M$ is a smooth closed manifold and that $\nu_0$ is supported on $\diff^r(M)$, the group of $C^r$ diffeomorphisms of $M$, for some $r>1$.   We will further typically assume that $\nu_0$ is either finitely supported or has countable  support and satisfies certain moment and entropy conditions (see \cref{sec:hypiid}).

\subsubsection{Coding random compositions} 
Let $C = \supp\{\nu_0\}$ and write  $\Sigma:= C^\Z$ and  $\widehat \Sigma:= C^\N$.  Equip both $\Sigma$ and $\what \Sigma$ with the left shift $\sigma$.  
Given $\omega\in \widehat \Sigma$, write $\omega = (f_0, f_1, f_2, \dots)$.  For $n\ge 1$, we write $$f_\omega^{(n)}:=f_{n-1}\circ \dots \circ f_{0}$$
for the $n$-fold iterate given by the word $\omega$.  In particular, $f_\omega:=f_\omega^{(1)}= f_0$.  We define $f_\omega^{(0)}$ to be the identity map.  
Given $\omega\in \Sigma$  write $\omega = (\dots f_{-2}, f_{-1},f_0, f_1, f_2, \dots)$.  Similarly define $f_\omega^{(n)}$ for $n\ge 0$; for $n<0$, define 
$$f_\omega^{(n)} := f_{n}^{-1}\circ \dots  \circ f_{-1}^{-1}.$$

\subsubsection{Skew extensions and sample measures}\label{sec:skewext}
Let $\what F\colon \what \Sigma\times M\to  \what \Sigma\times M $ and $F\colon \Sigma\times M \to  \Sigma\times M$ denote the skew-product dynamical systems  defined by  $$(\omega, x) \mapsto (\sigma(\omega), f_\omega(x)).$$  
It is well known that a probability measure $\mu_0$ on $M$ is $\nu_0$-stationary if and only if the product measure ${\nu_0}^\N \times \mu_0$ on $\what \Sigma \times M$ is $\what F$-invariant.  We say $\mu_0$ is ergodic as a  $\nu_0$-stationary measure if and only if ${\nu_0}^\N \times \mu_0$ is ergodic as a $\what F$-invariant measure.  

The product measure ${\nu_0}^\Z \times \mu_0$ on $ \Sigma \times M$ need not be $F$-invariant.  However, the limit $$\td \mu = \lim_{n\to \infty} F^n_*({\nu_0}^\Z \times \mu_0)$$
exists in the (compact) space of probability measures on $\Sigma\times M$ that project to $\nu_0^\Z$  (see e.g.\ \cite[Prop.~I.1.2.]{MR1369243}); by construction $\td \mu$ is $F$-invariant and projects to $\nu_0^\Z$.  
The measure $\td \mu$ is   further characterized as follows: $\td \mu$ is the  unique $F$-invariant probability measure  on $\Sigma\times M$ whose projection to 
$\widehat \Sigma\times M$ is 
$\hat\mu={\nu_0^\N}\times {\mu_0}$ (c.f.\ \cite[Prop.\  I.1.2]{MR1369243}).  In particular,   
 \begin{enumerate}
 \item  $F\colon (\Sigma\times M , \td \mu) \to (\Sigma\times M, \td \mu)$ is the (ergodic-theoretic) natural extension of $\what F\colon (\widehat \Sigma\times M , \hat \mu) \to (\widehat \Sigma\times M , \hat  \mu),$ and 
 \item $\td \mu$ is ergodic if and only if $\mu_0$ is ergodic.  
 \end{enumerate}

We let $\{\td \mu_\omega\}$ denote the (unique mod $\nu_0^\Z$) family of conditional measures relative to the partition of $\Sigma\times M$ into fibers of the projection $\Sigma\times M\to \Sigma$.  In the literature, the measures $\{\td \mu_\omega\}$ are often referred to as the \emph{sample measures} of the random dynamical system.  
By identifying $M$ with $\{\omega\}\times M$, we abuse notation and view each $\td\mu_\omega$  either as a probability measure on $M$ or as a probability measure on the fiber $\{\omega\}\times M$ of  $\Sigma\times M$.  

The natural extension $\td \mu$ associated to a stationary measure $\mu_0$ has the following additional geometric property: the sample measures $\omega\mapsto \td \mu_{\omega}$ are independent of the future.  Namely, given $\omega\in \Sigma$, write $$\Sigma^u_{\mathrm{loc}}(\omega) :=\{\omega'\in \Sigma \text{ with } \omega_k=\omega'_k \text{ for all }  k< 0\}$$
equipped with the measure $\nu_0^\N$ (on non-negatively indexed terms). 
Then: \begin{equation}\label{eq:statprop}
\text{\itshape
for $\nu_0^\Z$-a.e.\ $\omega\in \Sigma$ and $\nu_0^\N$-a.e.\ $\omega'\in \Sigma^u_{\mathrm{loc}}(\omega),$ we have  $\td \mu_\omega = \td \mu_{\omega'}$.}
\end{equation}
\subsection{Smooth ergodic theory and pointwise dimensions}
We now emphasize that $M$ is assumed to be a closed manifold equipped with a Riemannian metric.  Let $\nu_0$ be supported on $\diff^r(M)$ for some $r>1$.  Let $\mu_0$ be an ergodic, $\nu_0$-stationary probability measure.  

\subsubsection{Random Oseledec's theorem}
Throughout we will assume the maps $f\mapsto \log \|f\|_{C^r}$ and $\log \|f\inv\|_{C^r}$ are $L^1({\nu_0})$.  
We thus may apply   Oseledec's theorem (c.f.\ \cite[Thm.\ III.1.1]{MR884892} for a formulation in the setting of random dynamics): 
given ${\nu_0^\N}$-a.e.\ $\omega$, ${\mu_0}$-a.e.\ $x$, and every $0\neq v\in T_xM$, 
the limit $$\lim_{n\to \infty} \frac 1 n \log \|D_xf^{(n)}_\omega v\|$$ exists and is finite.  As $T_xM$ is finite dimensional and as $\mu_0$ was assumed ergodic, this limit takes finitely many values 
$$\lambda_1<\lambda_2<\dots<\lambda_{p}$$ 
that are independent (mod $\mu_0$) of the choice of $x$ or (mod $\nu_0^\N$) of the choice of $\omega$.  
The numbers  $\lambda_1<\lambda_2<\dots<\lambda_{p}$ are the \emph{Lyapunov exponents} of $\mu_0$.  

We let $s = \max \bigl(\{0\} \cup \{i: \lambda_i<0\} \bigr )$.   
As we focus on the case of hyperbolic stationary measures (with exponents of each sign), we will typically assume $1\le s<p$.  

\subsubsection{Random Pesin manifolds} 
Given $\omega\in \what \Sigma$ and $x\in M$, let 
$$W^s_{\omega}(x):=\left\{y\in   M:  \limsup_{n\to \infty}  \frac 1 n \log d_M( f_{ \omega}^{(n)}(x),f_{ \omega}^{(n)}(y)) <0 \right\}.$$
Similarly,  given $\omega\in  \Sigma$ and $x\in M$, let 
$$W^u_{\omega}(x):=\left\{y\in   M:  \limsup_{n\to \infty}  \frac 1 n \log d_M( f_{ \omega}^{(-n)}(x),f_{ \omega}^{(-n)}(y)) <0 \right\}.$$
Applying a fibered version of Pesin's stable manifold theorem (c.f.\ \cite[Thm.\  III.3.2]{MR1369243} for the case of $C^{1+\mathrm{Lip}}
$ fibered systems or discussion and references in \cite{MR4599404} for the setting of $C^{1+\text{H\"older}}$-fibered systems) 
we obtain for  $\nu_0^\N\times \mu_0$-a.e.\ $(\omega,x)\in \what \Sigma\times M $ that $W^s_{\omega}(x)$ is a connected, injectively immersed, $C^r$ submanifold of $M$.  Similarly, for  
$\td \mu$-a.e.\ $(\omega,x)\in  \Sigma\times M $, the set $W^u_{\omega}(x)$ is a connected, injectively immersed, $C^r$ submanifold of $M$.  

We also write $$\text{$W^s(\omega,x):= \{\omega\}\times W^s_\omega(x)$ and $W^u(\omega,x):= \{\omega\}\times W^u_\omega(x)$}$$ for the corresponding \emph{fiberwise stable} and \emph{fiberwise unstable} subsets of the   skew extensions.  We note that $W^s(\omega,x)$ is well-defined for both $(\omega,x)\in \what \Sigma\times M$ and $ (\omega,x)\in  \Sigma\times M$.  On the other hand, $W^u(\omega,x)$ is defined only for  $ (\omega,x)\in  \Sigma\times M$.

\subsubsection{Stable and unstable dimensions}
Given $(\omega,x)\in \what \Sigma\times M$, a $(\nu_0^\N\times\mu_0)$-measurable partition $\hat \xi^s$ of $\what \Sigma\times M$ is \emph{subordinate to fiberwise $W^s$-manifolds} if for $(\nu_0^\N\times\mu_0)$-a.e.\ $(\omega, x)\in \what \Sigma\times M$, the atom  $\hat \xi^s(\omega, x)$ is a precompact subset of $W^s(\omega,x)$ that  contains an open neighborhood of $x$ in $W^s(\omega,x)$.  

We may fix $\omega\in \what \Sigma$ and instead consider the partition $\hat\xi^s_\omega$ given by restricting $\hat \xi^s$ to  $\{\omega\}\times M$.  
For $\nu_0^\N$-a.e.~$\omega$, the partition $\hat\xi^s_\omega$ is a $\mu_0$-measurable partition of $M$ with the property that for $\mu_0$-a.e.~$x$, the atom $\hat\xi^s_\omega(x)$ is precompact in $W^s_\omega(x)$ and contains an open (in immersed topology) neighborhood of $x$ in $W^s_\omega(x)$.

We similarly define a measurable partition $ \bar \xi^u$ of $\Sigma\times M$ to be  \emph{subordinate to fiberwise $W^u$-manifolds} if  for $\td \mu$-a.e.\ $(\omega, x)\in  \Sigma\times M$, the atom  $ \bar \xi^u(\omega,x)$ is a precompact subset of $W^u(\omega,x)$ that  contains an open neighborhood of $(\omega,x)$ in $W^u(\omega,x)$.

Given $x\in M$, write $B(x,r)$ for the Riemannian ball of radius $r$ centered at $x$ in $M$.  Given $\omega\in\what  \Sigma$ or $\omega\in \Sigma$, let $B\left((\omega,x),r\right):= \{\omega\}\times B(x,r)$ denote the corresponding Riemannian ball in the fiber over $\omega$.
 By \cite[Prop.~2.5]{BRH} \cite{MR2032494} and ergodicity, there are non-negative quantities  $\hat d^s(\mu_0)$ and $\bar d^u(\mu_0)$, such that for $\nu_0^\N\times \mu_0$-a.e.~$(\omega,x)\in \what \Sigma\times M$ and $\td \mu$-a.e.~$(\omega',x')\in  \Sigma\times M$, the following limits hold:
\begin{equation*}\hat d^s(\mu_0):= \lim_{r\to 0} \dfrac {\log \left((\mu_0)^{(\hat \xi^s_{\omega})_x}(B(x,r))\right)}{\log r}\end{equation*}
\begin{equation} \label{eq:dbarstatmsr} \bar d^u( \mu_0) := \lim_{r\to 0} \dfrac {\log \left(\td \mu^{ \bar \xi^u}_{(\omega',x')}\left(B\left((\omega',x'),r\right)\right)\right)}{\log r}
\end{equation}

Because  $\td \mu^{ \bar \xi^u}_{(\omega',x')}$ is supported on $\{\omega'\}\times M$, we have, alternatively, 
\begin{equation*}  \bar d^u( \mu_0) := \lim_{r\to 0} \dfrac {\log \left(\td \mu^{ \bar \xi^u}_{(\omega',x')}\left(\Sigma\times B\left(x',r\right)\right)\right)}{\log r}.
\end{equation*}
Recall that $\hat \xi^s$ is a measurable partition of $\hat \Sigma\times M$ and that $\tau\colon \Sigma\times M\to\hat \Sigma\times M$ is the natural factor map.  We write $\xi^s= \tau\inv \hat \xi^s$.  Then atoms of $\xi^s$ are of the form $$\xi^s(\omega,x):= \{(\omega',y): \text{$y\in \hat \xi(\tau(\omega),x) $ and 
$\omega_i'=\omega_i$ for all $i\ge 0$}
\}.$$
Since $\tau_*(\td \mu)^{ \xi^s}_{(\omega,x)}= (\hat \mu)^{ \hat \xi^s}_{(\tau(\omega),x)}$ a.s., we have, alternatively, that for $\td \mu$-a.e.~$(\omega,x)$,
\begin{equation}\label{eqhatstatmsr}\hat d^s(\mu_0):=\lim_{r\to 0} \dfrac {\log \left(\td \mu^{ \xi^s}_{(\omega,x)}(\Sigma\times B(x,r))\right)}{\log r}.\end{equation}
Note the asymmetry in the geometry of the partitions $\overline \xi^u$  and $\xi^s$ in \eqref{eq:dbarstatmsr} and  \eqref{eqhatstatmsr} is intentional.

\subsection{Statement of results for i.i.d.~random walks}
\subsubsection{Hypotheses for i.i.d.~random walks}\label{sec:hypiid}
To move slightly outside the setting of random walks determined by finitely supported measures $\nu_0$, we will assume our measure $\nu_0$ satisfies the following properties:
\begin{enumerate}
\item (Discreteness) ${\nu_0}$ is a discrete measure; namely there is a countable set $C\subset \diff^{1+\beta}(M)$ with ${\nu_0}(C)=1$.
\item(Finite entropy)  $H({\nu_0}):= -\sum_{f\in C} {\nu_0}(f) \log {\nu_0}(f) <\infty$.
\item\label{inthyp} (Finite $\log$-moment) The assignments 
$f\mapsto \log  \|f\|_{C^{1+\beta}}$ and $f\mapsto \log  \|f\inv\|_{C^{1+\beta}}$ are $L^1({\nu_0})$.
\end{enumerate}

\subsubsection{Main result for  i.i.d.\ random walks}\label{sec:iid}
Throughout, we will always assume $\mu_0$ is ergodic and hyperbolic.  In particular,  the  Lyapunov exponents $$\lambda_1<\lambda_2<\dots<\lambda_{p}$$ listed without multiplicity are a.s.~independent of $(\omega,x)$; moreover, no Lyapunov exponent is zero.

Recall we let $s= \max \bigl(\{0\} \cup \{i: \lambda_i<0\} \bigr )$.  When $s=0$ all exponents are positive.  The classical invariance principle of Crauel \cite[Thm.~5.1]{MR1204858} (see also \cite{MR2651382}) implies $f_*\mu_0 = \mu_0$ for every $f$ in the support of $\nu_0$.  In particular, the measure $\td \mu$ in \cref{sec:skewext} is of the form $\td \mu = \nu_0^\Z\times \mu_0$.  In particular, the sample measure $\td\mu_\omega$ coincides (mod $\nu_0$) with $\mu_0$.  Since $\mu_0$ has all  Lyapunov exponents positive, a standard exercise shows that sample measures  $\td\mu_\omega$ are finitely many atoms (for a.e.~$\omega$).  Thus   the exact dimensionality, dimension formula (with dimension 0), and almost product structure follow immediately.   We thus assume $1\le s$.  On the other hand, when $s= p$, all Lyapunov exponents are negative; exact dimensionality of the stationary measure $\mu_0$ and the relationship between dimension, Lyapunov exponents, and (Furstenberg) entropy  follow from the main results (c.f.~\cite[Cor.~2.7]{BRH}) of \cite{BRH}.  

Thus, for the remainder we will assume  that $1\le s<p$.

\begin{theorem}[Exact dimensionality of hyperbolic stationary measures] \label{thm:mainiid} 
For $r>1$, let $\nu_0$ be a probability measure on $\diff^{r}(M)$ satisfying the hypotheses of \cref{sec:hypiid}.  
Let $\mu_0$ be an ergodic, hyperbolic $\nu_0$-stationary measure with $1\le s<p$.  
Then for $\mu_0$-a.e.~$x\in M$, 
\begin{equation}
\label{eq:exacdimiid}
\lim_{r\to 0} \frac{\log \mu_0(B(x,r))}{\log r}= \hat d^s(\mu_0)+  \bar d^u(\mu_0).\end{equation}
In particular, 
$\mu_0$ is exact dimensional and $\dim(\mu_0)=\hat d^s(\mu_0)+  \bar d^u(\mu_0).$
\end{theorem}
\subsubsection{Almost product structure}\label{sss:APS}
The exact dimensionality  in \cref{eq:exacdimiid} combined with the formula in terms of stable and unstable dimensions implies that hyperbolic stationary measures exhibit a certain ``almost-product'' structure.

\begin{corollary}\label{cor:iidaps}
In the setting of \cref{thm:mainiid}, for every $\delta>0$, $\td \mu$-a.e.~$(\omega,x)\in \Sigma\times M$, and all sufficiently  (depending on $(\omega,x)$) small $r>0$, 
\begin{equation*}
\begin{split}
r^{\delta}\td \mu^{\overline \xi^u}_{(\omega,x)} \bigl(B(x,r)\bigr) & \hat \mu^{\hat \xi^s}_{(\tau(\omega),x)} \bigl(B(x,r)\bigr)
 \le \td \mu\bigl(\Sigma\times B(x,r)\bigr) \\
 &\quad  \quad  \le r^{-\delta} \td \mu^{\overline \xi^u}_{(\omega,x)} \bigl(B(x,r)\bigr) \hat \mu^{ \hat \xi^s}_{(\tau(\omega),x)} \bigl(B(x,r)\bigr)
\end{split}
\end{equation*}
\end{corollary}

\subsection{Discussion of asymmetry in \cref{thm:mainiid,cor:iidaps,cor:Young}}\label{sss:asymmetry}
The asymmetry in definitions of the dimensions $\hat d^s(\mu_0)$ and $  d^u(\mu_0)$  in \cref{thm:mainiid,cor:iidaps} may seem mysterious at first.  

We justify this as follows.  A $\nu_0$-stationary measure $\mu_0$ naturally gives rise to a non-invertible skew extension $\what F\colon \what \Sigma \times M\to  \what \Sigma \times M$ preserving the ergodic probability measure $\nu_0 ^\N\times \mu_0$.  
We expect decay of conditional measures along both stable and unstable manifolds to contribute to the dimension of the stationary measure  $\mu_0$; since unstable manifolds require invertible dynamics to define, it becomes natural to pass to the two-sided skew shift $F\colon \Sigma \times M\to \Sigma \times M$ preserving the ergodic probability measure $\td \mu$ (the natural extension of $\what F$).  

Let $\tau \colon \Sigma\to \what \Sigma$ be the natural factor map.  In $\Sigma\times M$, let $\sigC$ denote the $\sigma$-algebra  generated by  preimages under the map $\Sigma\times M \to \what \Sigma$; then atoms of $\sigC$ are sets of the form $$
\tau\inv (\tau(\omega))\times M.$$
Similarly, let $\wtd\sigA$ denote the $\sigma$-algebra  generated by  preimages under the map $\Sigma\times M \to \what \Sigma\times M$.

The $\mu_0$-measure of a  ball $B(x,r)$ in $M$ corresponds to the $\td \mu^\sigC_{(\omega,x)}$-measure of the saturated ball $\tau\inv \{\tau(\omega)\}\times B(x,r)$ for $\td \mu$-a.e.~$(\omega,x)\in \Sigma\times M$.
Atoms of $\sigC$ are saturated by combinatorial local stable manifolds (atoms of $\wtd \sigA)$ but are not saturated by combinatorial local unstable manifolds.  Thus relative to unstable manifolds, the conditional $\td \mu^\sigC_{(\omega,x)}$-measure of the saturated ball $\tau^{-1}\{\tau(\omega)\} \times B(x,r)$ decays like the unstable decay of sample measures; this is precisely the standard conditional entropy.   

On the other hand, relative to stable manifolds, the conditional $ {\td \mu}^\sigC_{(\omega,x)}$-measure of the saturated ball $\tau^{-1}\{\tau(\omega)\}\times B(x,r)$ decays like the projection of the sample measures.  This is precisely the difference of refinement entropy $\maxent_{\td \mu}(F\mid \sigC)$ of the two-sided skew system (see \cref{Def:refent} below for definition) and the standard entropy $h_{\td \mu}(F,\Xi(\wtd \sigA))$ relative to $\Xi(\wtd \sigA)$; in the language of \cite{BRH}, this difference is precisely the transverse entropy $\transent_{ \hat  \mu}(\what F\mid \what \sigF)$ of the one-sided skew system.  The difference in entropies controlling decay of measure along unstable vs stable manifolds explains the asymmetry in the dimensions contributing to the total dimension of $\mu_0$.

We note that while the stable contribution to the dimension of $\mu_0$ is best explained in terms of transverse entropy 
$\transent_{ \hat  \mu}(\what F\mid \what \sigF)$ 
 of the one-sided skew system, we do not ever  use the transverse entropy in this paper.  

\subsection*{Use of AI}
ChatGPT supplied the ideas for the Anosov example (Example~\ref{ex:anosov-graph}) and both circle examples 
(Examples~\ref{ex:dimform} and~\ref{ex:attracting-two-cycle}).  
It also suggested the limit \eqref{eq:limitvsliminf} to obtain stronger conclusions (those currently stated) in  \cref{prop:entcal}\ref{entcal2}--\ref{entcal4}.  

We used ChatGPT to understand and check certain arguments in \cite{MR819557,MR1709302,MR2218998} and to audit and check our arguments, proofread, and check consistency of notation and cross-references, and to format some latex.  

The formulation of the results, new arguments, and adaptation of existing cited literature are due to the authors.  The examples were carefully edited by the authors; all other writing was done by the authors.

\subsection*{Acknowledgements} 
AB was partially supported by NSF Grants DMS-1752675 and DMS-2400191.
HL was supported in part by AMS-Simons travel grant, a KIAS Individual Grant (HP104101) via the June E Huh Center for Mathematical Challenges at Korea Institute for
Advanced Study, and Sanghyun Kim's  Mid-Career Researcher Program (RS-2023-00278510) through the National Research Foundation funded by the government of Korea.

Both authors benefited from discussions with Federico Rodriguez Hertz and Lin Shu and are grateful for hospitality from IHP,  KIAS, and NU~for hosting the authors at various periods of this collaboration.

\section{Formulation of results for general skew-products} \label{sec:general}
The result in \cref{thm:mainiid} is new even in the setting of finitely generated, i.i.d. random walks.  However, it becomes natural to generalize the setting and consider more general models of random dynamics that allow for non-independent (though still identically distributed) random walks.  This generalization also allows for the study of fibered skew product dynamics admitting systems of stable holonomies.

\subsection{Model of randomness and associated skew-product(s)}\label{sec:skewsyst}
 In this section, we define notation that will be used throughout the paper.
\subsubsection{Measurable models of driving dynamics}

Let $(\widehat{\Omega}, \widehat{\Bsc}, \widehat{\nu})$ be a standard probability space and let  $\widehat{\theta}\colon\widehat{\Omega}\to \widehat{\Omega}$ be a measurable, $\widehat{\nu}$-preserving transformation.  We will assume throughout that $\widehat{\nu}$ is ergodic.  Moreover, our results are primarily of interest when  $\widehat{\theta}$ is non-invertible.   (Indeed, when $\widehat{\theta}$ is invertible and our skew extensions are assumed $C^2$, our   results essentially  follow  from \cite{MR2218998}.)

Denote by  $\tau\colon (\Omega, \Bsc, \nu)\to (\widehat{\Omega},\widehat{\Bsc},\widehat{\nu})$  the natural extension of $\widehat{\theta}\colon\widehat{\Omega}\to \widehat{\Omega}$ equipped with the ergodic, $\nu$-preserving, invertible transformation $\theta\colon\Omega\to \Omega$ factoring (via $\tau$) onto $\hat \theta$.  
Let $\Asc\subset \Bsc$ denote the ($\theta$-)sub-invariant sub-$\sigma$-algebra $\Asc=\tau^{-1}\widehat{\Bsc}$. 
By uniqueness of the natural extension, we have  $$\Asc^{+}:=\bigvee_{i=-\infty}^{0} \theta^{-i}\Asc=\Bsc.$$

We will often  abuse  notation and use $\omega$ to represent  an element of $\Omega$ or $\what\Omega$.  

\subsubsection{Backward base dynamical system} 
Recall that a  \emph{one-sided generator} $\what\Qcal$ for  $(\what\Omega, \what\Bsc, \what\nu,\what\theta)$ is a measurable partition $\what \Qcal$ with  finite entropy, $H_{\what\nu}(\what\Qcal)<\infty$, such that $$\what \Qcal ^-:=\bigvee_{j=0}^{\infty}\what\theta^{-j}\what\Qcal =\widehat \ptpart$$
where  $\widehat \ptpart$ denotes the point partition on $\widehat{\Omega}$.
Note that if $\what \Qcal$ is a one-sided generator then $$h_{\hat \nu} (\hat \theta) = h_{\hat \nu} (\hat \theta, \what \Qcal)\le  H_{\what\nu}(\what\Qcal)<\infty.$$

When we assume that $(\widehat{\Omega}, \widehat{\Bsc}, \widehat{\nu},\what\theta)$  admits a one-sided generator $\what\Qcal$, the inverse map of the natural extension admits a natural factor dynamical system whose construction we outline: Let $\Qcal=\tau\inv\what\Qcal$ be the partition of $\Omega$ by preimages  of $\what \Qcal$ under $\tau$. 
Construct the  sub-$\sigma$-algebra  $\sigA'$ of  $\Bsc$ by  $$\Asc'=\sigma\left(\bigvee_{j=1}^{\infty}\theta^{j}\Qcal\right).$$
Let  $(\widecheck\Omega,\widecheck\Bsc,\widecheck\nu)$ denote the quotient measure space of $( \Omega, \Bsc, \nu)$ under $\Asc'$.  
Since $\theta \Asc' =\left(\theta\inv\right)\inv \Asc'\subset \Asc'$, 
$\theta^{-1}$ induces a measure-preserving dynamical system $$\widecheck\theta\colon(\widecheck\Omega,\widecheck\Bsc,\widecheck\nu)\to (\widecheck\Omega,\widecheck\Bsc,\widecheck\nu)$$ as a factor measure-preserving system under $\widecheck \tau$.  
Note that the factor system $(\widecheck\Omega,\widecheck\Bsc,\widecheck\nu,\widecheck\theta)$ has a one-sided generator $\widecheck\Qcal=\widecheck\tau\theta\Qcal$ satisfying $\widecheck\tau\inv\widecheck\Qcal=\theta\Qcal$.

\subsubsection{Random dynamics and induced skew products} Let $M$ be a closed smooth manifold equipped with its Borel $\sigma$-algebra $\Bsc_{M}$.  Fix $\beta>0$ and 
let $\widehat{\Omega}\to \Diff^{1+\beta}(M)$, $\omega \mapsto  f_{\omega}$,  be a $\what\Bsc$-measurable function.  (Here, we view $ \Diff^{1+\beta}(M)$ as a subset of the Polish space of $C^1(M,M)$.)

Write $X=\Omega\times M$ and $\widehat{X}=\widehat{\Omega}\times M$. We abuse notation and also write  $\tau\colon X\to \what{X}$, $\tau\colon (\omega,x)\mapsto (\tau(\omega),x)$,  for the natural projection.  This induces the $\Asc$-measurable map $\Omega \to \Diff^{1+\beta}(M)$, $\omega \mapsto f_\omega :=f_{\tau(\omega)}$.
Let $F\colon X\to X$ and $\widehat{F}\colon\widehat{X}\to \widehat{X}$ be the skew products given by $F(\omega,x)=\left(\theta(\omega), f_{\omega}(x)\right)$ and  $\what F(\omega,x)=\left(\hat \theta(\omega), f_{\omega}(x)\right)$.    We view $\what F$ as an identically distributed (rel.~$\hat \nu$) but  non-independent random dynamical system generated by $\omega \mapsto  f_{\omega}$.  (Alternatively, we view $\what F$ as the quotient  under stable holonomies of a fibered dynamical system $F\colon \Omega\times M\to \Omega\times M$ admitting holonomies along local stable manifolds represented by atoms of $\sigA$.)

Given $\omega\in \Omega$ and $n\in \Z$, we write \begin{equation*} f_{\omega}^{(n)}= \begin{cases}
f_{\theta^{n-1}(\omega)}\circ \dots \circ f_{\omega} & n\ge 1\\
\id & n=0\\
\bigl(f_{\theta \inv \omega}\circ \dots \circ f_{\theta^n(\omega)}\bigr)\inv &
n<0\\
\end{cases}.\end{equation*}
Then, given $n\in \Z$, we have  $F^n(\omega,x) = \bigl(\theta^n(\omega), f^{(n)}_\omega (x) \bigr)$. 
We use similar conventions for $f^{(n)}_{  \omega}$ for $  \omega\in \widehat \Omega$ and $n\ge 0$.  

 We denote by $\pi_1 \colon X\to \Omega$,  $\widehat{\pi}_1\colon \widehat{X}\to \widehat{\Omega}$, 
 $\pi_2 \colon X\to M$, and  $\widehat{\pi}_2\colon \widehat{X}\to M$ the natural projections.

We define the following sub-$\sigma$-algebras on $X$ and $\what X$.  
\begin{enumerate}
\item On $X$, $\sigF= \pi_1\inv \sigB$ is  the sub-$\sigma$-algebra of sets saturated by fibers of $\pi_1$.  Each atom of $\sigF$ is of the form $\{\omega \}\times M$. Similarly, on $\what X$, let $\widehat \sigF= \hat \pi_1\inv \what \sigB$.  
\item $\sigC= \pi_1\inv \sigA = \tau\inv \what \sigF$ is the sub-$\sigma$-algebra on $X$ of sets formed by  saturating sets in $\sigA$ by fibers of $\pi_1$.   Note that  $\sigF= \sigC^+$.  
\item Let $\wtd \sigA$ denote the product $\sigma$-algebra $\sigA\otimes \sigB_M$.  
\end{enumerate}

When $(\what\Omega,\what\Bsc, \what \nu, \what\theta)$ has a one-sided generator $\what \Qcal$, we let $\widetilde\Qcal \coloneqq \pi_{1}^{-1}\Qcal$. Observe that $\widetilde{\Qcal}^{\Zbb} =\Xi(\Fsc)$ and $\widetilde\Qcal^{-}=\Xi(\Csc)$.

\subsubsection{Backward skew product}
Suppose now that $(\widehat{\Omega}, \widehat{\Bsc}, \widehat{\nu}, \what\theta)$  admits a one-sided generator $\what\Qcal$.

Recall our random dynamics is generated by a measurable map  $\widehat{\Omega}\to \Diff^{1+\beta}(M)$, $\omega \mapsto  f_{\omega}$.  At times, we will make the following additional assumption:  
\begin{equation}\label{Qmsble} \text{\itshape $\widehat{\Omega}\to \Diff^{1+\beta}(M)$, $\omega \mapsto  f_{\omega}$ is $\sigma(\what\Qcal)$-measurable.}  \end{equation}

This implies the induced map $\Omega\to \Diff^{1+\beta}(M)$, $\omega\mapsto f_{\omega}$, is $\sigma(\Qcal)$-measurable. 
Define  $\Omega\to\Diff^{1+\beta}(M)$, $\omega\mapsto g_\omega$ where $$g_\omega \coloneqq f_{\theta\inv\omega}\inv.$$
Then  $\omega\mapsto g_{\omega}$ is $\sigma(\theta\Qcal)$-measurable. 
Since $$\Asc'=\sigma\left(\bigvee_{j=1}^{\infty}\theta^{j}\Qcal\right)=
\sigma(\theta \Qcal^+)\supset \sigma(\theta\Qcal),$$
the map $\omega\mapsto g_\omega$ is $\Asc'$-measurable and thus  factors through a map $\wcheck\Omega \to \Diff^{1+\beta}(M)$, $\wcheck\omega \to g_{\wcheck\omega}$.  
Under assumption \eqref{Qmsble}, we define the \emph{backward skew product}: let $\widecheck F \colon \widecheck X\to \widecheck X$ be $\widecheck F (\widecheck\omega,x)= (\widecheck\theta\widecheck\omega, g_{\wcheck\omega} (x))$.
We have a natural factor map $\widecheck\tau\colon X \to \widecheck X$, satisfying  $\widecheck\tau\circ F\inv(\omega,x)=\widecheck F \circ\widecheck\tau (\omega,x)$. We also write natural projection maps as $\wcheck\pi_{1}:\wcheck X \to \wcheck \Omega$ and $\wcheck\pi_{2}:\wcheck X \to M$. At times, we also denote the map $(\omega,x)\mapsto (\wcheck\tau(\omega),x)$ by $\wcheck\tau$ if there is no confusion. 

Define the following $\sigma$-algebras: 
\begin{enumerate}
\item  Let $\wcheck \Fsc =\wcheck\pi_{1}^{-1} \wcheck \Bsc$.
\item Let $\Csc' = \pi_{1}^{-1} \Asc' = \wcheck\tau^{-1}\wcheck \Fsc$. Note that $\Fsc= \Csc'^{-}$. 
\item Let $\widetilde{\Asc'}=\Asc'\otimes \Bsc_{M}$. 
\end{enumerate} 
On $X$, we also define $\wtd{{\Qcal}'} = \pi_{1}^{-1}\wcheck\Qcal$. Note that $\wtd{\Qcal'}^{\Zbb}=\Xi(\Fsc)$ and $\wtd{\Qcal'}^{+}=\Xi(\Csc')$.

\subsubsection{Invariant measures and disintegrations} \label{natext}
Throughout, we let  $\hat \mu$ be an ergodic, $\what F$-invariant probability measure on $\what X$ with $(\hat \pi_1)_*\hat \mu=\hat \nu$ where $\hat \pi_1\colon\what  \Omega\times M\to \what \Omega$ is the projection.

Recall $\tau\colon   \Omega\times M\to \what \Omega\times M$, $\tau(\omega,x)= (\tau(\omega),x)$ is the projection.  
Using that $(\Omega,\nu)$ is the natural extension of $(\what \Omega, \hat \nu)$, there exists a unique $F$-invariant probability measure $\mu$ on $ \Omega\times M$ with $\tau_*\mu = \hat \mu$; moreover, $(\pi_1)_*\mu=\nu$.  
Similarly, given an $F$-invariant probability measure $\mu$  on $X$ with $(\pi_1)_*\mu=\nu$, the measure $\hat \mu = \tau_*\mu$ is  $\what F$-invariant probability measure on $\what X$ with $(\hat \pi_1)_*\hat \mu=\hat \nu$.  Moreover, $\mu$ is ergodic if and only if $\hat \mu$ is.

Given  an ergodic, $F$-invariant probability measure $\mu$ on $X$ with $(\pi_1)_*\mu=\nu$,
we let $\{\mu_\omega\}$ denote a family of conditional measures relative to the fibers of the projection $X\to \Omega.$
As above, we refer to $\{\mu_\omega\}$ as the \emph{sample measures} of $\mu$.

Suppose that   $\what\Qcal$ is a one-sided generator for   $(\what\Omega, \what\Bsc, \what\nu,\what\theta)$ and if  $\omega\mapsto f_\omega$ is $\sigma(\what\Qcal)$-measurable.  Given  an ergodic, $F$-invariant probability measure $\mu$ on $X$ with $(\pi_1)_*\mu=\nu$, we write $\check \mu = \check \tau_* \mu$.  Then $\check \tau_* \mu$ is an ergodic, 
$\widecheck F$-invariant probability measure on $\widecheck X$ with $(\check \pi_1 )_* \check\mu  = \check \nu.$

\subsubsection{Measurability relative to a one-sided generator implies local constantness}\label{sec:reduction}
Although  we will not use this structure explicitly, we note that analysis of the skew product structure  in the setting of \eqref{Qmsble} reduces to the case of a locally constant skew product over a countable full shift.  Indeed, suppose that $\what\Qcal$ is a one-sided generator for   $(\what\Omega, \what\Bsc, \what\nu,\what\theta)$  
indexed as 
$\what \Qcal = \{\what  Q_1, \what Q_2, \what Q_3, \dots, \}$.  
Define the  one-sided full shift $\what \Sigma_{\what \Qcal}:= \{\what Q_1,\what Q_2,\dots \}^\N$.  

Given $\omega\in \what  \Omega,$ write $\bar \omega\in \Sigma_{\what \Qcal}$ for the \emph{$\what \Qcal$-itinerary of $\omega$};  that is, for $k\ge 0$, 
$$\bar \omega_k=\what Q_i \iff \hat \theta^k(\omega)\in \what Q_i.$$
Since $\what\Qcal$ is a one-sided generator for $(\what\Omega, \what\Bsc, \what\nu,\what\theta)$, the map $\omega\mapsto\bar\omega$ is injective $\bmod \hat \nu$.  Let $\bar \nu$ denote the image of $\hat \nu$ under the map $\omega\mapsto\bar\omega$; let $\bar \mu$ be the image of $\hat \mu$ under the map $(\omega,x) \mapsto(\bar\omega,x)$.
Finally, let $f_{\bar \omega} = f_{\omega}$.  If  $\omega\mapsto f_\omega$ is $\sigma(\what\Qcal)$-measurable, then $f_{\bar\omega}$ depends only on the $0$th coordinate $\bar \omega_0$  of $\bar \omega$, corresponding to the atom of $\what \calQ$ containing $\omega$.  
Define the locally constant skew-extension $\bar F\colon \what \Sigma_{\what \calQ}\times M\to  \what \Sigma_{\what \calQ} \times M$ by the map $$\bar F(\bar \omega, x) = (\sigma(\bar\omega), f_{\bar \omega}(x)).$$

\begin{proposition}\label{prop:osgimpliesshift}
Suppose that  $\what\Qcal$ is a one-sided generator for $(\what\Omega, \what\Bsc, \what\nu,\what\theta)$ and that $\omega\mapsto f_\omega$ is $\sigma(\what\Qcal)$-measurable.  Then the map $(\omega,x) \mapsto(\bar\omega,x)$ is a measurable conjugacy  
between measure-preserving dynamical systems $$F\colon (\what \Omega\times M,\hat \mu)\to (\what \Omega\times M,\hat \mu)$$ and 
$$\bar F\colon (\Sigma_\Qcal \times M,\bar \mu)\to (\Sigma_\Qcal \times M,\bar\mu).$$
\end{proposition}
Since the map $(\omega,x) \mapsto(\bar\omega,x)$ is the identity in the second  coordinate, the projections of the measures $\hat \mu$ and $\bar \mu$ to $M$ coincide.  

We note that while the skew product reduces to a locally constant skew-extension over a countable full shift, the base measure $\bar \nu$ need not have the structure of a Bernoulli product measure $\nu_0^\N$.   Moreover, the sample measures associated to $\bar \mu$ need not exhibit any invariance under changing the future as exhibited by the sample measures associated to stationary measures in \eqref{eq:statprop} of \cref{sec:skewext}.  Thus, measures invariant under skew-products satisfying \cref{Qmsble} are in a more general class than the stationary measures for i.i.d.~random walks.

\subsubsection{Conditionally future independent measures}\label{sec:condfIP}
We consider a subclass of  invariant measures for locally constant skew extensions that is slightly more general than the natural extension of i.i.d.~random walks.  This is less general than the measures invariant under locally constant cocycles in \cref{sec:reduction} since we impose the  geometric feature described by \eqref{eq:statprop} of \cref{sec:skewext}.  

Let $\Sigma=A^\Z$ be a full shift on a countable alphabet $A$.  Given the left-shift $\sigma$ on $\Sigma$, let $\nu$ be an ergodic, $\sigma$-invariant probability measure with finite entropy $h_\nu(\sigma)<\infty$.  
Let $\pi_i\colon \Sigma\to A$ denote the $i$th coordinate projection.  Given $I\subset \Z$, let $\sigG_I$ denote the sub-$\sigma$ algebra on $\Sigma$ generated by $\{\pi_i:i\in I\}$.  

We consider a skew extension $F\colon \Sigma\times M\to\Sigma\times M$ of the form $(\omega,x)\mapsto (\sigma(\omega), f_\omega(x))$.   As in the setting of i.i.d.\ random walks, we assume the skew-extension is locally constant; that is the function $\Sigma\to \diff^r(M)$,  $\omega\mapsto f_\omega$ depends only on the 0th coordinate $\omega_0$ of $\omega$.  
Given a locally constant skew extension $F$, we say that an $F$-invariant probability measure $\mu$ on $\Sigma\times M$ projecting to $\nu$ is \emph{conditionally future-independent} if the sample measures $\omega\mapsto\mu_\omega$ are $\sigG_{<0}$-measurable. Equivalently, conditioned on the strict past of the driving process, the present fiber state is independent of the present and future driving symbols. 
\begin{lemma}
If the random dynamical system $F$ is locally constant then for any probability measure $\mu_0$ on $M$,   any limit point of the sequence  $$\mu_N = \frac 1 N \sum_{n=0}^{N-1}F^n_* (\nu\times \mu_0)$$ is an $F$-invariant, conditionally future-independent probability measure.  
\end{lemma}

\subsection{Standing hypotheses}\label{sec:standing ass skew}
For  our main applications below, we will make the following finite entropy assumption: 
\begin{equation*}h _{\hat \nu}(\hat \theta, \ptpart):= H_{\hat \nu}(\ptpart\mid \hat \theta\inv  \ptpart)<\infty \tag{H}
\label{eq:finent}
\end{equation*}
and the integrability assumption:
\begin{equation*} 
  \label{eq:L1} 
\omega\mapsto \log ^+ \|f_\omega \|_{C^{1+\beta}}, \ \  \omega\mapsto \log ^+ \|f_\omega\inv \|_{C^{1+\beta}}  \in L^1(\hat \nu).
\tag{L1}
\end{equation*}

At times, we will also assume the existence of a compatible \emph{backwards random dynamical system}: 
\begin{align*}\label{eq:BD}
  \begin{split}
&\text{\itshape  $\hat \theta\colon (\what\Omega,\hat \nu)\to (\what\Omega,\hat \nu) $ admits a one-sided generator $\what\Qcal$, and} \\
&\text{\itshape the map $\what\Omega\to \Diff^{1+\beta}(M)$, $\omega\mapsto f_{\omega}$, is $\sigma(\what \Qcal)$-measurable.} 
\end{split}
\tag{BD}
\end{align*}

Recall that   hypothesis \eqref{eq:BD} implies $h_{\hat \nu} (\hat \theta)<\infty$ and so hypothesis \eqref{eq:BD} implies hypothesis \eqref{eq:finent}.    Also recall from \cref{prop:osgimpliesshift}  that hypothesis \eqref{eq:BD} implies that the shift $F\colon \Omega\times M\to \Omega \times M$ can be recoded as a locally constant skew extension over a (countable) full shift.

\subsection{Smooth ergodic theory of skew extensions}
We briefly recall several facts in smooth ergodic theory of the skew extension before stating the main theorem.  Further details will appear in \cref{sec:lyapunov}.

First, we assume \eqref{eq:L1} of \cref{sec:standing ass skew}.
Let $ \mu$ be an ergodic, $ F$-invariant probability measure on $ X$ with $( \pi_1) _*  \mu =  \nu$.  

\subsubsection{Lyapunov exponents and Oseledec's splitting}
We  apply Oseledec's theorem to the fiberwise derivative cocycle $(  \omega,x)\mapsto D_xf_{ \omega}^{(n)}$ over the action of $  F$  on $(  X,   \mu)$ and obtain Lyapunov exponents (listed without multiplicity)  
$$\lambda_1<\lambda_2<\dots <\lambda_p.$$
Throughout the paper, we will assume that $\lambda_1<0<\lambda_p$ and $\lambda_i\neq 0$ for all $i$. We find the natural number $s\ge 1$ so that \[\lambda_1<\dots<\lambda_s <0 < \lambda_{s+1}<\dots <\lambda_p.\] 
For $\mu$ almost every $(\omega,x)$, let $$E^s_{(\omega , x)} := \left\{v\in T_xM: \limsup_{n\to \infty}  \frac 1 n \log \|D_x f_{\omega }^{(n)}v\|<0 
\right\}$$
and  $$E^u_{(\omega , x)} := \left\{v\in T_xM: \limsup_{n\to \infty}  \frac 1 n \log \|D_x f_{\omega }^{(-n)}v\|<0 
\right\}$$
be the \emph{stable} and \emph{unstable Oseledec's subspaces} associated with the fiber dynamics.  

We note:
\begin{enumerate}
\item  $E^s_{(\omega , x)}$ is defined using only future iterates of $F$. In particular, the map $(\omega,x)\mapsto E^s_{(\omega,x)}$ is $\Asc\otimes \Bsc_{M}$-measurable. It follows that $E^s_{(\what\omega,x)}$ is defined for $\hat \mu$-a.e.~$(\what\omega,x)\in \what X$.
\item  $E^u_{(\omega,x)}$ is defined using only past iterates of $F$. In particular, if we assume \eqref{eq:BD},  the map $(\omega,x)\mapsto E^u_{(\omega,x)}$ is $\Asc'\otimes \Bsc_{M}$-measurable.  In this case, we can define the \emph{backward stable Oseledec's subspaces} for the backwards random dynamics: for $\wcheck\mu$-a.e.~$(\wcheck\omega,x)\in \wcheck X$, \[E^{s}_{(\wcheck\omega,x)}\coloneqq \left\{v\in T_{x}M : \limsup_{n\to \infty}\frac{1}{n} \log \|D_{x}g_{\wcheck\omega}^{(n)} v\|<0\right\}.\] In this case, we observe that, as $\wcheck F$ on $\wcheck X$ is the factor of $F^{-1}$ on $X$, the map $(\omega,x)\mapsto E^{u}_{(\omega,x)}$ factors through $(\wcheck\omega,x)\mapsto E^{s}_{(\wcheck \omega,x)}$.  That is, for $\mu$-a.e.~$(\omega, x) \in X$, $$E_{(\omega,x)}^{u} = E^{s}_{(\wcheck\tau(\omega),x)}.$$
\end{enumerate}
Since we assume throughout the paper $\lambda_1<0<\lambda_p$, we always have $\dim(E^s_{(\omega , x)})\ge 1$ and  $\dim(E^u_{(\omega , x)})\ge 1$ for $\mu$-a.e.~$(\omega,x)\in X$.  

\subsubsection{Random stable and unstable manifolds}\label{sss:randomstable}
Given $(\omega,x)\in  X$, write 
$$W^s_{\omega}(x):=\left\{y\in   M:  \limsup_{n\to \infty}  \frac 1 n \log d_M( f_{ \omega}^{(n)}(x),f_{ \omega}^{(n)}(y)) <0\right\}$$
and 
$$W^u_{\omega}(x):=\left\{y\in   M:  \limsup_{n\to \infty}  \frac 1 n \log d_M( f_{ \omega}^{(-n)}(x),f_{ \omega}^{(-n)}(y)) <0\right\}$$
for the \emph{stable} and \emph{unstable sets} in $M$ through $x$ for the fiberwise dynamics determined by $\omega$. 
We observe that the stable set is defined using only forwards dynamics $f_\omega^{(n)}$, $n\ge 0$.  In particular, since  $\omega\mapsto f_\omega^{(n)}$ is $\sigA$-measurable for every  $n\ge 0$,  $W^s_{\what\omega}(x)$ is defined for $(\what\omega,x)\in \what X$ and given $(\omega,x)\in X$, we have \[ W^s_{\omega}(x)=W^s_{\tau(\omega)}(x) \quad\textrm{for}\quad (\omega,x)\in X.\]

Similarly, under assumption \eqref{eq:BD}, we let 
	\[W^{s}_{\wcheck\omega}(x)=\left\{y \in M : \limsup_{n\to \infty} \frac 1 n \log d_{M}(g_{\wcheck\omega}^{(n)} x, g_{\wcheck\omega}^{(n)}y ) < 0 \right\}\] be the \emph{backward stable set} determined by $\wcheck \omega \in \wcheck \Omega$ through $x$ in $M$.  Since the unstable set is defined using only backward dynamics $f_\omega^{(n)}$, $n< 0$ and since, under assumption \eqref{eq:BD}, $\omega\mapsto  g_\omega^{(n)} = f_\omega^{(-n)}$ is $\sigA'$-measurable, 
	 we have 
	 \[ W^{u}_{\omega}(x)=W^{s}_{\wcheck\tau(\omega)}(x) \quad \textrm{for} \quad(\omega,x)\in  X.\]

Applying a fibered version of Pesin's stable manifold theorem (c.f.\ \cite[Thm.\  III.3.2]{MR1369243} for the case of $C^{1+\mathrm{Lip}}
$ fibered systems or discussion and references in \cite{MR4599404} for the setting of $C^{1+\text{H\"older}}$-fibered systems) 
we obtain the following.  
\begin{lemma} \label{lem:pesin}
Assume $f_{\omega}\in \Diff^r (M)$ with $r>1$ for $\nu$-a.e.~$\omega\in \Omega$ and satisfies the integrability assumption \eqref{eq:L1} of \cref{sec:standing ass skew}.  Then, for $\mu$-a.e.\ $(\omega,x)\in X$ and $\sigma\in \{s,u\}$,  the set $W^\sigma_{\omega}(x)$ is an injectively immersed, $C^{1+\beta}$ (in fact $C^r$ if the $C^r$ norms in   \eqref{eq:L1} are $\log$-integrable) submanifold tangent to $E^\sigma_{(\omega,x)}$ at $x$.
\end{lemma}

For  $\sigma\in \{s,u\}$, we also write $$W^\sigma{(\omega,x)}:= \{\omega\}\times W^\sigma_{\omega}(x)$$ 
for the corresponding  \emph{fiber stable set} or  \emph{fiber unstable set} in $X$ through $(\omega,x)$.  We view $W^\sigma{(\omega,x)} $  as an injectively immersed submanifold of $\{\omega\}\times M$.  The sets $W^\sigma {(\omega,x)}$ form atoms of a (typically non-measurable) partition of $X$. We similarly define the \emph{fiber stable set} $W^{s} (\what \omega,x)$ in $\what X$ and, assuming  \eqref{eq:BD}, the \emph{backwards fiber stable set}  $W^{s}(\wcheck\omega,x)$ in $\wcheck X$. 

\subsection{Subordinate and Markov partitions; notation for induced partitions}\label{sec:submarkov}
We note that the partitions into full fiber stable or full fiber unstable manifolds are typically non-measurable.  
To calculate dimension and entropy quantities, we define measurable partitions of $\what X$, $ X$, and---assuming  \eqref{eq:BD}---$\wcheck X$ with certain geometric and dynamical properties.

\subsubsection{Definitions and notation}\label{sec:notanddefnparti}
\begin{definition}[Subordinate partitions]\label{def:subord} 
\
\begin{enumerate}
\item
 A measurable partition $ \what \xi^s$  of $(\what X,\hat \mu)$ is \emph{subordinate} to the partition into fiber stable manifolds if for $\what \mu$-a.e.\ $(\what\omega, x)\in \what X$ 
the following hold:
\begin{enumerate}
\item $\what\xi(\what\omega,x)\subset W^s{(\what\omega,x)}$;  
\item each atom $\what \xi(\what\omega,x)$ contains an open (with respect to the immersed topology in $M$) neighborhood of $x$ in $W^s{(\what\omega,x)}$;
\item each atom $\what\xi(\what\omega,x)$ is precompact (with respect to the immersed topology in $M$)  in $W^s{(\what\omega,x)}$.
\end{enumerate}
\item For $\sigma\in \{s,u\}$, we analogously define a measurable partition $ \overline  \xi^{\sigma}$  of the two-sided skew extension $(  X,  \mu)$ to be  \emph{subordinate} to the partition into fiber stable manifolds or fiber unstable manifolds.

\item If we further assume \eqref{eq:BD}, then we also analogously define a measurable partition $\wcheck\xi^{s}$ of $(\wcheck X, \wcheck \mu)$ to be \emph{subordinate} to the partition into  backwards fiber stable manifolds.

\end{enumerate}

\end{definition}

\begin{definition}[Sub-invariant partitions] \
\

\begin{enumerate}
\item We say a measurable partition $\zeta$ of $(\what X,\hat \mu)$  is \emph{sub-invariant} or \emph{(forwards) Markov} if $\what F^{-1}\zeta\prec\zeta$. We define similarly, \emph{sub-invariant} or \emph{(forwards) Markov} measurable partitions under $F$ of $(X,\mu)$.  Assuming   \eqref{eq:BD}, we similarly define   \emph{sub-invariant} or \emph{(forwards) Markov} measurable partitions under $\wcheck F$ of $(\wcheck  X,\check \mu)$.

\end{enumerate}

\end{definition}

\begin{definition}[Notation for the pull-back of a subordinate partition]\label{def:pullback}\

\begin{enumerate}
\item For a measurable partition $\what\xi^s$ 
of $(\what X, \hat \mu)$ subordinate  to fiber stable manifolds, 
we define  \[\xi^{s}\coloneqq\tau^{-1} \what\xi^{s}.\]
Note that atoms of $\xi^s$ are saturated by atoms of the $\sigma$-algebra $\sigA$.  
By further conditioning on the fiber, we have that  \[\overline{\xi}^{s} =\xi^{s} \vee \Xi (\Fsc)\] is a measurable partition of $(X,\mu)$ subordinate  to  fiber stable manifolds.  

\item Similarly, assuming \eqref{eq:BD}, 
given a measurable partition $\wcheck \xi^s$ 
of $(\wcheck  X, \wcheck  \mu)$ subordinate  to backwards fiber stable manifolds,   
we define  \[\xi^{u}\coloneqq\check \tau^{-1} \wcheck \xi^{s}.\]
As above, $\xi^u$ is saturated by $\sigA'$; however, by further conditioning on the fibers we have that  \[\overline{\xi}^{u} =\xi^{u} \vee \Xi (\Fsc)\] is a measurable partition of $(X,\mu)$ subordinate  to  fiber unstable manifolds.

\end{enumerate}
\end{definition}

\subsubsection{Existence of sub-invariant subordinate partitions} 

Following the construction in \cite{MR693976} (see also discussion in \cite{MR819557}) sub-invariant measurable partitions subordinate to our foliations exist under our standing hypotheses.  
See also discussion in \cite[(4.3)]{MR968818} or \cite[Prop.~IV.2.1]{MR1369243} for constructions in the setting of i.i.d.~random dynamics.  An outline of their construction also appears in the setting of general random dynamics in \cite[\S5.3]{BRH}.

 \begin{lemma}\label{lem:LS}
Assume hypotheses \eqref{eq:L1} of \cref{sec:standing ass skew}.  
Then there exists a sub-invariant (for $\what F$) measurable partition $\hat  \xi^s$ of $(\what X, \hat \mu)$ subordinate  to the partition into fiber stable manifolds.   
Similarly, there exist measurable partitions $\overline  \xi^s$ (resp.~$\overline  \xi^u$) of $(  X, \mu)$ subordinate  to the partition into fiber stable  (resp.~unstable manifolds) that are sub-invariant for $  F$ (resp.~$F\inv$). 

Assume hypotheses \eqref{eq:L1} and \eqref{eq:BD} of \cref{sec:standing ass skew}.  Then there exists a measurable partition $\check  \xi^s$ of $(\wcheck X, \wcheck \mu)$ subordinate  to the partition into backwards fiber stable manifolds that is sub-invariant for $\check F$.  
\end{lemma}

\subsubsection{Generating properties of sub-invariant subordinate partitions}\label{subsec:subord}
Let $\hat \xi^s$, $\overline  \xi^s$, $\overline  \xi^u$ and (assuming \eqref{eq:BD} of \cref{sec:standing ass skew}) $\check  \xi^s$ be as in \cref{lem:LS}.  Let $\xi^s= \tau\inv \hat \xi^s$ and  (assuming \eqref{eq:BD} of \cref{sec:standing ass skew}) $\xi^u = \check \tau\inv \check  \xi^s$ be as in the notation of \cref{def:pullback}.  

Then 
\begin{enumerate}
\item  \label{stable3} \label{partition3}$\bigvee_{n=0}^\infty F^{n}(  \xi^s)=\bigvee_{n=0}^\infty F^{n}( \tau\inv \hat \xi^s)$ is the point partition on $ X$;
	\item $\Xi(\sigC) \prec \xi^s$ and  \label{stable5}$\bigvee_{n=0}^\infty F^{n}(  \xi^s)$ is the point partition on $X$;
	\item $\bigvee_{n=0}^\infty F^{n}( \overline \xi^s)$ is the point partition on $X$;
	\item $\bigvee_{n={-\infty}}^0 F^{n}( \overline \xi^u)$ is the point partition on $X$;
	\item $\Xi(\sigC') \prec \xi^u$ and  \label{unstable5}$\bigvee_{n=-\infty}^0 F^{n}(  \xi^u)$ is the point partition on $X$;
	\item $\bigvee_{n={-\infty}}^0  F^{n}( \xi^u)=\bigvee_{n={-\infty}}^0  F^{n}( \check \tau\inv \check \xi^s)$ is the point partition on $ X$;
\end{enumerate}

\subsection{Stable and unstable pointwise dimensions}
Given $x\in M$ and $r>0$, write $B(x,r)$ for the Riemannian ball of radius $r$ centered at $x$ in $M$. 

Given an ergodic,  $F$-invariant probability measure  $\mu$ on $X$  with $(\pi_1)_*\mu=\nu$, we define various stable and unstable dimensions of $\mu$ relative to various partitions of the space $X$ with dynamical properties.

\subsubsection{Stable and unstable dimensions}
\begin{proposition}\label{prpo:sudimmfull} Assume 
the integrability  \eqref{eq:L1} 
hypotheses 
of \cref{sec:standing ass skew}. Let $\mu$ be an ergodic,  $F$-invariant probability measure on $X$  with $(\pi_1)_*\mu=\nu$.  Let $\hat \xi^s, \bar\xi^s$, and $\bar \xi^u$ be subordinate partitions as in \cref{def:subord}.   Then the following hold:
\begin{enumerate}
\item Assume the entropy  \eqref{eq:finent} hypotheses of \cref{sec:standing ass skew}.  For $\hat \mu$-a.e.~$(\what\omega,x)\in \what X$ the following limit exists and is independent of $(\what\omega,x)$:
\begin{equation*} \lim_{r\to0} \frac{\log \what \mu_{(\what\omega,x)}^{\hat \xi^{s}}\left(\hat \pi_{2}^{-1}\left({B}(x,r)\right)\right)}{\log r}\eqqcolon  \hat d^{s}.
\end{equation*}
\item \label{propitme:sdimfull} For $ \mu$-a.e.~$(\omega,x)\in   X$ the following limit exists and is independent of $(\omega,x)$: 
\begin{equation*} \lim_{r\to0}\frac{\log\mu_{(\omega,x)}^{\bar{\xi}^{s}}\left( \pi_{2}^{-1}\left(B(x,r)\right)\right)}{\log r}\eqqcolon \overline{d}^{s}.
\end{equation*}

\item   \label{propitme:udimfull} For $ \mu$-a.e.~$(\omega,x)\in   X$ the following limit exists and is independent of $(\omega,x)$: 
\begin{equation*} \lim_{r\to0}\frac{\log\mu_{(\omega,x)}^{\bar{\xi}^{u}}\left(\pi_{2}^{-1}\left(B(x,r)\right)\right)}{\log r}\eqqcolon \overline{d}^{u}.
\end{equation*}

\end{enumerate}
If we further assume the hypotheses \eqref{eq:BD}  
of \cref{sec:standing ass skew},  then with $ \wcheck\xi^s$  subordinate partition as in  \cref{def:subord},
\begin{enumerate}[resume]
\item \label{propitme:udimfullBD} for $\wcheck\mu$-a.e.~$(\wcheck\omega,x)\in \wcheck X$ the following limit exists and is independent of $(\wcheck\omega,x)$:
\begin{equation*} \lim_{r\to0}\frac{\log\wcheck\mu_{(\wcheck\omega,x)}^{\wcheck\xi ^ s } \left(\wcheck \pi_{2}^{-1}\left(B(x,r)\right)\right)}{\log r}\eqqcolon \wcheck{d}^{u}
\end{equation*}
\end{enumerate}
The values $ \what d^{s}$,$  \overline{d}^{s}$,  $\overline{d}^{u}$, and $\wcheck d^{u}$ are independent of the choice of subordinate partitions.
\end{proposition}

We note that for  $ \mu$-a.e.~$(\omega,x)\in  X$, we may alternatively write  \begin{equation*} \lim_{r\to0} \frac{\log \mu_{(\omega,x)}^{\xi^{s}}\left(\pi_{2}^{-1} \left(  {B}(x,r) \right) \right)}{\log r} = \what d ^s.
\end{equation*}
If we further assume (\ref{eq:BD}), then  for $\mu$-a.e.~$(\omega,x)\in X$ we alternatively have 
\begin{equation*} \lim_{r\to0} \frac{\log \mu_{(\omega,x)}^{\xi^{u}} \left(\pi_{2}^{-1} \left( B(x,r)\right)\right)}{\log r}=\wcheck d^u \end{equation*}

\begin{remark}[Unstable dimension of conditionally future independent invariant measures]
Note that in the setting of locally constant skew extensions over a full shift, for conditionally future independent invariant measures $\mu$ as in \cref{sec:condfIP}, we have (identifying all fibers with $M$ or applying $(\pi_2)_*)$ $\mu^{\xi^u}_{(\omega,x)}=\mu_{(\omega,x)}^{\bar\xi^{u}}$ for $\mu$-a.e.\ $(\omega,x)\in X$, hence, $\bar d^{u}=\wcheck{d}^{u}$
\end{remark}

\begin{remark}[Entropy assumptions]
\label{rem:FWRIeq}
In 
 \cref{prpo:sudimmfull}\eqref{propitme:sdimfull}--\eqref{propitme:udimfullBD}, no entropy assumption on the base system is needed.  This is because the  relevant entropies used in the proofs are already finite.
 Indeed,   the integrability  \eqref{eq:L1}   hypotheses 
of \cref{sec:standing ass skew} 
 combined with the fiberwise Ruelle inequality \cite{MR1314494} implies the entropies leading to the dimensions  in \eqref{propitme:sdimfull} and \eqref{propitme:udimfull} (the conditional metric entropy conditioned on the fiber $\sigma$-algebra $\sigF$) are finite.  
For  \cref{prpo:sudimmfull}\eqref{propitme:udimfullBD},  hypothesis \eqref{eq:BD} of \cref{sec:standing ass skew} directly implies finiteness of the base entropy $h_{\nu}(\theta)<\infty$;  the fiberwise Ruelle inequality combined with the Abramov--Rohlin formula (or \cref{lem:natext} below) implies finiteness $h_\mu(F)<\infty$.  
\end{remark}

\subsubsection{Stable and unstable dimensions relative to a sub-invariant partition}
We introduce slightly more general notations of pointwise dimensions by conditioning further on a sub-invariant partition.  This extra abstraction adds only additional notation and bookkeeping in the proofs; it has the advantage that it allows us to prove results simultaneously and might be of use in future work.

Throughout, we will consider a measurable partition $\zeta$ of $(X,\mu)$ with the following properties:
\begin{enumerate}
\item $\zeta $ is sub-invariant: $F^{-1}\zeta\prec \zeta$.
\item $h_{\mu}(F,\zeta\vee \Xi(\Csc)) <\infty$. 
\end{enumerate}

\begin{proposition}[{c.f.~\cite[Prop.~4.8]{BRH}}]\label{prpo:sudimm} Assume 
the integrability  \eqref{eq:L1}  hypotheses of \cref{sec:standing ass skew}. Let $\mu$ be an ergodic,  $F$-invariant probability measure on $X$  with $(\pi_1)_*\mu=\nu$.   Let $\zeta$ be a sub-invariant measurable partition of $(X,\mu).$ 

\begin{enumerate} 
\item Suppose $h_{\mu}(F,\zeta\vee \Xi(\Csc)) <\infty$.  Then, for  $\mu$-a.e.~$(\omega,x)\in X$ the following limit exists and is independent of $(\omega,x)$: 
\begin{equation*} \lim_{r\to0} \frac{\log \mu_{(\omega,x)}^{\xi^{s}\vee\zeta}\left(\pi_2^{-1}\left(  {B}(x,r) \right) \right)  }{\log r}\eqqcolon \dim (\mu\mid \xi^{s}\vee\zeta)\eqqcolon \hat d_{\zeta}^{s}.
\end{equation*}

\item \label{eq:dimmF}For $\mu$-almost every $(\omega,x)\in X$ the following limits exist and are independent of $(\omega,x)$:
\begin{equation*} \lim_{r\to0}\frac{\log\mu_{(\omega,x)}^{\overline{\xi}^{u}\vee\zeta^{+}}\left( \pi_{2}^{-1} {B}(x,r) \right)}{\log r}\eqqcolon \dim(\mu\mid\overline{\xi}^u \vee \zeta^+) \eqqcolon\overline{d}^{u}_{\zeta^+}
\end{equation*}
and 
\begin{equation*} \lim_{r\to0}\frac{\log\mu_{(\omega,x)}^{\overline{\xi}^{s}\vee\zeta}\left( \pi_{2}^{-1} {B}(x,r) \right)}{\log r}\eqqcolon \dim(\mu\mid\overline{\xi}^s \vee \zeta) \eqqcolon\overline{d}^{s}_{\zeta}.
\end{equation*}
\item \label{eq:dimmBD} Suppose hypothesis (\ref{eq:BD}) holds.
 Then for $\mu$-almost every $(\omega,x)\in X$, the following limit exists and is independent of $(\omega,x)$:
\begin{equation*} \lim_{r\to0}\frac{\log\mu_{(\omega,x)}^{{\xi}^{u}\vee\zeta^{+}}\ \left(\pi_{2}^{-1} {B}(x,r) \right)}{\log r}\eqqcolon \dim(\mu\mid\xi^u \vee \zeta^+) \eqqcolon\wcheck{d}^{u}_{\zeta^+}.
\end{equation*}
\end{enumerate}
\end{proposition}

Again, no entropy assumption is needed in 
\cref{prpo:sudimm}\eqref{eq:dimmF} and \eqref{eq:dimmBD} because  the fiberwise Ruelle inequality implies the relevant entropy is finite.  See \cref{rem:FWRIeq}.  

\begin{remark}
The abstraction of allowing for various sub-invariant  $\zeta$ in \cref{prpo:sudimm} and our main result \cref{thm:prod} below may be of use in  future results.  
Even for our main applications in \cref{mainskewcor}, allowing for different $\zeta$ allows us to obtain results simultaneously.  For our main applications, we take $\zeta$ to be either:
\begin{enumerate}
	\item $\zeta = \Xi(\sigF)$, the partition of $X$ into fibers.
	\item $\zeta = \Xi(\sigC)$, the preimage of the partition of $\what X$ into fibers under $\tau\colon X\to \what X$.
	\item $\zeta= \{X,\emptyset\}$, the trivial partition of $X$.  
\end{enumerate}
\end{remark}

\subsection{Main result: exact dimensionality and dimension formulas} 
Having established the above notation, we are now in a position to formulate the main abstract result of the paper.  All applications discussed above are special cases.  
We again note as in \cref{sec:iid} that the results are not new when the fiberwise Lyapunov exponents are all of the same sign.  Thus, in our proofs of the main results, we will always assume the existence of a  fiberwise Lyapunov exponent of each sign. 

\subsubsection{Main abstract theorem}
With the notation in \cref{prpo:sudimm}, we have the following.  
\begin{theorem}\label{thm:prod} Assume the integrability (\ref{eq:L1}) hypotheses of \cref{sec:standing ass skew}. Let $\mu$ be an ergodic, $F$-invariant probability measure with $(\pi_1)_*\mu=\nu$. 
Assume no fiberwise Lyapunov exponent $\lambda_i$ of $F$ with respect to $\mu$ is zero.  
\begin{enumerate}[label*={(\arabic*)}, ref= {(\arabic*)} ]
\item\label{thm:prod:forward} 
Let $\zeta$ be a  sub-invariant measurable partition on $X$ with
$\Xi(\sigC)\prec \zeta$ and 
 $h_{\mu}(F,\zeta) <\infty$.  
Then 
for $\mu$-a.e.\ $(\omega,x)$,
\begin{equation*}\lim_{r\to 0}\frac{\log\mu_{(\omega,x)}^{\zeta}\left(\pi_2^{-1}\left({B}(x,r)\right)\right)}{\log r }= \overline{d}^{u}_{\zeta^+} + \hat d_\zeta ^s.\end{equation*}

\item\label{thm:prod:backward} Further assume  hypotheses \eqref{eq:BD} of \cref{sec:standing ass skew} and let $\zeta$ be  a  sub-invariant measurable partition on $X$ with
 $h_{\mu}(F,\zeta \vee \Xi(\sigC) ) <\infty$ and $h_{\mu}(F\inv ,\zeta^+\vee \Xi(\Csc')) <\infty$.  
Then 
for $\mu$-a.e.\ $(\omega,x)$,
\begin{equation*}\lim_{r\to 0}\frac{\log\mu_{(\omega,x)}^{\zeta}\left(\pi_2^{-1}\left({B}(x,r)\right)\right)}{\log r }= \widecheck{d}^{u}_{\zeta^+} + \hat d_\zeta ^s.\end{equation*}
\end{enumerate}

\end{theorem}

As discussed in \cref{sss:APS} one can derive an almost product structure from \cref{thm:prod}.

After establishing preliminaries in \cref{sec:4} and \cref{sec:5}, the proof of \cref{thm:prod} will follow directly from 
\cref{prop:UB,prop:LB}, following the notation in \cref{sec:6notation}.

\subsubsection{Main result for general random dynamics} 
Specializing \cref{thm:prod} to the  cases $\zeta=\Xi(\Csc),$ $\zeta=\Xi(\Fsc)$, or the trivial partition when also assuming  \eqref{eq:BD}, we obtain the following: 
\begin{theorem}\label{mainskewcor}
Assume  the integrability \eqref{eq:L1} and finite entropy hypotheses \eqref{eq:finent} of \cref{sec:standing ass skew}. Let $\mu$ be an ergodic, $F$-invariant probability measure with $(\pi_1)_*\mu= \nu$.  Assume no fiberwise Lyapunov exponent $\lambda_i$ of $F$ with respect to $\mu$ is zero.  Then the following hold:  
\begin{enumerate}
\item\label{Fscpt} For $\mu$-a.e.\ $(\omega,x)\in X$, we have 
\[\lim_{r\to0} \frac{\log\mu_{\omega}\left(\pi_{2}^{-1}\left(B(x,r)\right)\right)}{\log r} =\overline{d}^{u}+ \overline{d}^{s}.\]

\item\label{Cscpt} For $\hat \mu$-a.e.\ $(\omega,x)\in \what X$, we have 
\[\lim_{r\to0} \frac{\log \hat \mu_{\omega}\left(\what \pi_{2}^{-1}\left(B(x,r)\right)\right)}{\log r} =\overline{d}^{u}+\what{d}^{s}.\]

\end{enumerate}
If we further assume hypotheses  \eqref{eq:BD} of \cref{sec:standing ass skew},  then 
\begin{enumerate}[resume]
\item\label{trivialpt} for $\mu$-a.e.\ $(\omega,x)\in X$, we have 
\[\lim_{r\to 0}\frac{\log\mu\left(\pi_{2}^{-1}\left(B(x,r)\right)\right)}{\log r }= \wcheck d ^u+  \what{d}^s.\]
\end{enumerate}
\end{theorem}

Indeed, for $\mu$-a.e.\ $(\omega,x)\in X$ we have that $\mu_\omega= \mu_{(\omega,x)}^\sigF$ and so \cref{mainskewcor}\eqref{Fscpt}
follows from \cref{thm:prod} with $\zeta = \Xi(\sigF)$.  
Similarly, for $ \mu$-a.e.\ $(\omega,x)\in  X$
$$\hat \mu_{\tau(\omega)}\left(\what \pi_{2}^{-1}\left(B(x,r)\right)\right)= 
\mu_{(\omega,x)}^\sigC\left( \pi_{2}^{-1}\left(B(x,r)\right)\right),$$
and so \cref{mainskewcor}\eqref{Cscpt} follows from \cref{thm:prod} with $\zeta = \Xi(\sigC)$.  
\cref{mainskewcor}\eqref{trivialpt} follows from \cref{thm:prod} with $\zeta$ the trivial partition.

\begin{remark}[Relation to  \cite{MR2218998}]
We remark that when the maps $f_\omega$ are assumed to be $C^2$ diffeomorphisms (and assuming $\log \|f_\omega\|_{C^2}$ and $ \log\|f_\omega\inv \|_{C^2}$ are $L^1(\nu)$), the dimension formula for the fiberwise conditional measures in \cref{mainskewcor}\eqref{Fscpt} already follows from the main result of \cite{MR2218998}.  Our proof shows the results of \cite{MR2218998} still hold when the generators are $C^{1+\beta}$.    Indeed, if our original base dynamics  $\hat \theta\colon  (\what \Omega, \hat \nu)\to   (\what \Omega, \hat \nu)$ had been invertible, then $F\colon X\to X$ and $\what F\colon \what X\to \what X$ coincide.  We recover the results of  \cite{MR2218998} (for random dynamics with $C^{1+\beta}$ generators) by considering $\zeta= \Xi(\sigF)$ the invariant partition into fibers of $X$.  
We note, however, that the conclusions of \cref{mainskewcor}\eqref{Cscpt} and  \cref{mainskewcor}\eqref{trivialpt} compute the dimensions of certain projections of measures on $X$ whereas the  results in \cite{MR2218998} only compute the sample dimension of the sample measures.  
\end{remark}

\subsection{An invariance principle}
Recall that $\what F$-invariance of $\hat \mu$ gives the following ``invariance on average'' property for the fiber measures $\hat \mu_\omega$:
$$\int (f_{\omega'})_* \hat \mu_{\omega'} \, d \hat \nu^{\hat \theta\inv \what \sigB}_\omega(\omega')= \hat \mu_{\hat \theta (\omega)}.$$
We might instead consider the following stronger invariance property  of the family of measure $\{\hat \mu_\omega\}$ for the one-sided skew system: \begin{equation}\label{eq:IP?}(f_\omega)_*\hat \mu_\omega = \hat \mu_{\hat \theta (\omega)} \text{  \ \ for $\hat \nu$-a.e.~$\omega\in \what \Omega.$} \end{equation}

The Furstenberg entropy $\furstent_{\hat \mu}(\what F\mid \what \sigF)$, defined in \cref{sss:furst} below,  is a (non-negative) quantity that roughly quantifies the defect of \cref{eq:IP?}; in particular, vanishing of the Furstenberg entropy, $\furstent_{\hat \mu}(\what F\mid \what \sigF)=0,$ is equivalent to \cref{eq:IP?}.  
A number of related results in the literature, referred to as {\it invariance principles}, imply vanishing of the Furstenberg entropy.  

As shown in \cite{BRH},  in our settings the Furstenberg entropy $\furstent_{\hat \mu}(\what F\mid \what \sigF)$ can  be computed   entirely in terms of the difference between certain contributions to  entropy along the filtration into stable manifolds.  We summarize this in the following:

\begin{theorem}[{Invariance principle, \cite[Thm.~4.17]{BRH}}]\label{thm:INVhyp}
Let $\hat \nu$ satisfy  \eqref{eq:finent} and \eqref{eq:L1} of \cref{sec:standing ass skew}.
Let $\hat \mu$ be an ergodic, $\what F$-invariant probability measure on $\what X$ with $(\hat \pi_1)_*\hat \mu= \hat \nu$ and  let $\mu$ be the unique ergodic $F$-invariant probability measure  on $X$ guaranteed by \cref{natext}.  
Then the following are equivalent:
\begin{enumerate}
\item $\furstent_{\hat \mu}(\what F\mid \what\sigF)= 0$;
\item $\overline{d}^{s}=\hat {d}^{s}$;
\item for $\nu$-a.e.\ $\omega$,  $\mu^\sigF_\omega = \hat \mu^{\what \sigF}_{\tau( \omega)}$ (viewed as probability measures on $M$);
\item for $\hat \nu$-a.e.~$\omega\in \what \Omega$, $(f_\omega)_* \hat \mu_\omega = \hat \mu_{\hat \theta(\omega)}$.
\end{enumerate}
\end{theorem}

When the measure $\mu$ is ergodic and hyperbolic, let
$\dim (\mu)$ denote the a.s. constant value of $$\lim_{r\to0} \frac{\log\mu_{\omega}\left(\pi_{2}^{-1}\left(B(x,r)\right)\right)}{\log r} $$
and let $\dim (\hat \mu)$ denote the a.s. constant value of 
\[\lim_{r\to0} \frac{\log \hat \mu_{\omega}\left(\what \pi_{2}^{-1}\left(B(x,r)\right)\right)}{\log r}\]
Using that $\dim(\mu) =\overline{d}^{u}+\overline{d}^{s}$ and 
$\dim(\hat \mu) =\overline{d}^{u}+\what{d}^{s}$, we obtain the following invariance principle in terms of pointwise dimension for ergodic hyperbolic  measures invariant under random dynamics.  

\begin{corollary}[Invariance principle for hyperbolic measures]
Let $\hat \nu$ satisfy 
\eqref{eq:finent} and \eqref{eq:L1} of \cref{sec:standing ass skew}.
Let $\hat \mu$ be an ergodic, $\what F$-invariant probability measure on $\what X$ with $(\hat \pi_1)_*\hat \mu= \hat \nu$ and  let $\mu$ be the unique ergodic $F$-invariant probability measure   on $X$ guaranteed by \cref{natext}.  Assume that no Lyapunov exponent of $\hat \mu$ is zero.

Then $\dim^\sigF(\mu)\le \dim ^{\what \sigF}(\hat \mu)$; moreover  equality holds if and only if the equivalent criteria in \cref{thm:INVhyp} hold.  
\end{corollary}

\subsection{Application: fibered systems admitting stable holonomies}
Our abstract setup is motivated by a standard reduction for fibered systems admitting (equivariant) systems of stable holonomies.  We briefly describe the setup.  

Suppose that $g\colon Z\to Z$ is a  $C^{1+\beta}$ diffeomorphism  of a compact manifold (or homeomorphism of a compact metric space).  Let $\pi\colon Y \to Z$ be a continuous fiber bundle over $Z$ with compact fibers; the fibers of $\pi$ are assumed to be $C^{1+\beta}$ manifolds but the fibration structure is only assumed continuous.  Let $F\colon Y\to Y$ be a homeomorphism  extending   $g\colon Z\to Z$; that is, $\pi\circ F= g\circ \pi$.  
Given $z\in Z$, let $Y_z= \pi\inv (z)$ denote the fiber over $z$ and let $f_z\colon Y_z \to Y_{g(z)}$ be the transition map; we assume that $f_z$  is a $C^{1+\beta}$ diffeomorphism for every $z$ and that the map $z\mapsto f_z$ is continuous.  After a Borel trivialization with bounded $C^{1+\beta}$-norm, we may assume $Y= Z\times M$.  Thus, we view each $f_z\in \Diff^{1+\beta}(M)$ with $z\mapsto f_z$ bounded and measurable.  

Let $\nu$ be an ergodic, $g$-invariant probability measure on $Z$ (with $h_\nu(g)<\infty$ if $g$ is only assumed continuous). Let $\eta$ be a measurable partition of $(Z,\nu)$ that is sub-invariant and generating; that is $g\inv\eta\prec \eta$ and $g^n\eta\to \ptpart$.  In practice, $\eta$ will be a measurable partition subordinate to a system of (fast)  stable manifolds for $g$ with respect to $\nu$; see \cref{sec:submarkov} for details.  When $g\colon Z\to Z$ is Anosov (or topologically hyperbolic), $\eta$ will be a measurable partition subordinate to stable manifolds.  See discussions in \cite[\S~2.3]{MR2651382}

We assume the existence of holonomy maps: for $\nu$-a.e.~$z$ and for $\nu^{\eta}_z$-a.e.~$z',z''$ in the atom $ \eta(z)$ of $\eta$, there is a $C^{1+\beta}$ diffeomorphism $h_{z,z'}\colon Y_z\to Y_{z'}$ satisfying $$f_{z'}\circ h_{z,z'}= h_{g(z),g(z')}\circ f_z$$
and 
$$h_{z,z'}=h_{z'',z'}\circ h_{z,z''}.$$
We assume that maps $z,z'\mapsto h_{z,z'}$ are sufficiently measurable and bounded so that the passage to the quotient skew-shift preserves integrability.  
The existence of such holonomy maps  is well studied, especially when the base dynamics $g\colon Z\to Z$ is uniformly hyperbolic and when the fiber dynamics $f_z$ is sufficiently fiber-bunched.  See discussions in \cite[\S~2.4]{MR2651382}.

Let $\mu$ be an ergodic $F$-invariant probability measure on $Y$ with $\pi_*\mu= \nu$.  
Let $\{\mu_z\}$ denote a family of conditional measures of $\mu$ relative to the projection by $\pi$.  
Given $\nu$-a.e.~$z\in Z$, let $$\hat \mu_z:= \int \left(h^s_{z',z}\right)_*\mu_{z'} \, d \nu_z^\eta(z').$$

\begin{theorem}
Let $\mu$ be an ergodic $F$-invariant probability measure on $Y$ with $\pi_*\mu= \nu$.  
Suppose every fiberwise Lyapunov exponent of the bundle map $F\colon Y\to Y$ w.r.t. $\mu$ is non-zero.  
Then there are $0\le \overline d, \hat d\le \dim (M)$ such that 
\begin{enumerate}
\item for $\nu$-a.e.~$z\in Z$, the measure $\mu_z$ is exact dimensional with $\dim (\mu_z) = \overline d$;
\item for $\nu$-a.e.~$z\in Z$, the measure $\hat \mu_z$ is exact dimensional with $\dim (\hat \mu_z) = \hat d$;
\item $\overline d\le  \hat d$ with equality if and only if $\mu_z= \hat \mu_z$ for $\nu$-a.e.~$z\in Z$.  
\end{enumerate}

\end{theorem}

\subsection{Examples}
We present three examples showing how the hypotheses and certain dimensions in \cref{mainskewcor} are required.
These  are adapted from discussions with ChatGPT with writing edited by the authors.  

\begin{example}[A (locally) constant cocycle]\label{ex:anosov-graph}
Let $f$ be a smooth Anosov diffeomorphism on a 2-torus $M= \T^2$.  
Fix a Markov partition $\Sigma$ and let $p\colon \Omega \subset \Sigma^{\Zbb} \to M$ be a symbolic coding induced by this Markov partition (so that $\Omega$ is a sub-shift of finite type).   
Let $\sigma\colon \Omega\to\Omega$ be the shift map. 
Fix an ergodic shift-invariant probability measure $\beta$ on $\Omega$. 
Consider the skew-product dynamics $F$ on  $X\coloneqq \Omega\times M$,  $F(\omega,x)=(\sigma(\omega),f (x))$.    Note that $f$ has no $\omega\in \Omega$ dependence.  In particular, two-sided skew system $F\colon X\to X$ is induced from a one-sided skew system $\hat F\colon \hat X\to \hat X$ satisfying \eqref{eq:BD}  
of \cref{sec:standing ass skew}.
 
Let $\Phi\colon \Omega \to X$, $\omega\mapsto  (\omega,p(\omega))$ be the graph of the semiconjugacy. 
Let $\mu=\Phi_{*}\beta$ be the $F$-invariant measure on the graph of $p$.
We retain our previous notation.   
Since $\mu$ is a measure on the graph of $p$, for $\beta$-a.e.~$\omega\in \Omega$, $\mu_{(\omega,x)}^{\Fsc}=\delta_{x}$ with $x=p(\omega)$ and so the fiberwise stable and unstable dimensions vanish, $\bar{d}^{s}=0$ and  $\bar{d}^{u}=0$.  
Since $\pi_2 \circ \Phi=p$, $(\pi_{2})_{*}\mu=p_{*}\beta$.  
If we consider $\beta$ with positive entropy, then $\dim (\pi_2)_*\mu>0$.  In fact, we can pick $\beta$ with $p_*\beta $ the SRB measure for $f\inv$.  
Fixing the future symbols restricts $p(\omega)$ to a local stable plaque. For this choice of $\beta$, the conditional measures on these plaques are absolutely continuous with respect to arclength.  By \cref{mainskewcor}\eqref{Cscpt}, we have  $\hat d^{s}=1$.  
Since $p_*(\beta)$ is ergodic $f$-invariant and SRB (for $f\inv$), we have 
 the a.s.~constant value of the limit  $$\lim_{r\to0} \frac{\log\mu\left(\pi_{2}^{-1}\left(B(x,r)\right)\right)}{\log r}= 
\dim(p_*\beta)>1,$$ and so 
$\wcheck d ^u>0$ by \cref{mainskewcor}\eqref{trivialpt}.  
\end{example}

\begin{example}[Holonomies incompatible with projections] 
\label{ex:dimform}
We present an example of a skew system admitting both stable and unstable holonomies for which the dimension formula in \cref{mainskewcor}\eqref{trivialpt} fails.  That is, while the system admits both stable and unstable holonomies, the unstable holonomies  are incompatible with the projection $\pi_2$.  This shows the local constantness in \cref{prop:osgimpliesshift} is needed, at least for the dimension formula of  \cref{mainskewcor}\eqref{trivialpt}.  This example was  provided by ChatGPT with light editing by the authors.

Let $\theta$ be the full shift on $\Omega=\{0,1\}^{\Zbb}$ equipped with the ergodic, invariant $(1/2,1/2)$-Bernoulli measure $\nu$.  
Write
\[
t(\omega)=\sum_{j\ge 0}\frac{\omega_j}{2^{j+1}},
\qquad
b(\omega)=\sum_{j\ge 1}\frac{\omega_{-j}}{2^j}
=\frac{\omega_{-1}}{2}+\frac{\omega_{-2}}{4}+\frac{\omega_{-3}}{8}+\cdots.
\]
Thus $b(\omega)$ is the real number whose binary digits are the past symbols of $\omega$, starting with the most recent symbol $\omega_{-1}$.
For $t\in\Rbb$, let $R_t$ denote the rotation $x\mapsto x+t\pmod 1$ on $M=\Rbb/\Zbb$. Let $\varphi^s$ be the flow of $\dot x=-\sin(2\pi x)/(2\pi)$. Note that $0$ is an attracting fixed point and $1/2$ is a repelling fixed point for $\phi^s$.  
Fix $a,\epsilon>0$ with $a+\epsilon<\log 2$ and set
\[
f_\omega=R_{t(\theta\omega)}\circ
\varphi^{a+\epsilon\omega_0 t(\theta\omega)}\circ R_{-t(\omega)}.
\]
The graph $x=t(\omega)$ is invariant under the two-sided skew-shift. Let $\mu$ denote the image of $\nu$ on this graph.  Since $(\varphi^s)'(0)=e^{-s}$, the (unique) fiber Lyapunov exponent is $-a-\epsilon/4$.  We note that the cocycle $\omega\mapsto f_\omega$ is not locally constant; indeed  on the cylinder $\{\omega_0=0\}$, writing $t=t(\omega)$, we have $f_\omega(0)=2t+\varphi^a(-t)$ modulo one. 
We may check that $\omega\mapsto f_\omega(0)$ takes uncountably many values on every conull set of $\omega$ on the cylinder $\{\omega_0=0\}$.  (E.g. the derivative of $t\mapsto 2t+\varphi^a(-t)$ is $2-(\varphi^a)'(-t)>0$ for $t$ small.)    

Because $t(\omega)$ depends only on the future coordinates, the identity provides a system of stable holonomies.  
Thus, the two-sided system is induced from a construction as in \cref{sec:skewsyst}.  
For unstable holonomies, take $\omega,\omega'$ with $\omega_j=\omega'_j$ for all $j<0$. 
The unstable holonomy from the fiber over $\omega$ to the fiber over $\omega'$ is
\[
H^u_{\omega,\omega'}
=R_{t(\omega')}\circ\varphi^{2\epsilon b(\omega)(t(\omega')-t(\omega))}\circ R_{-t(\omega)}.
\]
The difference $t(\omega')-t(\omega)$ in the flow time is taken in $\Rbb$.
To check equivariance, we verify
\[
H^u_{\omega,\omega'}=f_{\theta^{-1}\omega'}\circ
H^u_{\theta^{-1}\omega,\theta^{-1}\omega'}\circ f_{\theta^{-1}\omega}^{-1}.
\]
Indeed, expanding we have 
\begin{align*}
H^u_{\omega,\omega'}&=\bigl(R_{t(\omega')}\circ\varphi^{a+\epsilon\omega'_{-1}t(\omega')}\circ R_{-t(\theta^{-1}\omega')}\bigr)\\
&\quad\circ\bigl(R_{t(\theta^{-1}\omega')}\circ
\varphi^{2\epsilon b(\theta^{-1}\omega)(t(\theta^{-1}\omega')-t(\theta^{-1}\omega))}
\circ R_{-t(\theta^{-1}\omega)}\bigr)\\
&\quad\circ\bigl(R_{t(\theta^{-1}\omega)}\circ\varphi^{-a-\epsilon\omega_{-1}t(\omega)}\circ R_{-t(\omega)}\bigr)\\
&=R_{t(\omega')}\circ
\varphi^{\epsilon(\omega_{-1}+b(\theta^{-1}\omega))(t(\omega')-t(\omega))}
\circ R_{-t(\omega)}\\
&=R_{t(\omega')}\circ\varphi^{2\epsilon b(\omega)(t(\omega')-t(\omega))}
\circ R_{-t(\omega)}=H^u_{\omega,\omega'},
\end{align*}
where the last line uses $\omega_{-1}+b(\theta^{-1}\omega)=2b(\omega)$. In particular, $H^u_{\omega,\omega'}$ sends $t(\omega)$ to $t(\omega')$.

On the support of the measure $\{(\omega,t(\omega)):\omega\in \Omega\}$, the holonomies commute though we can check they do not commute outside of the support.  We directly check the conclusion of \cref{mainskewcor}\eqref{trivialpt} fails here.  Indeed, for each $\omega\in \Omega$, the fiberwise conditional measure $\mu_\omega =\delta_{t(\omega)}$ is a point mass.  Since $t(\omega)$ depends only on the non-negative coordinates of $\omega$, fiberwise conditional measures $\omega\mapsto \mu_\omega$ are invariant under stable holonomy.  Also, since the fiber Lyapunov exponents are all negative, the invariance principle (cf.~\cite{MR2651382,MR1204858}) ensures the fiberwise conditional measures $\omega\mapsto \mu_\omega$ are invariant under unstable holonomy.  In particular, we induce one-sided future and past skew shifts by quotienting along holonomies.  
In particular, if we define the backwards dynamical system by the quotient along unstable holonomies then,
 following the notation of \cref{mainskewcor}\eqref{trivialpt}, $\hat d^s=\wcheck d^u=0$.  On the other hand, $(\pi_2)_*\mu=t_*\nu$ is Lebesgue measure on $M=S^1$.  For $\zeta$ the trivial partition, the dimension formula in \cref{mainskewcor}\eqref{trivialpt} would assert $1=0+0$.  
Note that  setting  $\epsilon=0$, we have 
$$H^u_{\omega,\omega'}=\bigl(R_{t(\omega')}\circ R_{-t(\omega)}\bigr)$$
which depends only on the future and so stable and unstable holonomies commute.  We still have $\hat d^s=\wcheck d^u=0$ and $(\pi_2)_*\mu=t_*\nu$ is Lebesgue measure on $M=S^1$.

The explanation for the  defect of \cref{mainskewcor}\eqref{trivialpt} is that the projection $\pi_2$ is  taken along stable holonomies but not along unstable holonomies.  Thus, the projection $\pi_2$ is still, in a sense, transverse to the combinatorial unstable holonomies leading to the increase in dimension.  
\end{example}

\begin{example}[A non-exact dimensional projection]\label{ex:attracting-two-cycle}
Let $\theta$ be the full shift on $\Omega=\{0,1\}^{\Zbb}$ equipped with the ergodic, invariant $(1/2,1/2)$-Bernoulli measure $\nu$ and let $t(\omega)=\sum_{j\geq 0}\omega_j2^{-j-1}$ be as in Example~\ref{ex:dimform}.  
Fix a smooth circle diffeomorphism $g$ with an attracting two-cycle $0\mapsto \tfrac12\mapsto 0$ with $|(g^2)'(0)|<1$. 
Let $C_t\colon S^1\to S^1$, for  $0\leq t\leq 1$, be a smooth family of circle diffeomorphisms such that $C_t(0)=0$ and $t\mapsto C_t(\tfrac12)$ is a diffeomorphism onto an arc. Set
\[
f_\omega=C_{t(\theta\omega)}\circ g^{\omega_0\omega_1}\circ C_{t(\omega)}^{-1}.
\]
Let $\mu$ denote the measure on $\Omega\times S^1$ that   projects to $\nu$ and whose fiber measure over $\omega$ is $$\mu_\omega=\frac 1 2 \delta_0 + \frac 1 2 \delta_{C_{t(\omega)}(\tfrac12)}.$$ 
Then $\mu$ is $F$-invariant.  
The dynamics exchanges the two branches exactly when $\omega_0=1$ and $\omega_1=1$.

Again, we note that stable holonomies are the identity.  For $\omega,\omega'$ with $\omega_j=\omega'_j$ for all $j<0$, the unstable
holonomy in the original circle coordinates is
\[
H^u_{\omega,\omega'}
= C_{t(\omega')} \circ
  g^{\omega_{-1}(\omega'_0-\omega_0)} \circ
  C_{t(\omega)}^{-1}.
\]
Now consider $\omega= (\omega_-,\omega_+)$ and $\omega'= (\omega'_-,\omega'_+)$, opposite corners of a rectangle.  
We study the holonomy along the rectangle 
$$
H^s_{(\omega'_-,\omega_+),\,(\omega_-,\omega_+)}\circ 
H^u_{(\omega'_-,\omega'_+),\,(\omega'_-,\omega_+)}\circ H^s_{(\omega_-,\omega'_+),\,(\omega'_-,\omega'_+)}
 \circ H^u_{(\omega_-,\omega_+),\,(\omega_-,\omega'_+)}.
$$
Since the stable holonomies are the identity, the return map is
\begin{align*}
&H^u_{(\omega'_-,\omega'_+),\,(\omega'_-,\omega_+)}
 \circ H^u_{(\omega_-,\omega_+),\,(\omega_-,\omega'_+)} \\
&\qquad =
 C_{t(\omega)} \circ
 g^{(\omega_{-1}-\omega'_{-1})(\omega'_0-\omega_0)}
 \circ C_{t(\omega)}^{-1}.
\end{align*}
If $\omega_{-1}\neq \omega'_{-1}$ and 
$\omega'_{0}\neq \omega_0$, the iterate of $g$ is $g^{\pm 1}$; thus for such pairs $\omega,\omega'$, the composition of holonomy 
maps exchanges the atoms $0$ and $C_{t(\omega)}(\tfrac12)$.
A  Hopf argument shows the ergodic components of $\mu$ are invariant under stable and unstable holonomies, thus $\mu$ is ergodic.

Since $g$ has an attracting 2-cycle $0\mapsto \frac 1 2 \mapsto 0$, the Lyapunov exponent of $\mu$ is negative.  The conditional measures are invariant under stable and unstable holonomies.  In particular, 
the one-sided quotient dynamics have atomic fiber conditional measures that are invariant under holonomies and so $\hat d^s=\wcheck d^u=0$. On the other hand, $t(\omega)$ is uniformly distributed on $[0,1]$. The constant branch projects to $0$, while the moving branch projects to a measure equivalent to arclength on the arc traced by $C_t(\tfrac12), 0\le t\le 1$.  Thus $(\pi_2)_*\mu$ gives mass $\tfrac12$ to $0$ and mass $\tfrac12$ to the arc $I=\{C_t(\tfrac12):0\le t\le1\}$, with pointwise dimensions of zero at $0$ and of one at every point of $I$.  In particular, $(\pi_2)_*\mu$ is not exact dimensional, and the conclusions of  \cref{mainskewcor}\eqref{trivialpt} fail.  

\end{example}

\subsection{Comparison with prior results} 
For ergodic hyperbolic invariant measures $\mu$ invariant under $C^{2}$ diffeomorphisms, an upper bound on pointwise dimension of the form $\overline d_\mu\le d^s+d^u$ (with a related upper bound in the case $\mu$ has zero Lyapunov exponents)   is already established by  Ledrappier and Young in 
\cite[Thm.~F and \S12]{MR819557}.

For $C^{1+\alpha}$ diffeomorphisms, Barreira, Pesin, and Schmeling established the lower pointwise bound  $ d^s+d^u\le \underline d_\mu$  for ergodic hyperbolic  invariant measures (as well as an alternate argument for the upper bound) in \cite{MR1709302} and showed that ergodic hyperbolic measures admit an asymptotic almost product structure.  
Certain inclusions (especially \cite[(8) and (13)]{MR1709302}  leading to \cite[Prop.~5]{MR1709302}) seem to  not be justified by the cited construction since atoms of the cited partition $\partP$ might have no interior or leaf-wise interior and so it is unclear why balls in stable manifolds are contained in $\partP$ under the cited construction.

For  invertible $C^2$ random dynamical systems over a measurable base,  Liu and Xie study the dimension of the fiberwise conditional measures (e.g.~sample measures for stationary measures associated to i.i.d.~random walks)  for  ergodic invariant measures that  are hyperbolic along the fibers in \cite{MR2218998} and establish exact dimensionality with the expected dimension formula.  
In particular, they derive the lower bound by adapting the counting argument in \cite{MR1709302} and   replacing  \cite[Prop.~5]{MR1709302} with more detailed counting and conditional measure estimates.  
A new dimension type estimate (measuring decay of stable Riemannian balls and forwards cylinder sets that see only slow decay along stable manifolds) appears in \cite[Prop.~8.1]{MR2218998} which supplies the quantitative input needed in the replacement \cite[Lem.~6.1]{MR2218998} adapting the proof of the upper bound from \cite[Lem.~6]{MR1709302}.   We avoid this by conditioning on the entire past and invoking (a slight extension of)  \cite[ Lem.~12.4.1]{MR819557}.
One counting argument in \cite[pp.~3768--3769]{MR2218998} seems partly incomplete, since a partition atom could intersect an uncontrolled number of balls so Besicovitch bounds on multiplicity of covers by balls do not, a priori, bound the associated count of partition atoms.  We give a more detailed justification of the geometry of iterated partitions relative to stable and unstable manifolds  in \cref{sss:prebesicovitch} (using continuity of stable and unstable manifolds when restricted to sets of large measure),
to justify a Besicovitch argument in the proof of \cref{lem:compareN} using new arguments  \cref{countpart} and \cref{lem:approxa}\ref{approx4:stable-separation} and \ref{approx4:unstable-separation}.

We also highlight  the approach to the lower bound given by Ben Ovadia and Rodriguez Hertz in  \cite{MR4756121}.  Here, the authors define a notion of neutralized metric entropy and show it coincides with Brin--Katok local entropy (which coincides  with usual metric entropy for ergodic measures).  Using neutralized metric entropy, they 
establish the lower bound on pointwise dimension of an ergodic invariant measure  by covering by Bowen balls and establishing Besicovitch-type multiplicity properties of covers by Bowen balls (with arbitrarily slow exponential growth rather than boundedness of  multiplicity) which  provides a density lemma to control counting of certain covers.  Their approach recovers a lower bound of the form $ d^s+d^u\le \underline d_\mu$ even in the case the measure $\mu$ has  zero Lyapunov exponents.

In a setting closer to ours, Shu established exact dimensionality for ergodic hyperbolic invariant measures for $C^2$ endomorphisms with no singularities (i.e.~local diffeomorphisms) in \cite{MR2657815}.  Shu follows the approach of Liu and Xie for the lower bound; Shu also gives an upper bound for pointwise dimension first following the approach of Liu and Xie  in \cite[Eq.~(3.11)]{MR2657815} for hyperbolic measures 
followed by the approach of Ledrappier and Young in \cite[Lem.~4.1]{MR2657815} (which also gives a bound for measures with zero Lyapunov exponents).  

In our proof of exact dimensionality and our formulas for the pointwise  dimension, we mostly follow the counting arguments in \cite{MR2218998,MR1709302} for the lower bound and adapt the arguments of \cite{MR2657815} for the upper bound.   We could follow the approach of \cite{MR4756121} to lower bound pointwise dimension but choose instead to adapt \cite{MR2218998} as that approach follows similar arguments used to establish the upper bound (following  \cite[\S12]{MR819557}).  

One novelty of our arguments here and also in the paper \cite{BRH} is that we never use local entropy along stable or unstable foliations and never compute entropy in terms of the decay (relative to conditional measures) of stable or unstable Bowen balls.  This allows us to treat contributions to entropy from the combinatorial and manifold directions simultaneously and to simplify certain arguments from \cite{MR2218998,MR819557}.

\section{Entropy  in abstract ergodic theory} 
We review notions of conditional metric entropy in the setting of abstract measure-preserving transformations.  We also present a number of elementary properties, including an elementary proof of a Shannon--McMillan--Breiman theorem for (standard) conditional entropy.  

\subsection{Standing notation}
Let $(X,\sigB, \mu)$ be a standard probability space.  Let $T\colon (X,\mu)\to (X,\mu)$ be a measure-preserving transformation.  
We say a $\sigma$-algebra $\sigA$ is \emph{sub-invariant} if $T\inv \sigA\subset \sigA$ and \emph{\siname} if $T\inv \sigA= \sigA$. 
Let $\zeta$ be a measurable partition of $X$.  We similarly say $\zeta$ is sub-invariant if $T\inv \zeta\prec \zeta$ and invariant if $T\inv \zeta = \zeta.$

 Given a measurable partition  $\xi$ of $(X,\mu)$ and non-empty $S\subset \N$,  write 
\begin{enumerate}
\item $T\inv \xi:= \{T\inv (C): C\in \xi\}$;
\item  $\xi^{\emptyset} := \trivpart$;
\item $ \xi^S:= \bigvee_{i\in S} T^{-i} \xi$
 and so $\xi^{[0,n)}= \bigvee_{j=0}^{n-1}  T^{-j}(\xi) $; 
\item $\xi^- :=\xi^{[0,\infty)}= \bigvee_{j=0}^\infty  T^{-j}(\xi)$. 
\end{enumerate}
Also write \begin{enumerate}[resume]
\item $\xi^\emptyset =\trivpart$. 
\end{enumerate}
If $T$ is an automorphism and $S\subset \Z$,  we  also write 
\begin{enumerate}[resume]
\item $ \xi^S:= \bigvee_{i\in S} T^{-i} \xi$; 
\item $\xi^\Z= \bigvee_{j=-\infty}^\infty  T^{j}(\xi)$;
\item $\xi^+= \bigvee_{j=0}^\infty  T^{j}(\xi)$. 
\end{enumerate}
We note that for any measurable partition $\xi$, the $\sigma$-algebra  $\sigma(\xi^-)$ is sub-invariant and, if $T$ is an automorphism, we have that $\sigma(\xi^\Z)$ is \siname.  
We use similar notation for preimages of sub-$\sigma$-algebras. 

Given a sub-$\sigma$-algebra $\sigA$ of $(X,\mu)$, we let $\{\mu_x^\sigA\}$ denote a family of conditional probability measures.  Given a measurable partition $\xi$ of $(X,\mu)$ we define the 
conditional information function $$I\bigl( \xi \mid \sigA\bigr)(x)=
- \log\bigl({ \mu }^{ \sigA}_{ x} \bigl(  \xi ( x) \bigr)   \bigr)  $$
and the conditional entropy $$H_\mu(\xi\mid \sigA) = \int I\bigl( \xi \mid \sigA\bigr)(x) \, d \mu(x).$$

\subsection{Notions of entropies}
For the remainder of this section, let $\sigF$ be a sub-invariant sub-$\sigma$-algebra  for $T\colon (X,\mu)\to (X,\mu)$.  
\subsubsection{Standard metric entropy}
Given a measurable partition $\xi$ of $(X,\mu)$, the \emph{entropy} of $T$ with respect to $\xi$ is $$h_\mu(T,\xi):= H_\mu(\xi \mid T\inv \xi^-).$$
When $H_\mu(\xi)<\infty$, we alternatively have 
\begin{equation}\label{eq:decaydefn} h_\mu(T,\xi)= \lim_{n\to \infty} \frac 1 n H_\mu(\xi^{[0,n)}).\end{equation}

\subsubsection{Standard conditional entropy}\label{standard conditional entropy}
Motivated by \eqref{eq:decaydefn}, we define (standard) conditional entropy relative to a sub-invariant sub-$\sigma$-algebra $\sigF$.  We note that unlike standard entropy, this will only be defined for measurable partitions $\xi$ of $(X,\mu)$ with $H_\mu(\xi\mid \sigF)<\infty$.  

Let $\xi$ be a measurable partition of $(X,\mu)$ with $H_\mu(\xi\mid \sigF)<\infty$.  
We define the   (standard) \emph{conditional entropy} of $T$ with respect to $\xi$ conditioned on $\sigF$ to be $$h_\mu(T,\xi \mid \sigF):= \lim_{n\to \infty} \frac 1 n H_\mu\left(\xi^{[0,n)}  \mid \sigF\right)=   \inf_{n} \frac 1 n H_\mu\left(\bigvee_{i=0}^{n-1} T^{-i}\xi \mid \sigF\right) $$
where the infimum in the final equality follows from subadditivity.  
We define the   (standard) \emph{conditional entropy} of $T$  conditioned on $\sigF$ to be $$h_\mu(T\mid \sigF):= \sup \{h_\mu(T,\xi \mid \sigF)\}$$
where the supremum is taken over all measurable partitions $\xi$ of $(X,\mu)$  with $H_\mu(\xi\mid \sigF)<\infty$.

We note that 
when  $\sigF$ is invariant (rather than only sub-invariant), we have an alternative computation of the (standard) conditional entropy.
\begin{lemma}[{c.f.~\cite[Fact 2.3.4]{MR2809170}}]\label{lem:invariantSTDisTRANS}
If $\sigF$ is invariant and if $H_\mu(\xi\mid \sigF)<\infty$ then 
$$h_\mu(T,\xi \mid \sigF)=H_\mu (\xi\mid T\inv \xi^- \vee \sigF).$$
Moreover, if $T$ is an automorphism and if $\sigF$ is invariant  then $$h_\mu(T,\xi \mid \sigF)=h_\mu(T\inv,\xi \mid \sigF).$$
\end{lemma}
For the second conclusion in \cref{lem:invariantSTDisTRANS}, note  by invariance of $\mu$ and of $\sigF$ that 
$$H_\mu(\xi^{[0,n)} \mid \sigF) =H_\mu ( T^{n-1}(\xi^{[0,n)}) \mid \sigF) =H_\mu ( \xi^{(-n,0]} \mid \sigF).$$

\begin{corollary}\label{cor:STDentALT}
If $\sigF$ is invariant then 
$$h_\mu(T \mid \sigF)=\sup \{H_\mu (\xi \mid T\inv \xi^- \vee \sigF)\}$$
where the supremum is taken over all measurable partitions of $(X,\mu)$.  
\end{corollary}
\begin{proof}
We need  only show that if  $h_\mu(T \mid \sigF)<\infty$ and if $\xi$ is an arbitrary measurable partition of $(X,\mu)$, then 
$$H_\mu (\xi \mid T\inv \xi ^-\vee \sigF)\le h_\mu(T \mid \sigF).$$
For this, consider a sequence of finite partitions $\eta_k$ with $\eta_k\nearrow \xi$ as ${k\to\infty}.$
Then by \cref{lem:invariantSTDisTRANS}
\begin{align*}
H_\mu (\xi \mid T\inv \xi^- \vee \sigF)&= \lim _{k\to \infty} H_\mu (\eta_k \mid T\inv \xi ^-\vee \sigF)\\&\le 
\lim _{k\to \infty} H_\mu (\eta_k \mid T\inv \eta_k^- \vee \sigF)\le h_\mu(T \mid \sigF).  \qedhere
\end{align*}
\end{proof}

\begin{remark}
In \cite{BRH}, given  a \emph{sub-invariant} $\sigF$, the quantity  $$H_\mu (\xi\mid T\inv \xi ^-\vee \sigF)$$  was called the transverse conditional entropy of $T$ with respect to $\xi$ conditioned on $\sigF$ and the quantity $$\sup _{\xi } H_\mu (\xi\mid T\inv \xi ^-\vee \sigF)$$  was called the transverse conditional entropy of $T$  conditioned on $\sigF$.  (Here the supremum is over all measurable partitions $\xi$.) 
However, in the present paper, we will only ever use the expression $H_\mu (\xi\mid T\inv \xi ^-\vee \sigF)$   in the setting where $\sigF$ is (completely) invariant; thus the concept of transverse conditional entropy will not be used in the remainder.  
\end{remark}

\subsubsection{Refinement entropy relative to $\sigF$}\label{Def:refent}
We present another notion of entropy relative to a sub-invariant sub-$\sigma$-algebra $\sigF$.  Following \cite{BRH}, we define the \emph{refinement entropy} of $T$ relative to  $\sigF$  to be
$$\maxent_\mu(T\mid \sigF):= \sup \{ h_\mu(T,\xi\vee \Xi( \sigF))\}$$
where the supremum is taken over all measurable partitions $\xi$ of $(X,\mu).$ 
Alternatively, we have $$\maxent_\mu(T\mid \sigF)= \sup_{\sigF \subset \sigma(\xi)}  h_\mu(T,\xi).$$

\subsubsection{Furstenberg entropy}\label{sss:furst}
Recall that given two probability measures $\nu$ and $\eta$ on a Polish space  $Y$, 
the  Lebesgue--Radon--Nikodym decomposition  of $\nu$ relative to $\eta$ is the decomposition   $$\nu= \rho \eta + \nu_s$$ where $ \rho :=\frac{d\nu}{d\eta}\in L^1( \eta)$  is the Radon--Nikodym derivative of the component of $\nu$ absolutely continuous with respect to $\eta$   and $\nu_s\perp \eta$ is the component of $\nu$ singular with respect to $\eta$.  
Note the positive part of the function $\log \bigl({\frac{d\eta}{d\nu}}\bigr) $ is always $L^1(\nu)$.  
The \emph{Kullback--Leibler divergence} between $\nu$ and $\eta$ is the quantity 
\begin{equation*}\KL(\nu \parallel \eta)= -\int \log \left({\frac{d\eta}{d\nu}}\right) \, {d\nu} .\end{equation*}

Returning to $T\colon (X,\mu)\to (X,\mu)$ and a sub-invariant sub-$\sigma$-algebra $\sigF$, 
we define the  \emph{Furstenberg entropy} of $T$ relative  to  $\sigF$ and  with respect to  $\mu$ to be
$$\furstent_\mu(T\mid \sigF)=\int \KL(T_*\mu_x^\sigF\parallel \mu_{T(x)}^\sigF) \, d \mu(x).$$
We have $0\le \furstent_\mu(T\mid \sigF)\le \infty$.

\subsection{Properties of entropies}

We collect various properties of our entropy definitions that will be used throughout.

\begin{lemma}[{\cite[Thm.\  0.4.4]{MR1369243}}]\label{lem:inv-entropy}
Suppose that $T\colon (X,\mu)\to (X,\mu)$ is an automorphism of the measure space $(X,\mu)$.  Then for any measurable partition $\xi$ with $H_\mu(\xi \mid \sigF)<\infty$, we have $$h_\mu(T,\xi \mid \sigF) = h_\mu(T,\xi \mid \sigF^+).$$
\end{lemma}

We also have the following direct corollary of \cref{lem:invariantSTDisTRANS}
\begin{corollary}[{\cite[Fact 2.3.4]{MR2809170}}]\label{equalifinv}
Suppose  that $\sigF$ is \siname\ for $T\colon X\to X$.  
Then 
$$\maxFent  =  h_\mu(T\mid \sigF).$$
\end{corollary}

We also have the following additivity formula for standard (conditional) entropy.  
\begin{lemma}[{cf. \cite[Facts 2.4.2 + 2.3.7]{MR2809170} or \cite[7.7]{MR0217258}}]
\label{lem:natext}
Let $\sigF$ be a sub-invariant sub-$\sigma$-algebra.  
Let $\eta$ and $\xi$ be measurable partitions with $H_\mu(\xi \mid \sigF)<\infty$ and $H_\mu(\eta \mid \xi\vee \sigF)<\infty$.  Then 
$$h_\mu(T,\eta \vee \xi \mid \sigF) = h_\mu(T,\xi  \mid \sigF)  + h_\mu(T,\eta \mid  \xi^-\vee \sigF) . $$
Moreover, if  $T$ is an automorphism then 
$$h_\mu(T,\eta \vee \xi \mid \sigF) = h_\mu(T,\xi  \mid \sigF)  + h_\mu(T,\eta \mid  \xi^\Z\vee \sigF^+) . $$

If  $T$ is an automorphism and if $H_\mu(\eta\vee \xi\mid T\inv \xi^-)<\infty$
then $$h_\mu(T, \xi\vee \eta) = h_\mu(T, \xi) + H_\mu(\eta\mid T\inv \eta^-\vee \xi^\Z).$$
\end{lemma}

We have the following immediate corollaries of the unconditional version of \cref{lem:natext}.  
\begin{corollary}\label{cor:localentmax}
For any measurable 
partitions $\xi, \eta$ with $H_\mu(\eta\mid \xi)<\infty$,
$$h_\mu(T,\xi \vee \eta) \ge h_\mu(T,\xi).$$
Moreover if $T$ is an automorphism and if $\xi$ generates (i.e.~if $\xi^\Z=\ptpart)$ and if $h_\mu(T,\xi)<\infty$ then 
$$h_\mu(T,\xi \vee \eta) = h_\mu(T,\xi).$$
\end{corollary}

We also recall the standard fact: 
\begin{claim}\label{claim:entshift}
If $\zeta$ is a $T$-invariant measurable partition then for any measurable partitions $\xi,\eta$ and any $n\ge 0$,  
$h_\mu(T,\xi\vee\eta\vee \zeta) = h_\mu(T,T^{-n}\xi\vee\eta\vee \zeta).$
\end{claim}

\subsection{Conditional Shannon--McMillan--Breiman theorem}
Let $T\colon (X,\mu)\to (X,\mu)$ be an ergodic  measure-preserving system of the standard probability space.  
Let $\eta$ be a partition and $\sigF$ be a subinvariant sub-$\sigma$-algebra on $X$ with  $H_{\mu}(\eta\mid \sigF)<\infty$.  One can study the exponential decay rate:
\begin{equation}\label{eq:SMB??}\lim_{n\to \infty} -\frac{1}{n}\log \mu^{\sigF}_{x}\left(\eta^{[0,n)} (x) \right).\end{equation}
When $\sigF$ is the trivial $\sigma$-algebra, the classical Shannon--McMillan--Breiman theorem asserts \eqref{eq:SMB??} converges to the metric entropy $h_\mu(T, \eta)$ of $T$ with respect to $\eta$.

On \cite[p.~60]{MR884892}, Y.~Kifer attributes the following to B.~Weiss: ``{\it For any good notion of entropy, one
should be able to prove some kind of Shannon--McMillan--Breiman theorem.}'' 
For invariant  sub-$\sigma$-algebras $\sigF$, a version of the  Shannon--McMillan--Breiman for conditional metric entropy appeared many times including in \cite[Thm.~1.5]{MR884892} and \cite{MR1181382}.  

However, when $\sigF$ is assumed to be only sub-invariant, the only  statement and proof of a conditional Shannon--McMillan--Breiman known to the authors  
appears in  Downarowicz's book, \cite[Thms.~3.3.7 and B.0.1]{MR2809170}.

While Downarowicz's proof in \cite[Thm.~B.0.1]{MR2809170} is complete and fully self-contained, though complicated.  
 We present an elementary proof of the conditional Shannon--McMillan--Breiman theorem (assuming the Chung--Neveu maximal inequality,  martingale convergence theorems, and Maker's ergodic theorem for dominated sequences of functions \cite{MR2028}), following the standard proof of the classical   Shannon--McMillan--Breiman theorem presented in \cite[\S2.2]{MR0262464} or \cite[\S6.2]{MR1073173}.

\begin{proposition}[Conditional Shannon–McMillan–Breiman theorem] \label{SMB}
Assume that $\mu$ is an ergodic $T$-invariant probability measure.  
Let $\sigF$ be a subinvariant sub-$\sigma$-algebra on $X$ and let  $\eta$ be a measurable partition with  $H_{\mu}(\eta\mid \sigF)<\infty$.  Then for $\mu$-almost every $x\in X$,
\[\lim_{n\to \infty} -\frac{1}{n}\log \mu^{\sigF}_{x}\bigl(\eta^{[0,n)} (x) \bigr)=h_{\mu}(T,\eta \mid \sigF).\] 
\end{proposition}

\begin{proof}[Proof of \Cref{SMB}]
We first note that by passing to the natural extension, we may assume without loss of generality that $T$ is an automorphism.  (See \cite[Prop.~3.14]{BRH}) 
 
We expand the conditional information function
 \begin{align*}
 -\log\bigl({ \mu }^{ \sigF}_{ x} \bigl(  \eta ^{[0,n)}( x) \bigr)   \bigr)  
 &= I\bigl( \eta ^{[0,n)}\mid \sigF\bigr)(x)\\
  &= I\bigl( \eta \mid \sigF\bigr)(x)+  I\bigl(\eta ^{[1,n)}\mid  \eta \vee \sigF\bigr)(x)\\
    &= I\bigl( \eta \mid \sigF\bigr)(x)+  I\bigl(\eta ^{[0,n-1)}\mid T( \eta) \vee T( \sigF)\bigr)(T(x))\\
&= \sum_{j=0}^{n-1}  I\bigl(\eta \mid  T(\eta^{(-j,0]} )\vee T^j( \sigF)\bigr)(T^j(x)).
 \end{align*}
(Recall the convention that  $T(\eta^{(-0,0]})$ is the trivial partition).
For $j\in\{0,1,2,\dots\}$, let 
 \begin{align*}
 g_j(x) &:= I\bigl(\eta \mid  T(\eta^{(-j,0]} )\vee T^j( \sigF)\bigr)(x) \\
 g_*(x) &:= \sup_{j\ge 0}\left\{I\bigl(\eta \mid  T(\eta^{(-j,0]} )\vee T^j  (\sigF)\bigr)(x)\right\}\\
 g_\infty(x)& := I\bigl(\eta \mid  T(\eta^+ )\vee  \sigF^+\bigr)(x). 
  \end{align*}
By the (conditional) Chung--Neveu maximal lemma (see \cite[Lem.\  2.1]{MR0262464}), we have $$\int g_*\, d\mu\le H(\eta\mid \sigF)+1.$$  Thus, by martingale convergence,  we have dominated a.e.\ convergence $g_j \to g_\infty$.  
  Applying  Maker's ergodic theorem for dominated sequences of functions \cite{MR2028} and using that $\sigF^+$ is invariant and thus \cref{lem:invariantSTDisTRANS} applies, we have 
   \begin{align*}
\lim_{n\to\infty}-\frac 1 n \log\bigl({ \mu }^{ \sigF}_{ x} \bigl(  \eta ^{[0,n)}( x) \bigr)   \bigr)  &= \int  g_\infty \, d \mu = H_{\mu}(\eta \mid T\eta^{+} \vee \sigF^{+})\\&= {h}(T^{-1} ,\eta \mid \sigF^+)
\\&=h(T,\eta \mid \sigF^+) = h(T ,\eta \mid \sigF).\qedhere\end{align*}
 \end{proof}

We record the following special case of \cref{SMB} when $\sigF=\sigma(\xi)$ for a sub-invariant partition that is forwards generating.  
\begin{corollary}\label{smbzero}
In the setting of  \Cref{SMB}, suppose that $T$ is an automorphism.  Let $\xi$ be a sub-invariant measurable partition with $\xi^+=\ptpart$.  Let $\eta$ be a measurable partition with $H_\mu(\eta\mid \xi)<\infty$.
Then for $\mu$-almost every $x\in X$,
\[\lim_{n\to \infty} -\frac{1}{n}\log \mu^{\xi}_{x}\bigl(\eta^{[0,n)} (x) \bigr)=0.\]
\end{corollary}

\subsection{Further SMB-type theorems}
Repeating proofs  nearly identical to the proof of \cref{SMB}, we deduce a number of SMB-like estimates on the decay rates of atoms of iterated partitions.

First, we have the following whose proof  is given in  \cite[Lem.~9.3.1]{MR819557} or \cite[Lem.~3.9]{BRH}.  
\begin{lemma}\label{lemma:pastfinitepartabstract} 
Suppose that $T\colon (X,\mu)\to (X,\mu)$ is an ergodic, measurable automorphism. 
Let $\zeta$ be a  sub-invariant measurable partition of $(X,\mu)$ with $h_\mu(T, \zeta)<\infty$ and let  $ \partP$ be a measurable partition of $(  X,  \mu)$ with $H_\mu( \partP\mid  \zeta)<\infty$. 
Then for $ \mu$-a.e.\ $x$,
$$\lim_{n \to \infty} - \frac 1 n \log  { \mu }^{  \zeta}_{x} \left( (    \zeta \vee   \partP) ^{(-n,0]}(x) \right)  = h_{ \mu }(   T, \partP \vee \zeta ).$$
\end{lemma}

Similarly, we have the following whose proof is again nearly identical to the proof of \cref{SMB}; we present the full proof for completeness.
\begin{lemma}\label{lemma:pastfinitepartabstract-half} 
Suppose that $T\colon (X,\mu)\to (X,\mu)$ is an ergodic, measurable automorphism. 
Let $\zeta$ be a  sub-invariant measurable partition of $(X,\mu)$ with $h_\mu(T, \zeta)<\infty$ and let  $ \partP$ be a measurable partition of $(  X,  \mu)$ with $H_\mu( \partP\mid  \zeta)<\infty$. 
Then for $ \mu$-a.e.\ $x$,
$$\lim_{n \to \infty} - \frac 1 n \log  { \mu }^{ ( \zeta\vee \partP)^{[0,n)}}_{x} \left( (    \zeta \vee   \partP) ^{(-n,0)}(x) \right)  =h_{ \mu }(   T, \partP \vee \zeta ).$$
\end{lemma}
\begin{proof}
We expand 
\begin{align*} 
- \log  { \mu }^{ ( \zeta\vee \partP)^{[0,n)}}_{x} &\bigl( (    \zeta \vee   \partP) ^{(-n,0)}(x) \bigr)
\\& = I_{\mu}\bigl((\zeta\vee\Pcal)^{(-n,0)}\mid (\zeta\vee\Pcal)^{[0,n)}\bigr)(x)
\\& = \sum_{k=1}^{n-1} I_{\mu}\bigl(T^{k}(\zeta\vee\Pcal)\mid (\zeta\vee\Pcal)^{(-k,0)}\vee (\zeta\vee\Pcal)^{[0,n)}\bigr)(x)
\\& = \sum_{k=1}^{n-1} I_{\mu}\bigl(\zeta\vee\Pcal \mid  T^{-1}(\zeta\vee\Pcal)^{[0, n+k-1)}\bigr)(T^{-k}(x))
\\& = \sum_{k=1}^{n-1} I_{\mu}\bigl(\zeta\vee\Pcal \mid  T^{-1}\zeta\vee T^{-1}\Pcal^{[0, n+k-1)}\bigr)(T^{-k}(x)).
\end{align*}
For $n,j\in\{0,1,2,\dots\}$, let 
 \begin{align*}
 g_{n,0}(x) &:= 0\\
 g_{n,j}(x) &:= I\bigl(\zeta\vee\Pcal \mid  T^{-1}\zeta\vee T^{-1}\Pcal^{[0,n+ j-1)} \bigr)(x), \text{ for $j\ge1$} \\
 g_*(x) &:= \sup_{j\ge 0}\left\{I\bigl(\zeta\vee\Pcal \mid T^{-1}\zeta\vee T^{-1}\Pcal^{[0, j)}\bigr)(x)\right\}\\
 g_\infty(x)& := I\bigl(\zeta\vee\Pcal \mid  T\inv \zeta\vee T\inv \Pcal ^- \bigr)(x). 
  \end{align*}
By the Chung--Neveu maximal lemma (see \cite[Lem.\  2.1]{MR0262464}), we have $g_*\in L^1(\mu)$ and thus by martingale convergence,  we have dominated a.e.\ convergence $g_{n,j} \to g_\infty$ as $\min \{n,j\}\to \infty$.  
  Applying  Maker's ergodic theorem for dominated sequences of functions \cite{MR2028},
  we have
  \begin{align*}
 \lim_{n \to \infty} -\frac 1 n  &\log  { \mu }^{ ( \zeta\vee \partP)^{[0,n)}}_{x} \bigl( (    \zeta \vee   \partP) ^{(-n,0)}(x) \bigr)
\\  &=  \lim_{n \to \infty} \frac 1 n  \sum_{k=0}^{n-1} g_{n,k}(T^{-k}(x))
 = \int  g_\infty(x) \, d \mu(x) 
  = h_{ \mu }(   T, \partP \vee \zeta ).\qedhere
  \end{align*}\end{proof}

Finally, a simple version of the same arguments gives the following.   
\begin{claim}\label{lemma:pastfinitepartabstract-full} 
Suppose that $T\colon (X,\mu)\to (X,\mu)$ is an ergodic, measurable automorphism. 
Let $\zeta$ be a  sub-invariant measurable partition of $(X,\mu)$  and let  $ \partP$ be a measurable partition of $(  X,  \mu)$ with $H_\mu( \partP\mid  \zeta)<\infty$.  Let $\xi$ be a measurable partition with $T^n\xi\nearrow \ptpart$ and $h_\mu(T, \xi\vee \zeta)<\infty$.  
Then for $ \mu$-a.e.\ $x$,
$$\lim_{n \to \infty} - \frac 1 n \log  { \mu }^{(\xi \vee  \zeta\vee \partP)^{(-n,0]}}_{x} \left(  (  \xi\vee  \zeta \vee   \partP) ^{(0,n)}(x) \right)  =0.$$

\end{claim}
\begin{proof} We have 
\begin{align*}
 \frac 1 n I_{\mu}((\xi\vee &\zeta\vee\Pcal)^{(0,n)}\mid\ (\xi\vee \zeta\vee\Pcal)^{(-n,0]})(x)
\\& =\frac 1 n\sum_{k=1}^{n-1}I_{\mu}(T^{-k}(\xi\vee\zeta\vee\Pcal)\mid(\xi\vee \zeta\vee\Pcal)^{(-n,k
)})(x)
\\& =\frac 1 n \sum_{k=1}^{n-1}I_{\mu}(\xi\vee \zeta\vee\Pcal\mid T^{n+k-1} \xi\vee (\zeta\vee\Pcal)^{(-n-k,0)})(T^k(x))
\end{align*}
Since $T^n\xi\nearrow \ptpart$, the Chung--Neveu maximal inequality and Maker's ergodic theorem for dominated sequences of functions \cite{MR2028} imply this converges to $0$.  
\end{proof}

\subsection{Revisiting Lemma 12.4.1  of \cite{MR819557}}
A key estimate used to bound dimension from above in  \cite{MR819557} is \cite[Lem.~12.4.1]{MR819557}.  
The proof of \cite[Lem.~12.4.1]{MR819557} is  correct though somewhat compressed with direction of time  reversed in one hypothesis.   For completeness, we give a full proof of a slight extension of \cite[Lem.~12.4.1]{MR819557}.  

Throughout,  we assume that $T$ is an automorphism of $(X,\mu)$.    

\begin{proposition}[{\cite[Lem.~12.4.1]{MR819557}}] \label{LY12.4.1} Suppose that $T$ is an ergodic automorphism of a probability space $(X,\mu)$.
 Let $\zeta$ be a sub-invariant measurable partition with 
 $h_{\mu}(T \mid \zeta)<\infty$.   
 Let $\xi$ be a partition that satisfies 
 \begin{enumerate}
 \item $T\xi\prec \xi$ (i.e., $\xi$ is sub-invariant for $T^{-1}$) and
 \item $H_\mu(\xi\mid T \xi\vee \zeta^+) = h_\mu(T\inv \mid \zeta^+).$ 
\end{enumerate}
Let $ \calP$ be a measurable partition of $(  X,  \mu)$ with $H_\mu( \calP\mid  \zeta)<\infty$.   Then for $ \mu$-a.e.\ $x$,
$$\lim_{n \to \infty} - \frac 1 n \log  { \mu }^{ \zeta^+\vee \partP^+\vee \xi}_{ x} \left( (   \zeta \vee   \calP ) ^{[0,n)}( x) \right)  = h_{ \mu }(   T, \partP \mid \zeta^+).$$
 \end{proposition}
Recall that since $\zeta^+$ is invariant, we have $$ h_\mu(T\inv \mid \zeta^+)=\maxent _\mu(T^{\pm1} \mid \zeta^+) =h_\mu(T \mid \zeta^+)=h_\mu(T \mid \zeta).$$

 \begin{proof}[Proof of \Cref{LY12.4.1}]
We have 
 \begin{align*}
 -\log\bigl({ \mu }^{  \zeta^+\vee \partP^+\vee \xi}_{ x} &\bigl(  (\zeta\vee \calP) ^{[0,n)}( x) \bigr)   \bigr)  
 =: I\bigl( (\zeta\vee \calP) ^{[0,n)}\mid  \zeta^+\vee \partP^+\vee \xi\bigr)(x)\\
&  = I\bigl(\zeta^{[0,n)} \mid \zeta^+\vee \partP^+\vee \xi\bigr)(x) 
+ I\bigl(\calP ^{[0,n)}\mid \zeta^+\vee \partP^+\vee \xi \bigr)(x)\\
&  = 0+ I\bigl( \calP ^{[0,n)}\mid \zeta^+\vee \partP^+\vee \xi \bigr)(x)\\
&  =  I\bigl(\calP ^{[0,n)}\mid \zeta^+\vee \partP^+\vee \xi \bigr)(x).  
 \end{align*}
We further expand
 \begin{align*}
  I\bigl(  \calP ^{[0,n)}\mid \zeta^+\vee & \partP^+\vee \xi \bigr)(x)
  \\
  &= I\bigl(   \calP \mid \zeta^+\vee  \partP^+\vee \xi \bigr)(x) + I\bigl( \calP ^{[1,n)} \mid \zeta^+\vee  \partP^+\vee \xi \bigr)(x) \\
&= 0+  I\bigl(   T\inv  \calP ^{[0,n-1)} \mid \zeta^+\vee  \partP^+\vee \xi \bigr)(x) \\
&=  I\bigl(    \calP ^{[0,n-1)} \mid \zeta^+\vee  T \partP^+\vee T \xi \bigr)(Tx) \\
&= 
\sum_{j=1}^{n-1}  I\bigl(    \calP  \mid \zeta^+\vee  T \partP^+\vee T ^j\xi \bigr)(T^jx).
 \end{align*}

Let $\sigG := \sigma (\bigwedge_{j\ge 0}  ( \zeta^+\vee  T \partP^+\vee T ^j\xi) )$ and note that $\sigG$ is sub-invariant for $T\inv$.  
For $j\in\{0,1,2,\dots\}$, let 
 \begin{align*}
 g_0(x)&:= 0\\
 g_j(x) &:= I\bigl( \partP \mid\zeta^+\vee  T \partP^+ \vee T^{j}\xi \bigr) (x) , \text{ for $j\ge 1$}\\
 g_*(x) &:= \sup\{ g_j(x)\}\\
 g_\infty(x)& := I\bigl(\partP \mid   \sigG \bigr)(x). 
  \end{align*}
The Chung--Neveu maximal lemma (see \cite[Lem.\  2.1]{MR0262464}) implies 
$$\int \max_{0\le j\le N} g_j \, d \mu\le H_\mu(\partP\mid \zeta ^+)+1$$ for every $0\le N$.  Monotone convergence then gives $g_*\in L^1(\mu)$ and thus by martingale convergence,  we have dominated a.e.\ convergence $g_j \to g_\infty$.  
Applying  Maker's ergodic theorem for dominated sequences of functions \cite{MR2028}, we have the upper bound
\begin{align*}
\lim_{n\to \infty}  \frac 1 n  I\bigl(  \calP ^{[0,n)}\mid \zeta^+\vee & \partP^+\vee \xi \bigr)(x)\notag
= \int  g_\infty(x) \, d \mu(x)\notag\\
&= H_\mu (\partP \mid \sigG) \le H_\mu (\partP \mid\zeta^+\vee T \partP^+)=   h_\mu (T\inv, \partP \mid \zeta^+).\notag
\end{align*}

On the other hand, using \cref{claim:entshift,lem:natext}
\begin{align}
 H\bigl(\partP & \mid   \sigG \bigr)=\lim_{n\to \infty} H\bigl(\partP \mid  \zeta^+\vee  T \partP^+\vee T ^n\xi\bigr)\notag
 \\
  &=\lim_{n\to \infty}\left( H\bigl(T^{n-1}\xi \vee \partP \mid  \zeta^+\vee  T \partP^+\vee T ^n\xi\bigr) -
  H\bigl(T^{n-1} \xi  \mid  \zeta^+\vee   \partP^{+}\vee T^n \xi \bigr)
  \right)
  \notag
  \\
  &=\lim_{n\to \infty}\left( H\bigl(\xi \vee \partP \mid  \zeta^+\vee  T \partP^+\vee T \xi\bigr) -
  H\bigl(  \xi  \mid  \zeta^+\vee    \partP^{(-\infty,n)}\vee T \xi \bigr)
  \right) \notag \label{thisstep}
  \\&=
  H\bigl(\xi \vee \partP \mid  \zeta^+\vee  T \partP^+\vee T \xi\bigr) -
  H\bigl(  \xi  \mid ( \zeta^+\vee    \partP^{\Z}\vee T \xi)\bigr)
  \notag
  \end{align}
  
Take a sequence of finite partitions $\partQ_n\nearrow \xi$.  Then, by applying \cref{lem:natext} twice (with the dynamics $T\inv$), we have
\begin{align*} 
 H\bigl(\partP & \mid   \sigG \bigr)
  \\&=
  H\bigl(\xi \vee \partP \mid  \zeta^+\vee  T \partP^+\vee T \xi\bigr) -\lim_{n\to \infty}
  H\bigl(  \partQ_n  \mid ( \zeta^+\vee    \partP^{\Z}\vee T \xi)\bigr)
  \\&\ge
  H\bigl(\xi \vee \partP \mid  \zeta^+\vee  T \partP^+\vee T \xi\bigr) -\lim_{n\to \infty}
  H\bigl(  \partQ_n  \mid ( \zeta^+\vee    \partP^{\Z}\vee T \partQ_n^+)\bigr)
 \\&=
  h_\mu\bigl(T\inv, \xi \vee \partP \vee \zeta^+)- \lim_{n\to \infty}   \left ( h_\mu(T\inv, \partQ_n\vee \zeta^+\vee \partP) -  h_\mu(T\inv  , \zeta^+\vee \partP)
      \right)\\ 
 &= h_\mu\bigl(T\inv, \xi  \vee \zeta^+)+ H(\partP\mid T\partP^+\vee \xi^\Z\vee \zeta^+) \notag\\ &\quad\quad
- \lim_{n \to \infty} h_\mu(T\inv, \partQ_n \vee \zeta^+\vee \partP) + h_\mu(T\inv,  \zeta^+\vee \partP)   
\\&       \ge  h_\mu(T\inv,  \zeta^+\vee \partP)   =h_{ \mu }(   T\inv , \partP \mid \zeta^+).  \notag\qedhere
  \end{align*}
\end{proof}

\subsection{Abstract results on differentiation and dimension} \label{sec:diffanddim}
Throughout the sequel, we will consider $X=\Omega\times M$ where $\Omega$ is an abstract measure space and $M$ is a smooth manifold.   Let $\pi_2\colon X\to M$ denote the projection onto the second coordinate.  

Given a probability measure $\mu$ on $X$, we frequently restrict $\mu$ to subsets $\Gamma\subset X$ of positive (often nearly full) measure.  Applying Besicovitch differentiation on the manifold $M$, dimension properties of $\mu$ and $\restrict{\mu}{\Gamma}$ coincide when projected to $M$.  

\begin{proposition}\label{prop:Lebdensity}
Let  $X=\Omega\times M$, let $\mu$ be a probability measure on $X$, and let $\sigA$ be a sub-$\sigma$-algebra.  Let $\eta$ be a measurable partition of $X$ with either $\eta$ countable or $H(\eta\mid \sigA)<\infty$ and let $\Gamma\subset X$ be a subset with $\mu(\Gamma)>0$.  

Then, given any $0<\epsilon\le  1 $ there are $C\ge 1$, $r_0>0$, and $\Gamma'\subset \Gamma$ with $\mu(\Gamma')\ge \mu(\Gamma)-2\epsilon$ such that for every $(\omega,x)\in \Gamma'$ and for all $r<r_0$,
\begin{equation}\label{eq:lowerRND}\frac{1}{C} \mu^{\sigA}_{(\omega,x)}\bigl(\pi_{2}^{-1} (B(x,r)) \bigr)\le 
\mu^{\sigA}_{(\omega,x)}\bigl(\pi_{2}^{-1}\bigl( B(x,r) \bigr)\cap \Gamma \cap \eta(\omega,x)\bigr). 
 \end{equation}

In the case that $\eta\prec \Xi(\sigA)$, then $\Gamma'$ and $r_0$ can be chosen such  that for every  $(\omega,x)\in \Gamma'$ and all $r<r_0$, \eqref{eq:lowerRND}  holds with $C=\left(1-2 \epsilon^{-1}\bigl(1-\mu(\Gamma)\bigr)\right)\inv$ (provided $\mu(\Gamma)>1-\epsilon/2$).
\end{proposition}

\begin{proof}
Let $K(\omega,x) $ denote either $K(\omega,x) = \eta(\omega,x)\cap \Gamma$ or $K(\omega,x) =\Gamma$.
For $\mu$-a.e.~$(\omega',x')\in K(\omega,x)$ we have  
\[(\pi_{2})_{*} \bigl({\mu^{\sigA}_{(\omega',x')}} \restrict {{K(\omega',x')}}\bigr) \ll (\pi_{2})_{*}\mu^{\sigA}_{(\omega',x')}
.\]
Let $$g_{(\omega',x')}(y) = \frac{d (\pi_{2})_{*} \bigl({\mu^{\sigA}_{(\omega',x')}} \restrict {{K(\omega',x') }}\bigr)}{d(\pi_{2})_{*}\mu^{\sigA}_{(\omega',x')}}(y)$$
denote  the Radon--Nikodym derivative. 
Note that   $0< g_{(\omega',x')}(y)\le 1$ for $(\pi_{2})_{*} ({\mu^{\sigA}_{(\omega',x')}} \restrict {K(\omega',x') })$-a.e.~$y$. 
Moreover, 
for $(\pi_{2})_{*}\mu^{\sigA}_{(\omega',x')}$-a.e.~$y\in M$
\[\lim_{r\to 0} \frac{\mu^{\sigA}_{(\omega',x')} (\pi_{2}^{-1}(B(y,r))\cap K)}{\mu_{(\omega',x')}^{\sigA}(\pi_{2}^{-1}(B(y,r)))}=g_{(\omega',x')}(y).\]
The first conclusion follows by choosing an appropriate $C, r$, and  $\Gamma'\subset \Gamma$.  
For the second conclusion,  if $\eta \prec \Xi(\sigA)$ then by Markov's inequality,
\[ \mu\bigl(\{(\omega,x)\in \Gamma: g_{(\omega,x)}(x)\ge 1-\epsilon^{-1}\bigl(1-\mu(\Gamma)\bigr)\}\bigr) \ge \mu(\Gamma)- \epsilon. \qedhere \]
\end{proof}

Using that upper dimension is non-decreasing under averaging, we have the following as a special case of  
\cite[Lem.~12.1.2, p.~573]{MR819557}.  
\begin{proposition}\label{prop:integral} 
Let  $\mu$ be a probability measure on $X=\Omega\times M$ and let $\sigA$ and $ \sigC$ be sub-$\sigma$-algebras with $\sigC\subset \sigA$.  
Suppose  \[\limsup_{r\to0} \frac{\log (\mu_{(\omega,x)}^\sigC (\pi_2\inv B(x,r)))}{\log r} \le \Delta(\omega,x) \]  for some measurable function $\Delta\colon X\to [0,\infty)$ and $\mu$-a.e.~$(\omega,x)\in X$.   Then for $\mu$-a.e.~$(\omega,x)\in X$,
  \[\limsup_{r\to0} \frac{\log (\mu_{(\omega,x)}^\sigA (\pi_2\inv B(x,r)))}{\log r} \le \Delta(\omega,x). \] 
\end{proposition}

\section{Preliminaries: dynamics in Lyapunov charts and adapted partitions}\label{sec:lyapunov}
\label{sec:4}

Throughout this section, we retain the hypotheses and notations from \cref{sec:general}.  
In particular, $M$ is a closed smooth manifold and $\mu$ is an ergodic,  hyperbolic   $F$-invariant measure on $X=\Omega\times M$.  We assume the integrability  \eqref{eq:L1}  hypotheses of \cref{sec:standing ass skew} for the marginal $\nu= (\pi_1)_*\mu$.

Let $d_M$ denote a  fixed distance induced by a  choice of Riemannian metric on $M$.   Let $d=\dim M$.  On an invariant conull set the derivative cocycle has an Oseledec splitting $TM=\bigoplus_{i=1}^{p}E_{\omega}^i (x)$ with Lyapunov exponents $$\lambda_{1}<\dots <\lambda_{s} <0<\lambda_{s+1}<\dots < \lambda_{p}.$$  Recall that we  always assume that $1\le s<p$ so that the fiberwise derivative cocycle has an exponent of each sign.  

 Let \[\Rbb^{d}=\bigoplus_{i=1}^{p}\Rbb_i,\qquad \dim E_{\omega}^i(x)=\dim \Rbb_i,\qquad ||v||=\max_i ||v_i||\] where $v\in \Rbb^d$ is uniquely written $v=\sum_{i=1}^{p} v_i$ with $v_i\in \Rbb_i$ and the summand has its Euclidean norm. Write $$V_s =\bigoplus_{i\le s}  \Rbb_i,\quad \text{and} \quad  V_u =\bigoplus_{i>s} \Rbb_i.$$  Let $$\lambda_{\max}=\max_{i}|\lambda_i|, \quad \text{and} \quad \lambda_{\min}=\min_{i}|\lambda_i |.$$ Fix once and for all $$0<\epsilon_0 <1/100$$ sufficiently small (so that, in particular, 
  $\lambda_{\min} - 100 \epsilon_0 >0$).

\subsection{Two-sided Lyapunov charts}
For a normed vector space $W$, let $W(r)$ denote the open ball of radius $r$ centered at $0$ in $W$.
Recall our fixed $\epsilon_0>0$.  

\begin{proposition}[Two-sided Lyapunov charts]\label{Lyapcharts} There are a constant  $0<\hat \delta<1$, 
   a  $F$-invariant conull set $\Lambda_0\subset X$, a measurable function $\ell\colon \Lambda_0\to [1,\infty)$, $C^\infty $-charts $\Phi_{(\omega,x)}\colon \R^d(1)\to M$, and $C^{1+\beta}$ diffeomorphisms $G_{(\omega,x)}\colon \R^d\to \R^d$   for $(\omega,x)\in  \Lambda_0$ such that the following properties hold for all $(\omega,x)\in \Lambda_0$:
\begin{enumerate}\item (Properties of charts)
\begin{enumerate}
\item $\Phi_{(\omega,x)}\colon \R^d(1)\to M$ is a diffeomorphism onto its image with $\Phi_{(\omega,x)}(0) = x$; 
\item $\|\Phi_{(\omega,x)}\|_{C^1}\le 1$ and $\|\Phi_{(\omega,x)}\inv \|_{C^1}\le  \ell(\omega,x) $, that is, \[d(\Phi_{(\omega,x)}v,\Phi_{(\omega,x)}w)\le \|v-w\| \le \ell(\omega,x)d(\Phi_{(\omega,x)}v,\Phi_{(\omega,x)}w);\]  
\item $D_0\Phi_{(\omega,x)}\R_i = E^i_{\omega}( x)$;
\item $\ell(F^n(\omega,x))\le   e^{|n|\epsilon_0}\ell(\omega,x)$ for $n\in \Z$.  
\end{enumerate}
\item (Properties of globalized dynamics in charts) 
\begin{enumerate}
\item $G_{ (\omega,x)}(0) = 0$.
\item For $v\in \R^d$ with $\|v\|<\hat \delta$, $$G_{(\omega,x)}(v) = \Phi_{F(\omega,x)}\inv \circ f_{\omega}\circ \Phi_{(\omega,x)}(v).$$
and $$G_{ F^{-1}(\omega,x)}\inv(v)  = \Phi_{F\inv(\omega,x)}\inv \circ f_{ \theta^{-1}(\omega)}\inv \circ \Phi_{(\omega,x)}(v).$$
\item $D_0 G$ preserves each $\Rbb_i$. Moreover, for $v\in \R_i$, $$e^{\lambda_i- \epsilon_0}\|v\|\le \| D_0G_{ (\omega,x)}v\|\le e^{\lambda_i+ \epsilon_0}\|v\|.$$
\item $\|G_{(\omega,x)}- D_0G_{ (\omega,x)}\|_{C^1}\le \epsilon_0$ and  $\|G_{ (\omega,x)}\inv - D_0G_{ (\omega,x)}\inv \|_{C^1}\le \epsilon_0$.
\item $\Hol(DG_{(\omega,x)})<C$ and $ \Hol(DG_{(\omega,x)}\inv)\le C$ for some uniform $C$ independent of $\omega$ and $x$.    
\end{enumerate}
\item (Measurability) $(\omega,x)\mapsto \Phi_{(\omega,x)}$ and  $(\omega,x)\mapsto G_{(\omega,x)}$ are measurable.  
\end{enumerate}
\end{proposition}

We refer to the family $\{\Phi_{(\omega,x)}\}$ in \cref{Lyapcharts} as a family of  (two-sided) $(\epsilon_0, \ell)$-Lyapunov charts. We fix this family of charts for the remainder of the paper.

\subsection{Stable and unstable manifolds relative to Lyapunov charts}

Write 
\[G_{(\omega,x)}^{(n)}=\begin{cases} G_{F^{n-1}(\omega,x)}\circ \dots \circ G_{F(\omega,x)}\circ G_{(\omega,x)} & n>0 \\  G_{(\omega,x)}^{(0)}=\textrm{Id},  & n=0 \\ G_{F^{n}(\omega,x)}^{-1}\circ \dots \circ G_{F^{-1}(\omega,x)}^{-1}& n<0 \end{cases}.\] 

Recall the stable and unstable sets defined in \cref{sss:randomstable}.  Using the dynamics in \cref{Lyapcharts}, we justify the conclusion of \cref{lem:pesin}.

\begin{lemma}[Parametrization of stable manifolds]\label{lem:localstable}
For every  $(\omega,x)\in \Lambda_0$,  there is a  $C^{1+\text{H\"older}}$ function $$h^s_{(\omega,x)} \colon V_{s}  \to V_{u}$$  
with the following properties: 

\begin{enumerate}
\item $h^s(0) = 0$ and $D_0h^s = 0$.  
\item There is a constant $c>0$ (independent  of $(\omega,x)$) such that  \[\Hol(Dh^s_{(\omega,x)}) <c.\]
\item $\|Dh^s_{(\omega,x)}\|\le 1/4$ and $h^s_{(\omega,x)}$ is bounded.\end{enumerate}
Let $W_{(\omega,x)}^s$ denote the graph of $h^s_{(\omega,x)} \colon V_{s}  \to V_{u}$.  Then 
\begin{enumerate}[resume]
\item $G_{(\omega,x)}(W_{(\omega,x)}^s )=W_{F(\omega,x)}^s$.  
\item For $u,v\in  W^s_{(\omega,x)}$ 
	$$e^{\lambda_1 - 2 \epsilon_0} \|v-u\| \le  \|G_{ (\omega,x)} (v) -G_{ (\omega,x)}(u)\| \le e^{\lambda_s+ 2\epsilon_0} \|v-u\|.$$
\item $  W^s_{(\omega,x)}$ coincides with the set of all $u\in \R^d$ for which  
$$ \limsup_{n\to \infty}\frac 1 n \log \|G_{  (\omega,x)}^{(n)}(u)\|< 0$$ 
\end{enumerate}
Moreover, 
the family $\{h^s_{(\omega,x)}\}$ depends 
measurably on $(\omega,x)\in   X$.  
\end{lemma}

For $(\omega,x)\in \Lambda_0$ and $0<\delta\le \hat \delta$, write  $$W^s_{(\omega,x), \delta}:= \{u\in  W^s_{(\omega,x)}: \|u\|< \delta\}$$ and
\begin{equation*}
W^s_{ \loc,\delta}(\omega,x):=  \Phi_{ (\omega,x)}( W^s_{(\omega,x),  \delta}).\end{equation*}
We note that the choice of local stable manifold $W^s_{ \loc,\delta}(\omega, x)$ depends on both the future and past itinerary  $\omega$ as it is defined relative to  two-sided charts.  Also note that $W^s_{ \loc,\delta}(\omega,x)$ is an embedded disk in $M$.  

Given the above, we can justify \cref{lem:pesin} as follows: 
For $\mu$-a.e.\ $(\omega,x)$ (whose forward orbit is recurrent to a set on which the function $\ell\colon\Lambda_0\subset X\to (0,\infty)$ in \cref{Lyapcharts} is bounded from above) and every $y\in W^s_\omega(x)$, there is $n\ge 0$ such that 
$$f_\omega^{(n)} (y)\in W^s_{ \loc,\delta}(F^n(\omega,x)).$$
Then $W^s_{\omega}(x)$ becomes a nested union 
$$W^s_{\omega}(x)=\bigcup_{n=0}^\infty \left(f_\omega^{(n)}\right)^{\inv} \left( W^s_{ \loc,\delta}(F^n(\omega,x))\right).$$

We  similarly parameterize unstable manifolds relative to Lyapunov charts via graphs of functions $$h^u_{(\omega,x)} \colon V_{u}  \to V_{s}$$ and similarly write $W^u_{(\omega,x)}$ for the component of the graph of $h^u_{(\omega,x)}$,   $W^u_{(\omega,x), \delta}:= \{u\in  W^u_{(\omega,x)}: \|u\|\mathrel{<} \delta\}$ and
\begin{equation*}
W^u_{ \loc,\delta}(\omega,x):=  \Phi_{ (\omega,x)}( W^{u}_{(\omega,x),  \delta}).\end{equation*}

\subsection{Properties of dynamics relative to Lyapunov charts}

Let
\[
\delta_0=\frac{\hat\delta}{4e^{\lambda_{\max}+2\epsilon_0}}.
\]

We first record Lipschitz bounds on dynamics relative to charts.  

\begin{lemma}\label{lem:estiter}
Let $v\in \R^d$ and write $v=v^s+ v^u$.  
\begin{enumerate}
\item  If $\|v^u\|\ge \|v^s\|$, then $\|(G_{(\omega,x)}v)^u\|\ge \|(G_{(\omega,x)}v)^s\|$ and 
$$e^{\lambda_{s+1}-2\epsilon_0}\|v\| \le \|(G_{(\omega,x)}v)\|\le e^{\lambda_p+2\epsilon_0}\|v\|.$$
\item If $\|v^s\|\ge \|v^u\|$, then $\|(G_{ F^{-1}(\omega,x)}\inv v)^s\|\ge \|(G_{ F^{-1}(\omega,x)}\inv v)^u\|$ and 
$$e^{-\lambda_s-2\epsilon_0}\|v\|\le \|(G_{ F^{-1}(\omega,x)}\inv v)\|\le e^{-\lambda_1+2\epsilon_0}\|v\|.$$
\end{enumerate}
In particular, if $0<\|v\|\le \delta_0$ then 
\begin{enumerate}[resume]
\item $e^{-\lambda_{\max}-2\epsilon_0}\|v\|\le 
	\left\|\Phi_{F(\omega,x)}^{-1}\left(
	f_\omega \bigl( \Phi_{(\omega,x)} (v)\bigr)
	\right) \right\| \le e^{\lambda_{\max}+2\epsilon_0}\|v\|,$ and

\item $e^{-\lambda_{\max}-2\epsilon_0}\|v\|\le 
	\left\|\Phi_{F\inv(\omega,x)}^{-1}\left(
	f_\omega^{(-1)} \bigl( \Phi_{(\omega,x)} (v)\bigr)
	\right) \right\| \le e^{\lambda_{\max}+2\epsilon_0}\|v\|.$

	\end{enumerate}

\end{lemma}

\begin{corollary}\label{cor:stayinsmall}
	Fix  $\delta\le \delta_0$,  $(\omega,x)\in \Lambda_0$, and  $v\in\Rbb^{d}$.
	\begin{enumerate}
	\item If 
	\(||G^{(k)}_{(\omega,x)}v|| < \delta\) for all $-n\le k\le n$ then 
	\[v\in \R^d\bigl(\delta\exp\left(-n(\lambda_{\min}-3\epsilon_0)\right){\bigr)}.\]  
	\item If 
	\(||G^{(k)}_{(\omega,x)}v|| < \delta\) for all $0\le k\le n$ then 
	\[v\in \R^d\bigl(\delta\exp\left(-n(\lambda_{\min}-3\epsilon_0)\right)\bigr)\cup C^s(\delta).\]  
	\item If 
	\(||G^{(k)}_{(\omega,x)}v|| < \delta\) for all $-n\le k\le 0$ then 
	\[v\in \R^d\bigl(\delta\exp\left(-n(\lambda_{\min}-3\epsilon_0)\right){\bigr)}\cup C^u(\delta).\]  

\end{enumerate}

\end{corollary}

\begin{lemma}\label{lem:Scucs}
\label{lem:dynamicsoverflow}

For every $0<\delta\le \delta_0$, and every $(\omega,x)\in \Lambda_0$ the following hold: 
\begin{enumerate}
\item we have 
	\begin{align*}W^{s}_{(\omega,x),\delta}&=\left\{v\in \Rbb^{d}(\delta): \left\|G^{(n)}_{(\omega,x)}(v)\right\|<\delta \quad\textrm{ for all } n\ge 0 \right\}\\&=\left\{v\in \Rbb^{d}(\delta) : \left\|\Phi_{F^n (\omega,x)}^{-1} \circ f_{\omega}^{(n)} \circ \Phi_{(\omega,x)}(v)\right\|< \delta  \textrm{   for all   } n\ge 0\right\}\end{align*} and
\begin{align*}W^{u}_{(\omega,x),\delta}&=\left\{v\in\Rbb^{d}(\delta): \left\|G^{(n)}_{(\omega,x)}(v)\right\|<\delta\quad\textrm{ for all } n\le 0\right\}\\ 
&=\left\{v\in \Rbb^{d}(\delta): \left\|\Phi_{F^{n}(\omega,x)}^{-1}\circ f_{\omega}^{(n)} \circ \Phi_{(\omega,x)}(v)\right\|< \delta \textrm{ for all } n\le 0\right\}
\end{align*}

\end{enumerate}
In particular,
\begin{enumerate}[resume]
\item 
$W^{s}_{(\omega,x),\delta}=G^{-1}_{(\omega,x)}\left(W^{s}_{F(\omega,x),\delta}\right)\cap \Rbb^{d}(\delta)$

\item 
$W^{u}_{(\omega,x),\delta}=G_{{F^{-1}}(\omega,x)}\left(W^{u}_{F^{-1}(\omega,x),\delta} \right)\cap \Rbb^{d}(\delta)$

\end{enumerate}

\end{lemma}

\subsection{Properties of partitions adapted to Lyapunov charts} 

Having constructed a family of two-sided $({\epsilon_0},\ell)$-Lyapunov charts and established properties of the dynamics in charts and local stable and unstable manifolds relative to those charts, we construct a partition $\partP$ of $(X,\mu)$ whose properties will be used for the remainder of the paper.

For the statement and constructions, we also recall the choice of subordinate partitions defined in \cref{lem:LS} with properties in \cref{subsec:subord}.

\begin{proposition}\label{prop:outcomeofadapted}
For every $0<\delta\le \delta_0$, there exists a measurable partition $\partP$ of $(X,\mu)$ with $H_\mu(\partP)<\infty$ and a measurable function $n_0\colon X\to\N$ such that for a.e.~$(\omega,x)\in X$ and for $n\ge n_0(\omega,x)$, the following properties hold:
\begin{enumerate}[label= {(\arabic*)}, ref= {(\arabic*)}]
\item \label{adaptedok0} 
 For every $(\omega',y)\in(\Xi(\sigC)\vee \partP)^{[-n,n]}(\omega,x)$,
 $$f_{\omega'}^{(j)} (y)\in \Phi_{F^j(\omega,x)}(\R^d(\delta))\subset B_M(f_\omega^{(j)}(x),\delta)
 $$ for all $-n\le j\le n$.  

\item \label{adaptedok0BD} 
If we assume hypothesis \eqref{eq:BD} of \cref{sec:standing ass skew}  then for every $(\omega',y)\in \partP^{[-n,n]}(\omega,x)$,
 $$f_{\omega'}^{(j)} (y)\in \Phi_{F^j(\omega,x)}(\R^d(\delta))\subset B_M(f_\omega^{(j)}(x),\delta)
 $$ for all $-n\le j\le n$.

\item \label{adaptedok1} 
 For every $(\omega',y)\in(\Xi(\sigC)\vee \partP)^{[-n,0]}(\omega,x)\cap \xi^s(\omega,x)$,
 $$f_{\omega'}^{(j)} (y)\in W^{s}_{\loc,\delta}(\theta^{(j)}(\omega),f_\omega^{(j)}(x))$$ for all $-n\le j\le 0$.  

\item   \label{adaptedok2} If we assume hypothesis \eqref{eq:BD} of \cref{sec:standing ass skew} then for every $(\omega',y)\in \partP^{[-n,0]}(\omega,x)\cap \xi^s(\omega,x)$,
 $$f_{\omega'}^{(j)} (y)\in W^{s}_{\loc,\delta}(\theta^{(j)}(\omega),f_\omega^{(j)}(x))$$ for all $-n \le j\le 0$.  

\item \label{adaptedok5}  For any $(\omega',y)\in \partP^{[0,n]}(\omega,x)\cap \overline \xi^u(\omega,x)$, 
$$f_{\omega'}^{(j)} (y)\in W^{u}_{\loc,\delta}(\theta^{(j)}(\omega),f_\omega^{(j)}(x))$$ for all $0\le j\le n$.  
\item \label{adaptedok6}  If we assume hypothesis \eqref{eq:BD} of \cref{sec:standing ass skew} then for  any $(\omega',y)\in \partP^{[0,n]}(\omega,x)\cap  \xi^u(\omega,x)$, 
 $$f_{\omega'}^{(j)} (y)\in W^{u}_{\loc,\delta}(\theta^{(j)}(\omega),f_\omega^{(j)}(x))$$ for all $0\le j\le n$.  

\end{enumerate}
\end{proposition}

We will rigorously justify conclusions \ref{adaptedok0}--\ref{adaptedok2} of \cref{prop:outcomeofadapted} involving dynamics in ambient Lyapunov charts or along local stable manifolds relative to Lyapunov charts.   Conclusions \ref{adaptedok5} and \ref{adaptedok6} involving dynamics along local unstable manifolds relative to Lyapunov charts can be obtained by considering inverses.   Note that given a partition $\partP$ satisfying any conclusion of  \cref{prop:outcomeofadapted} and a partition $\partQ$ with $H_\mu(\partQ)<\infty$, the partition $\partP\vee \partQ$ still satisfies the same conclusion of \cref{prop:outcomeofadapted}.

\subsubsection{A finite entropy partition with dynamically adapted diameter} 
Before proving  \cref{prop:outcomeofadapted}, we first construct a  partition $\partR$ that will serve as the prototype for the partition $\partP$ 
in the statement of \cref{prop:outcomeofadapted}.  

Recall the local stable manifolds $W^s_{\loc,\delta}$ constructed in \cref{lem:localstable} relative to the Lyapunov charts whose properties are outlined in \cref{Lyapcharts}.
 Recall that $\xi^s=\tau\inv \hat \xi^s$.   Recall that for a.e.~$(\omega,x)\in \what X$,  the atom $\hat \xi^s(\omega,x)$ is precompact in $W^{s}(\omega,x)\subset \{\omega\}\times M$.

\begin{claim}\label{claim:localslice}
There is a measurable function $c^s\colon (X,\mu)\to (0,1]$ such that  
$$\pi_2\bigl(\bigl\{ (\omega',y)\in \xi^s(\omega, x) : d_M(x,y)\le c^s(\omega,x)\bigr\}\bigr) \subset  W^s_{\loc,\delta_0}(\omega,x).$$

Similarly, there is a measurable function $c^u\colon (X,\mu)\to (0,1]$ such that 
$$\pi_2\bigl(\bigl\{ (\omega',y)\in \overline \xi^u(\omega, x) : d_M(x,y)\le c^u(\omega,x)\bigr\}\bigr) \subset  W^u_{\loc,\delta_0}(\omega,x).$$
If we further assume the hypotheses \eqref{eq:BD}  
of \cref{sec:standing ass skew}, then $$\pi_2\bigl(\bigl\{ (\omega',y)\in  \xi^u(\omega, x) : d_M(x,y)\le c^u(\omega,x)\bigr\}\bigr) \subset  W^u_{\loc,\delta_0}(\omega,x).$$

\end{claim}

Fix a set $\Lambda\subset X$ with $\mu(\Lambda)>0$ on which the measurable function  $\ell\colon X\to (0,\infty)$ is  bounded from above by $1\le \ell_0<\infty$ and the measurable functions $c^s,c^u \colon X\to (0,1]$ are bounded from below by $0<c_0\le 1$.

For $(\omega,x)\in X$, let
\[\begin{aligned}
T_+(\omega,x)&=\inf\{m\ge 1:F^m(\omega,x)\in\Lambda\},\\
T_-(\omega,x)&=\inf\{m\ge1:F^{-m}(\omega,x)\in\Lambda\},\\
T(\omega,x)&=\max\{T_+(\omega,x),T_-(\omega,x)\},\\
n_+(\omega,x)&=\inf\{m\ge 0:F^m(\omega,x)\in\Lambda\},\\
n_-(\omega,x)&=\inf\{m\ge0:F^{-m}(\omega,x)\in\Lambda\}.
\end{aligned}\]

Set $q_0= \ell_0\inv \min\{c_0,\delta\}$ and define
\[r(\omega,x)=
\begin{cases}
q_0  e^{-T(\omega,x)(\lambda_{\max}+3{\epsilon_0})}&(\omega,x)\in\Lambda,\\
1 &(\omega,x)\notin\Lambda.
\end{cases}\]
Then $-\log r$ is in $L^1(\mu)$.  
By a result of Ma\~n\'e \cite[Lem.\  2]{MR627789} (see also \cite[Clm.~5.13]{BRH} for minor modifications to the more general setting), there is a measurable partition $\partR$ of $X$ with $H_\mu(\partR)<\infty$ and such that $$\diam (\pi_2(\partR(\omega,x))) \le  \frac 1 2 r(\omega,x)$$ 
for $\mu$-a.e.~$(\omega,x)$.

\subsubsection{Properties of the partition $\partR$}

The partition $\partR$ has the following dynamical property.
\begin{lemma}\label{Claim:controldelta}
For a.e.~$(\omega,x)\in X$ and for any $n, m\ge 0 $ such that either 
 $n\ge n_-(\omega,x)$ or $m\ge n_+(\omega,x)$ the following holds:
if $(\omega',y)\in \partR^{[-n,m]}(\omega,x)$ and if 
\begin{equation}\label{samedyn}f_{\theta ^j(\omega')} = f_{\theta ^j(\omega)},\quad \text{
for all $-n\le j<m$,} \end{equation}
then  for all $j\in \Z$ with $-n -T_-(F^{-n}(\omega,x))\le j\le m+ T_+(F^{m}(\omega,x))$,
\begin{equation}\label{dypart} f_\omega^{(j)}(y)\in \Phi_{F^j(\omega,x)}(\R^d(\ell_0 q_0)).\end{equation}
\end{lemma}

\begin{proof}
We take $(\omega,x)\in X$ that returns infinitely often to $\Lambda$ in both the future and in the past. Let  $$\cdots<r_{-1}<0\le r_0<r_1<\cdots$$ denote the distinct times for which $F^{r_j}(\omega,x)\in\Lambda$.
By the assumption, either  $-n\le r_{-1}$ or $r_0\le m$ 
and thus there are $p\le q$ for which 
$${-n-T_-(F^{-n}(\omega,x))}=r_{p-1}< -n \le r_p\le r_{q}\le  m< r_{q+1}={m+T_+(F^{m}(\omega,x))}.$$

To prove the lemma, consider any $j\in \Z$ with 
$$r_{p-1}\le j\le  r_{q+1}.$$ 
Take $p-1\le k \le q$ with $r_{p-1}\le r_k\le j\le r_{k+1}\le r_{q+1}$.

Note that either $$\text{$-n\le r_{k+1}\le m$ \quad  or \quad  $-n\le r_k\le m$}.$$ We assume the latter case (that $-n\le r_k\le m$) and show the conclusion \eqref{dypart} holds for all $r_k\le j\le r_{k+1}$.  
Since we assume $-n\le r_{k}\le m$ and $(\omega',y)\in \partR^{[-n,m]}(\omega,x)$,
we have 
$$F^{r_k}(\omega',y)\in \partR (F^{r_k}(\omega,x))$$
and from hypothesis  we have
$$f_{\omega}^{(r_{k})} (y)= f_{\omega'}^{(r_{k})} (y).$$
In particular, 
$$f_{\omega}^{(r_{k})} (y)= f_{\omega'}^{(r_{k})} (y)\in B(f_{\omega}^{(r_{k})} (x), r(F^{r_k}(\omega,x))).$$
Since $F^{r_k}(\omega,x)\in \Lambda$, we have 
$$ \Phi_{F^{r_k}(\omega,x)}\inv (f_{\omega}^{(r_{k})} (y))\in \R(\ell_0 q_0  e^{-T(F^{r_{k}}(\omega,x))(\lambda_{\max}+3{\epsilon_0})})\subset
\R(\ell_0 q_0  e^{-(r_{k+1}-r_{k})(\lambda_{\max}+3{\epsilon_0})}).
$$
From \cref{lem:estiter}, for every $r_k\le j\le  r_{k+1}$ we have 
$$ \Phi_{F^{j}(\omega,x)}\inv (f_{\omega}^{(j)} (y))\in 
\R(\ell_0  q_0  e^{-(r_{k+1}-r_{k})(\lambda_{\max}+3{\epsilon_0})} e^{(j-r_k)(\lambda_{\max}+2{\epsilon_0}) })\subset \R(\ell_0  q_0 )
$$
and the conclusion follows.

Similarly, if $-n\le r_{k+1}\le m$ then $$F^{r_{k+1}}(\omega',y)\in \partR (F^{r_{k+1}}(\omega,x))$$
and from hypothesis  we have
$$f_{\omega}^{(r_{{k+1}})} (y)= f_{\omega'}^{(r_{{k+1}})} (y).$$
A similar argument as above shows \eqref{dypart} holds for all $r_k\le j\le r_{k+1}$ 
\end{proof}

\subsubsection{Modifications to $\partR$} 
We modify $\partR$ if needed to obtain a partition $\partP$ where the hypothesis \eqref{samedyn} of \cref{Claim:controldelta} can be verified.  
Under the hypotheses in  \ref{adaptedok0}, \ref{adaptedok1}, or \ref{adaptedok5},  we let $\partP=\partR$.  
Then for $(\omega',y)\in (\Xi(\sigC)\vee \partP)^{[-n,m]}(\omega,x)$, we have 
$$f_{\theta^{j}(\omega)} =f_{\theta^{j}(\omega')} $$ for all $-n\le j< \infty$. 
In the case of hypotheses \ref{adaptedok5}, note that if $(\omega',y)\in \partP^{[0,n]}(\omega,x)\cap \overline \xi^u(\omega,x)$, then $\omega'=\omega$ by definition. 
 Similarly, under the hypotheses in  \ref{adaptedok0BD},  \ref{adaptedok2}, or \ref{adaptedok6}, we let $\partP=\partR\vee \wtd\partQ$. Then again, if   
$(\omega',y)\in \partP^{[-n,m]}(\omega,x)$, we have $$f_{\theta^{j}(\omega)} =f_{\theta^{j}(\omega')} $$ for all $-n\le j<m$.

\subsubsection{Proof  of \cref{prop:outcomeofadapted}}
We now   justify conclusions \ref{adaptedok0}, \ref{adaptedok0BD}, \ref{adaptedok1} and \ref{adaptedok2}.  Conclusions    \ref{adaptedok5} and \ref{adaptedok6}  are similar by passing to inverses.

\begin{proof}[Proof  of \cref{prop:outcomeofadapted}]
Conclusions \ref{adaptedok0} and \ref{adaptedok0BD} hold for a.e.~$(\omega,x)$ and for all  $n\ge  \max \{n_-(\omega,x),n_+(\omega,x)\}$ by  \cref{Claim:controldelta} using that $$f_{\omega}^{(j)} (y)= f_{\omega'}^{(j)} (y)$$
for all $-n\le j\le n$.

We now show conclusions \ref{adaptedok1} and  \ref{adaptedok2}. 

We first claim that if  $n\ge  n_-(\omega,x)$, if 
$f_{\theta ^j(\omega')} = f_{\theta ^j(\omega)}$ for all $-n\le j<0$,
and if  $(\omega',y)\in \partP^{[-n,0]}(\omega,x)\cap \xi^s(\omega,x)$ then   $y\in W^{s}_{\loc,\delta}(\omega,x)$.  
Indeed,  \cref{Claim:controldelta} implies for every $j\in \Z$ with $-n_-(\omega,x)\le j\le n_+(\omega,x)$ that
\begin{equation*}\label{eq:controllyapdyn}f_\omega^{(j)}(y)\in \Phi_{F^{j}(\omega,x)}(\R^d(\ell_0q_0))\subset  \Phi_{F^{j}(\omega,x)}(\R^d(\delta)).\end{equation*}  
Since Lyapunov charts are 1-Lipschitz,
\begin{equation}\label{eq:controllyapdynedn}
d_M(f_\omega^{(n_+(\omega,x))}(y),f_\omega^{(n_+(\omega,x))}(x))\le \ell_0q_0\le c_0.
\end{equation}
Since $\xi^s\prec F\xi^s$, if $(\omega',y)\in \xi^s(\omega,x)$ then $(\omega,y)\in \xi^s(\omega,x)$ and so
we have $$(\theta^{n_+(\omega,x)}(\omega), f_\omega^{({n_+(\omega,x)})}(y) )\in \xi^s(F^{n_+(\omega,x)}(\omega,x) ).$$
By  \eqref{eq:controllyapdynedn} and definition of $c_0$, we have $$
f_\omega^{({n_+(\omega,x)})}(y) \in W^{s}_{\loc,\delta}(F^{n_+(\omega,x)}(\omega, x)).$$
Repeated application of \cref{lem:dynamicsoverflow} implies $y\in W^{s}_{\loc,\delta}(\omega,x)$.  

Under the hypotheses in  \ref{adaptedok1},
for every $n\ge  n_-(\omega,x)$, every $-n\le j\le 0$, and every $(\omega',y)\in(\Xi(\sigC)\vee \partP)^{[-n,0]}(\omega,x)\cap \xi^s(\omega,x)$,
we have $f_{\omega'}^{(j)} (y) =f_{\omega}^{(j)} (y) \in \Phi_{F^{j}(\omega,x)}(\R^d(\delta))$.  
Since $y\in W^{s}_{\loc,\delta}(\omega,x)$, repeated applications of \cref{lem:dynamicsoverflow}  imply  $$f_{\omega'}^{(j)} (y)=f_{\omega}^{(j)} (y)\in W^{s}_{\loc,\delta}(\theta^{(j)}(\omega),f_\omega^{(j)}(x))$$ for all $-n\le j\le 0$.  
 Under the hypotheses in  \ref{adaptedok2}, the same holds if   
$(\omega',y)\in \partP^{[-n,0]}(\omega,x)\cap \xi^s(\omega,x)$.  \end{proof}

Write $C^u(\delta) = \{v\in \R^d(\delta) : \|v^u\|\ge \|v^s\|\}$ and $C^s(\delta) = \{v\in \R^d(\delta) : \|v^s\|\ge \|v^u\|\}$.

\subsubsection{Dynamical properties of $\partP$}
Fix a positive integer $a$ such that $a(\lambda_{\min}-3\epsilon_0)>1$.
Using  \cref{cor:stayinsmall}  and \cref{Claim:controldelta}, we obtain the following:
\begin{corollary}\label{lem:inball}  Fix $0<\delta\le \delta_0$.  
Let $\partP$ be as in \cref{prop:outcomeofadapted}.   Then there exists a measurable map $\overline{n}\colon X \to \Zbb_{+}$
such that for a.e.~$(\omega,x)$ and $n\ge \overline{n}(\omega,x)$, the following hold:
\begin{enumerate}[label= {(\arabic*)}, ref= {(\arabic*)}]
\item\label{itboth} $\pi_{2} \left(\left(\Xi(\sigC)\vee\Pcal\right)^{[-{an},{an}]}(\omega,x)\right)\subset \Phi_{(\omega,x) }\left( \R^d({e^{-n}})\right)\subset  B(x, {e^{-n}})$.

\item\label{itsn} $\pi_2 \left((\Xi(\sigC)\vee\Pcal)^{[-{an},0]}(\omega,x)\right)\subset 
\Phi_{(\omega,x) }\left( \R^d({e^{-n}})\cup C^u(\delta)\right)$

\item\label{itun}  $\pi_{2}\left((\Xi(\sigC)\vee \Pcal)^{[0,{an}]}(\omega,x)\right)\subset \Phi_{(\omega,x) }\left( \R^d({e^{-n}})\cup C^s(\delta)\right)$

\item \label{itboth2} If we assume hypothesis \eqref{eq:BD} of \cref{sec:standing ass skew},  then  $$\pi_{2} \left(\Pcal^{[-{an},{an}]}(\omega,x)\right)\subset \Phi_{(\omega,x) }\left( \R^d({e^{-n}})\right)\subset  B(x, {e^{-n}}).$$ 

\item If we assume hypothesis \eqref{eq:BD} of \cref{sec:standing ass skew}, then 
$$\pi_{2}\left( \Pcal^{[0,{an}]}(\omega,x)\right)\subset \Phi_{(\omega,x) }\left( \R^d({e^{-n}})\cup C^s(\delta)\right).$$
and $$\pi_{2}\left( \Pcal^{[-{an},0]}(\omega,x)\right)\subset \Phi_{(\omega,x) }\left( \R^d({e^{-n}})\cup C^u(\delta)\right).$$

\end{enumerate}

\end{corollary}

\subsection{Partition atoms in general position relative to other charts}\label{sss:prebesicovitch}

Recall we view the family of Lyapunov charts $(\omega,x)\mapsto \Phi_{(\omega,x)}$ as a measurable function $X\to \Emb^1(\R^d(1), M)$.  
Given $(\omega,x), (\omega'',z)\in \Lambda_0$, with $\Phi_{(\omega'',z)} (\R^d(1/2))\subset \Phi_{(\omega,x)} (\R^d( 1 ))$, let 
$$ \overline \Delta\bigl((\omega'',z); (\omega,x)\bigr)= \| \Phi_{(\omega,x)}\inv \circ  {\Phi_{(\omega'',z)}}\restrict{ \R^d(\frac 1 2)}-\Id \|_{C^1}.$$
Otherwise, define $\overline \Delta\bigl((\omega'',z); (\omega,x)\bigr)=+\infty$.  
Recall $\delta_0<\frac {\hat \delta}{4}$.   
Let $h^s_{(\omega,x)}$ and $h^u_{(\omega,x)}$ denote their restrictions to $V^s(2\delta_0)$ and $V^u(2\delta_0)$.  
Define  
\begin{align*}\Delta\bigl((\omega,x), (\omega'',z)\bigr)
&:=\overline \Delta\bigl((\omega'',z); (\omega,x)\bigr)+ \overline \Delta\bigl( (\omega,x);(\omega'',z)\bigr)\\
&\quad \quad +
\|h^s_{(\omega,x)}-h^s_{(\omega'',z)}\|_{C^1}+ \|h^u_{(\omega,x)}-h^u_{(\omega'',z)}\|_{C^1}.\end{align*}

From measurability of functions $(\omega,x)\mapsto \Phi_{(\omega,x)}$, $(\omega,x)\mapsto h^s_{(\omega,x)}$, and $(\omega,x)\mapsto h^u_{(\omega,x)}$ (and second countability of the Polish space $X$)
 we obtain the following.  
\begin{claim}\label{countpart}
For every $\epsilon'$, there exists a countable partition $\partK$ of $X$ such that for $(\omega,x), (\omega'',z)\in \Lambda_0$ with 
$ (\omega'',z)\in \partK(\omega,x)$, $\Delta\bigl((\omega,x), (\omega'',z)\bigr)<\epsilon'$. 
\end{claim}

Using that the charts and the parametrized local stable  and unstable manifolds are close on atoms of $\partK$ gives the following. Recall we fix $a\in \Z_+$ with $a(\lambda_{\min}-3{\epsilon_0})>1$.  
\begin{lemma}
\label{lem:chart-separation}
Let $\zeta$ be a sub-invariant measurable partition.  Assume either $\Xi(\sigC)\prec\zeta$ or that hypothesis \eqref{eq:BD} of \cref{sec:standing ass skew} holds.  
There is  $0<\epsilon'<1$ sufficiently small so that the following holds: 

Fix any $0<\delta\le \delta_0$ and partition $\partP$ as in \cref{prop:outcomeofadapted}.
For a.e.~$(\omega'',z)$, all sufficiently large $n$, 
and 
 $(\omega,x)\in  \Lambda_0\cap \partK(\omega'',z)$, if $y\in W^s_{\loc,\delta_0}(\omega,x)$ is such that
  \begin{itemize}
  \item $\|\Phi_{(\omega'',z)}\inv (x)\|\le  e^{-n}$
\item $\|\Phi_{(\omega'',z)}\inv (y)- \Phi_{(\omega'',z)}\inv (x)\|\ge 8e^{-n}$
\end{itemize}
then 
$$B(y,\ell(\omega'',z)\inv e^{-n})\cap \pi_2\bigl((\zeta\vee\partP)^{[-an,0]}(\omega'',z)\bigr)=\emptyset.$$
If instead $y\in W^u_{\loc,\delta_0}(\omega,x)$ then 
 $$B(y,\ell(\omega'',z)\inv e^{-n})\cap \pi_2\bigl((\zeta\vee\partP)^{[0,an]}(\omega'',z)\bigr)=\emptyset.$$
\end{lemma}
\begin{proof}
If $\epsilon'>0$ is sufficiently small, then the set
$${\Phi_{(\omega'',z)}\inv\bigl(W^{s}_{\loc,\hat \delta}(\omega,x)\bigr)\cap \R^d(2\delta_0) }$$coincides  with (the restriction to $\R^d (2\delta_0)$ of)  the graph of a $C^1$ function
$$\td h\colon V^s(1)\to V^u$$ 
with $\|D\td h\|\le \frac 1 3$.

Let $v=\Phi_{(\omega'',z)}\inv (y)$ and $v_0=\Phi_{(\omega'',z)}\inv (x)$.   Then, 
$\|v\|\ge 7e^{-n}$ and 
\begin{align*}
\|v^u \|&\le \|v_0^u\|+  \frac 1 3 \|v^s-v_0^s\|
\le  e^{-n}+ \frac 1 3 \|v^s\| + \frac 1 3 e^{-n}\\
&\le
\frac 43 e^{-n} +  \frac 1 3  \|v\|
\le
\frac 2 7 \|v\|+ \frac 1 3 \|v\|\\
&\le \|v\|
\end{align*}

We have
$$B(y,\ell(\omega'',z)\inv e^{-n})\subset \Phi_{(\omega'',z)}\left( v+ \R^d(e^{-n})\right).$$
If $w\in v+ \R^d(e^{-n})$ then 
$$\|w\|\ge \|v\|- e^{-n}\ge  7e^{-n}- e^{-n} 
=6e^{-n}>2\delta e^{-an(\lambda_{\min}-3\epsilon_0)}.$$
Moreover,
\begin{align*}
\|w^u\|&\le \|v^u\| + e^{-n}
\le \frac 2 7 \|v\|+ \frac 1 3 \|v\| +  e^{-n}\\
&\le \frac 3 7 \|v\|+ \frac 1 3 \|v\|\\
&< \|v\| - \frac 1 7 \|v\|\le \|w^s\|
\end{align*}
and so $w\notin C^u(\delta)$.  
The conclusion then follows from \cref{lem:inball}.   
For $y\in W^u_{\loc,\delta_0}(\omega,x)$, a similar argument holds.
\end{proof}

 \section{Preliminaries:  entropy estimates}\label{sec:entropies}\label{sec:5}
 
 We establish various entropy computations.  Our primary results are the computations in \cref{prop:entcal}.  
 
Recall the various sub-invariant partitions subordinate to fiber stable manifolds and related notation from \cref{sec:notanddefnparti}. Also, we fix a measurable sub-invariant partition $\zeta$ on $X$.  Throughout, we always assume $ \maxent_\mu(F\mid \zeta)<\infty$.

\subsection{Properties of  entropies}
We have a slight variant of \cite[Prop.~5.5]{BRH}. 
 \begin{claim}[{cf.~\cite[Prop.~5.5]{BRH}}]\label{claim:slicingC}
Let  $\zeta$ be  a  sub-invariant measurable partition on $X$.  
  Assume that either 
 \begin{enumerate}
 \item $\Xi(\Csc)\prec\zeta$ and $h_{\mu}(F, \zeta)<\infty$; or 
 \item hypothesis \eqref{eq:BD} of \cref{sec:standing ass skew} holds.
 \end{enumerate} 
 Then for any sub-invariant, measurable partition $\xi^s$ as in 
\cref{def:pullback} of  \cref{sec:submarkov},
 \[h^{\refine}_{\mu}(F\mid \sigma(\zeta))=h^{\refine}_{\mu}(F\mid \sigma(\zeta)\vee \Csc)=h_{\mu}(F,\zeta\vee\xi^{s}).\]

Suppose  $\Xi(\Csc)\prec\zeta$.  Then for any measurable partition $\overline \xi^u$ as in 
\cref{def:subord} of  \cref{sec:submarkov} that is sub-invariant for $F\inv$,
 \[h_{\mu}(F \mid \sigma(\zeta^+))=h_{\mu}(F\inv\mid \sigma(\zeta^+))=h_{\mu}(F\inv,\zeta^+\vee \overline \xi^u).\]
Suppose hypothesis \eqref{eq:BD} of \cref{sec:standing ass skew} holds.   Then for any measurable partition $ \xi^u$ as in \cref{def:pullback} of  \cref{sec:submarkov} that is sub-invariant for $F\inv$,
  \[h_{\mu}(F \mid \sigma(\zeta^+))= h_{\mu}(F\inv \mid \sigma(\zeta^+)) =h_{\mu}(F\inv,\zeta^+\vee\xi^{u}).\]
\end{claim}

We note that by the fiberwise Ruelle inequality and hypothesis \eqref{eq:L1}  of \cref{sec:standing ass skew},  $h^{\refine}_{\mu}(F\mid \sigma(\zeta))<\infty$ under either  hypothesis of \cref{claim:slicingC}.  

\begin{proof}   [Proof of \cref{claim:slicingC}]
In the case that $\Xi(\Csc)\prec\zeta$ and $h_{\mu}(F, \zeta)<\infty$, this is exactly \cite[Prop.~5.5]{BRH}.
In the case hypothesis \eqref{eq:BD} of \cref{sec:standing ass skew} holds, we need only show $$h^{\refine}_{\mu}(F\mid \sigma(\zeta))=h_{\mu}^\refine(F\mid \sigma(\zeta)\vee \Csc).$$
Note $h^{\refine}_{\mu}(F\mid \sigma(\zeta)\vee \Csc)\le h^{\refine}_{\mu}(F\mid \sigma(\zeta))$  holds by definition.
 Recall that $\wtd{\Qcal}$ denotes the lift of the one-sided generator in $\Omega$ and, in particular,  $\wtd{\Qcal}^- = \Xi(\Csc)$. 
  Using that $H_\mu(\wtd \partQ)<\infty$, for any measurable partition $\xi$ on $X$, we have 
\begin{align*}
h_{\mu}^{\refine}(F\mid \sigma(\zeta)\vee\Csc) &\ge h_{\mu}(F,\xi\vee\zeta\vee\Xi(\Csc))\\
&=h_{\mu}(F,\xi\vee\zeta \vee \wtd\Qcal ) \ge h_{\mu}(F,\xi \vee\zeta)
\end{align*}
where the final inequality holds, for instance, from \cite[(8.7)]{MR0217258}.  
Taking the supremum over all $\xi$ shows the reverse inequality.  

The unstable versions follow by passing from $F$ to $F\inv$, using that $\Csc^+ =\sigF$, and using that $\zeta^+$ is sub-invariant for $F\inv$.    
\end{proof}

\subsection{Entropy calculations}
Set \[h= h_{\mu}(F^{\pm1} \mid \sigma(\zeta^+))= h_{\mu}(F\mid \sigma(\zeta))\]
and \[\hat {h}=h_{\mu}^{\refine}(F\mid \sigma(\zeta)). \]

Recall that under our standing hypotheses, $h$ and $\hat h$ may be computed using \cref{claim:slicingC}.  
 After a slight refinement of the partitions in \cref{prop:outcomeofadapted} if needed, we obtain partitions whose iterations capture all entropies considered in the sequel.

\begin{proposition}\label{prop:entcal}
Let $\zeta$ be a sub-invariant measurable partition for $F\colon X\to X$. 
\begin{enumerate}[label= {(\arabic*)}, ref= {(\arabic*)}]
\item \label{entcal1} For every $\epsilon>0$, there exists a measurable partition $\partP$
of $(X,\mu)$ with $H(\partP)<\infty$ such that for $\mu$-a.e.\
$(\omega,x)\in X$,
\begin{equation}\label{eq:1_}
  \hat h + h - \epsilon
  \le \lim_{n\to\infty} -\frac{1}{n}
      \log\mu_{(\omega,x)}^\zeta
      \bigl((\zeta\vee\Pcal)^{(-n,n)}(\omega,x)\bigr)
  \le \hat h + h,
\end{equation}
\begin{equation}\label{eq:2_}
  \hat h - \epsilon
  \le \lim_{n\to\infty} -\frac{1}{n}
      \log\mu_{(\omega,x)}^\zeta
      \bigl((\zeta\vee\Pcal)^{(-n,0]}(\omega,x)\bigr)
  \le \hat h,
\end{equation}
and
\begin{equation}\label{eq:3_}
  h - \epsilon
  \le \lim_{n\to\infty} -\frac{1}{n}
      \log\mu_{(\omega,x)}^\zeta
      \bigl((\zeta\vee\Pcal)^{[0,n)}(\omega,x)\bigr)
  \le h.
\end{equation}

\item \label{entcal2}  If we further assume either $\Xi(\Csc)\prec\zeta$ or \eqref{eq:BD}
then, for every $\epsilon>0$, there exists a measurable partition $\partP$
of $(X,\mu)$ with $H(\partP)<\infty$ such that for $\mu$-a.e.\
$(\omega,x)\in X$,
\begin{equation}\label{eq:4}
 \lim_{n\to\infty} -\frac{1}{n}
         \log\mu_{(\omega,x)}^{\xi^s\vee\zeta}
         \bigl((\zeta\vee\Pcal)^{(-n,0]}(\omega,x)\bigr)= \hat h,
\end{equation}

\begin{equation}\label{eq:s0}
  \lim_{n\to\infty} -\frac{1}{n}
  \log\mu_{(\omega,x)}^{\xi^s\vee\zeta}
  \bigl((\zeta\vee\Pcal)^{[0,n)}(\omega,x)\bigr)= 0,
\end{equation}
and
\begin{equation}\label{eq:sboth}
     \lim_{n\to\infty} -\frac{1}{n}
         \log\mu_{(\omega,x)}^{\xi^s\vee\zeta}
         \bigl((\zeta\vee\Pcal)^{(-n,n)}(\omega,x)\bigr)
     = \hat h.
\end{equation}

\item \label{entcal3}  If we further assume   $\Xi(\Csc)\prec\zeta$ then for every $\epsilon>0$, there exists a measurable partition $\partP$
of $(X,\mu)$ with $H(\partP)<\infty$ such that for $\mu$-a.e.\
$(\omega,x)\in X$,
\begin{equation}\label{eq:ufuture}
    \lim_{n\to\infty} -\frac{1}{n}
         \log\mu_{(\omega,x)}^{\bar\xi^u\vee\zeta^+}
         \bigl((\zeta\vee\Pcal)^{[0,n)}(\omega,x)\bigr)
     = h,
\end{equation}
\begin{equation}\label{eq:u0}
  \lim_{n\to\infty} -\frac{1}{n}
  \log\mu_{(\omega,x)}^{\overline{\xi}^u\vee\zeta^+}
  \bigl((\zeta\vee\Pcal)^{(-n,0]}(\omega,x)\bigr)= 0,
\end{equation}
and
\begin{equation}\label{eq:uboth}
    \lim_{n\to\infty} -\frac{1}{n}
         \log\mu_{(\omega,x)}^{\overline{\xi}^u\vee\zeta^+}
         \bigl((\zeta\vee\Pcal)^{(-n,n)}(\omega,x)\bigr)
     = h.
\end{equation}

\item \label{entcal4}  If we further assume \eqref{eq:BD} then, for every $\epsilon>0$, there
exists a measurable partition $\partP$ of $(X,\mu)$ with
$H(\partP)<\infty$ such that for $\mu$-a.e.\ $(\omega,x)\in X$,
\begin{equation}\label{eq:ufutureBD}
    \lim_{n\to\infty} -\frac{1}{n}
         \log\mu_{(\omega,x)}^{\xi^u\vee\zeta^+}
         \bigl((\zeta\vee\Pcal)^{[0,n)}(\omega,x)\bigr)
     = h,
\end{equation}
\begin{equation}\label{eq:u0BD}
  \lim_{n\to\infty} -\frac{1}{n}
  \log\mu_{(\omega,x)}^{\xi^u\vee\zeta^+}
  \bigl((\zeta\vee\Pcal)^{(-n,0]}(\omega,x)\bigr)= 0,
\end{equation}
and
\begin{equation}\label{eq:ubothBD}
    \lim_{n\to\infty} -\frac{1}{n}
         \log\mu_{(\omega,x)}^{\xi^u\vee\zeta^+}
         \bigl((\zeta\vee\Pcal)^{(-n,n)}(\omega,x)\bigr)
     = h.
\end{equation}

\end{enumerate}
\end{proposition}

We note that applying \cref{prop:entcal}\ref{entcal3} with $F$ replaced by $F\inv$, 
we obtain an analogue of \cref{prop:entcal}\ref{entcal3} with  $\overline{\xi}^u$ replaced by $\overline{\xi}^s$, the directions of  time  as in \cref{prop:entcal}\ref{entcal2}, and  entropy $h$.

\subsection{Proof of \cref{prop:entcal}} 
The remainder of this section is devoted to the proof of \cref{prop:entcal}.
\subsubsection{Proof of \cref{prop:entcal}\ref{entcal1}}
We may select a measurable partition $\partP$ of $(X,\mu)$ with $H_\mu(\partP)<\infty$ and $$h-\epsilon/2\le h_\mu(F,\partP\mid \sigma(\zeta))\le h,\quad\quad \text{and} \quad\quad
\hat h-\epsilon/2\le h_\mu(F,\partP\vee \zeta)\le \hat h.$$
Then \eqref{eq:3_} follows from \Cref{SMB} and \eqref{eq:2_} follows from \cref{lemma:pastfinitepartabstract}.
For \eqref{eq:1_}, we have  
 \begin{align*}
  - \log\mu_{(\omega,x)}^\zeta
      \bigl(&(\zeta\vee\Pcal)^{(-n,n)}(\omega,x)\bigr)=  I_{\mu}((\zeta\vee\Pcal)^{(-n,n)}\mid\zeta)\\&=I_{\mu}((\zeta\vee\Pcal)^{[0,n)}\mid\zeta)+I_{\mu}((\zeta\vee\Pcal)^{(-n, 0)}\mid\zeta\vee (\zeta\vee\Pcal)^{[0,n)}).\end{align*}
and so  \eqref{eq:1_} follows from \cref{SMB,lemma:pastfinitepartabstract-half}.

\subsubsection{Proof of  \eqref{eq:4}, \eqref{eq:ufuture}, and \eqref{eq:ufutureBD}}
The upper bounds in \eqref{eq:4}, \eqref{eq:ufuture}, and \eqref{eq:ufutureBD} follow immediately from \cref{lemma:pastfinitepartabstract,cor:localentmax}.
Indeed, for \eqref{eq:4} we have 
$$-\frac{1}{n} \log \mu_{(\omega,x)}^{\xi^{s}\vee\zeta}\left((\zeta\vee\Pcal)^{(-n,0]}(\omega,x)\right)
\le - \frac 1 n \log  { \mu }^{  \xi^{s}\vee\zeta}_{(\omega,x)} \left( (    \xi^{s}\vee\zeta \vee   \partP) ^{(-n,0]}(\omega,x) \right)$$
By  \cref{lemma:pastfinitepartabstract,cor:localentmax,claim:slicingC}, the righthand side converges to 
$$ h_{ \mu }(   F, \partP \vee \xi^s \vee\zeta )= h_{ \mu }(   F, \xi^s \vee\zeta )=\hat h. $$
The upper bounds in \eqref{eq:ufuture} and \eqref{eq:ufutureBD} follow similarly by considering $F^{-1}$.

To justify the lower bounds in \eqref{eq:4}, \eqref{eq:ufuture}, and \eqref{eq:ufutureBD}, we apply \cref{prop:outcomeofadapted}. 
Consider  the lower bounds in \eqref{eq:4}.  Given $\delta>0$, let  $\partP$ be as guaranteed by  either conclusion \ref{adaptedok1} or \ref{adaptedok2} of \cref{prop:outcomeofadapted} (for this $\delta>0$.)  Set 
\[A_{\delta}=\left\{(\omega, x)\in X: \partP(\omega,x) \cap\tau\inv \bigl( \{\tau(\omega)\}\times W^{s}_{\loc,\delta}(\omega,x)\bigr)\subset \xi^s(\omega,x)\right\}.\]
From the defining properties of  $\xi^{s}$, we have  $\mu(A_{\delta})\to 1$ as $\delta\to 0$. In particular, we may select $\delta>0$ with  $\mu(A_{\delta})>0$.  

Given $(\omega,x)\in X$ 
and  given $n\ge 0$, let 
$R(\omega,x;n)=\sup \{r< n: F^{-r}(\omega,x)\in A_\delta\}$.  
 For all  $n>0$ sufficiently large, consider 
$$
(\omega',y)\in (\partP\vee \zeta)^{(-n,0]}(\omega,x) \cap \xi^s(\omega,x)\subset (\partP\vee \zeta)^{[-R(\omega,x;n),0]}(\omega,x) \cap \xi^s(\omega,x).
$$
By definition of $\xi^s$, $\tau(\omega) =\tau(\omega')$; similarly, using that either $(\omega',y)\in F^{n-1}\zeta(\omega,x) $ in the case  $\Xi(\Csc)\prec\zeta$ or 
that $(\omega',y)\in \partP^{(-n,0]}(\omega,x) $ in the case  \eqref{eq:BD}, 
$$f_{\theta ^j(\omega) }=f_{\theta ^j(\omega') } \quad \text{and} \quad  \tau(\theta^j(\omega) ) = \tau(\theta^j(\omega'))$$
for all $-n+1\le j\le 0$.  
  By  \cref{prop:outcomeofadapted},
$$f_{\omega'}^{(j)} (y)\in W^{s}_{\loc,\delta}(\theta^{(j)}(\omega),f_\omega^{(j)}(x))$$ for all $-R(\omega,x;n)\le j\le 0$.  
 By the definition of $A_\delta$,
$$F^{-R(\omega,x;n)}(\omega',y)\in (\xi^s\vee \zeta \vee \partP)(F^{-R(\omega,x;n)}(\omega,x))$$
and so
\begin{align*}
(\partP\vee \zeta)^{(-n,0]}(\omega,x) 
\cap \xi^s(\omega,x) 
\subset \bigr(\xi^s\vee \zeta\bigr)^{[-R(\omega,x;n),0]}(\omega,x). 
\end{align*}

Let $T(\omega,x)=\min \{n\ge 1: F^{-n}(\omega,x)\in A_\delta\}$. 
Let $R\colon A_\delta\to A_\delta$, $R(\omega,x)= F^{-T(\omega,x)}(\omega,x)$, denote the first (past) return map to $A_\delta$.  
Since $T\in L^1(\mu\restrict {A_\delta})$, for a.e.~$(\omega,x)\in A_\delta$, we have  $$\lim_{n\to \infty}\frac 1 n T(R^n(\omega,x)) = 0.$$ 
Thus, for a.e.~$(\omega,x)\in X$, 
\begin{equation}\label{eq:limitvsliminf}
\begin{split}
\liminf_{n\to \infty} \frac{R(\omega,x;n)}{n}
&\ge 
	\liminf_{n\to \infty} \frac{n - T(F^{-R(\omega,x;n)}(\omega,x))}{n}\\&= 1 - 	\limsup_{n\to \infty} \frac{ T(F^{-R(\omega,x;n)}(\omega,x))}{n}=1
\end{split}
\end{equation}
and so $ \frac{R(\omega,x;n)}{n}\to 1$.  
Let
\[g(\omega,x)=-\log\mu^{\zeta\vee\xi^{s}}_{(\omega,x)} \left((F(\zeta\vee\xi^{s}))(\omega,x)\right)=I_{\mu}(F(\zeta\vee\xi^{s})\mid\zeta\vee\xi^{s})(\omega,x).\] 
By the pointwise ergodic theorem and \cref{claim:slicingC}, 
\begin{align*}
 \liminf_{n\to \infty}-\frac{1}{n} &\log \mu_{(\omega,x)}^{\xi^{s}\vee\zeta} \left((\zeta\vee\Pcal)^{(-n,0]}(\omega,x)\right)\\
 &
\ge  \liminf_{n\to \infty}  -\frac{1}{n} \log \mu_{(\omega,x)}^{\xi^{s}\vee\zeta}\left(\bigr(
\xi^s\vee \zeta\bigr)^{[-R(\omega,x;n),0]
}(\omega,x)\right)
\\
&=   \lim_{n\to \infty} - \frac{R(\omega,x;n)}{n}  \frac{1}{R(\omega,x;n)} \log \mu_{(\omega,x)}^{\xi^{s}\vee\zeta}\left(\bigr(
\xi^s\vee \zeta\bigr)^{[-R(\omega,x;n),0]
}(\omega,x)\right)
\\
&=   \lim_{n\to \infty}  \frac{R(\omega,x;n)}{n}  \frac{1}{R(\omega,x;n)} \sum_{k=0}^{R(\omega,x;n)-1}g(F^{-k}(\omega,x))
\\
&=  h_{ \mu }(   F, \xi^s \vee \zeta ) =   \hat h.
\end{align*}

Replacing $F$ with $F^{-1}$, the same argument establishes the lower bounds in \eqref{eq:ufuture} and \eqref{eq:ufutureBD}.

\subsubsection{Proof of \cref{prop:entcal}\ref{entcal2}--\ref{entcal4}}

The existence and vanishing of the limits in \eqref{eq:s0}, \eqref{eq:u0}, and \eqref{eq:u0BD} follow from the generating properties in \cref{subsec:subord} of the partitions and by \cref{smbzero}.

The two-sided estimates \eqref{eq:sboth}, \eqref{eq:uboth}, and \eqref{eq:ubothBD} also now follow immediately.  
For instance we bound \eqref{eq:sboth} above and below by 
\begin{align*}
  -\frac 1 n \log\mu_{(\omega,x)}^{\xi^s\vee\zeta} &
         \big((\zeta\vee\Pcal)^{(-n,0]}(\omega,x)\bigr)\\
&\le   -\frac 1 n \log\mu_{(\omega,x)}^{\xi^s\vee\zeta} 
         \big((\zeta\vee\Pcal)^{(-n,n)}(\omega,x)\bigr)\\
         &\le     -\frac 1 n \log\mu_{(\omega,x)}^{\xi^s\vee\zeta} 
         \big((\xi^s\vee\zeta\vee\Pcal)^{(-n,n)}(\omega,x)\bigr)\\
                  &=  \frac 1 n I_{\mu}((\xi^s\vee\zeta\vee\Pcal)^{(-n,0]}\mid\xi^s\vee \zeta)(\omega,x)\\
&\quad\quad      +       \frac 1 n I_{\mu}((\xi^s\vee \zeta\vee\Pcal)^{(0,n)}\mid\ (\xi^s\vee \zeta\vee\Pcal)^{(-n,0]})(\omega,x).
\end{align*}
The lower bound is controlled by \eqref{eq:4}.  
For the upper bound, 
$$ \frac 1 n I_{\mu}((\xi^s\vee\zeta\vee\Pcal)^{(-n,0]}\mid\xi^s\vee \zeta)(\omega,x)\to \hat h$$
by \cref{lemma:pastfinitepartabstract,cor:localentmax,claim:slicingC}.  Then by    \cref{lemma:pastfinitepartabstract-full}
$$ \frac 1 n I_{\mu}((\xi^s\vee \zeta\vee\Pcal)^{(0,n)}\mid\ (\xi^s\vee \zeta\vee\Pcal)^{(-n,0]})(\omega,x)\to 0.$$ 

The proposition then follows.

\section{Proof of Main theorem}\label{sec:proofmain}

With the preliminaries in \cref{sec:lyapunov,sec:entropies}, we now present 
 the proof of \cref{thm:prod}.  After some notational conventions, the proof of both cases of \cref{thm:prod} follows from the same argument.  
 
Our arguments here are slightly different but mostly follow \cite[\S5--\S6]{MR2218998} and \cite[\S5--\S6]{MR1709302}.

\subsection{Standing notation}\label{sec:6notation}
Throughout, we assume  \eqref{eq:L1}  of \cref{sec:standing ass skew}.
Recall we fix a sub-invariant, measurable partition $\zeta$.  We assume either \begin{enumerate}
\item 
 $\Xi(\sigC)\prec \zeta$ and $h_{\mu}(F,\zeta\vee \Xi(\Csc)) <\infty$ (which implies $\maxent_\mu(F\mid \zeta)<\infty$), 
\item or that  hypothesis \eqref{eq:BD} of \cref{sec:standing ass skew} holds.
\end{enumerate}

Recall that the full measure set $\Lambda_0$ is the set of Oseledec generic points on which (two-sided) $(\epsilon_0,\ell)$-Lyapunov charts from \cref{Lyapcharts} are defined.
From \Cref{prpo:sudimm},  we select $$\Gamma_{0}\subset X$$ with $\mu(\Gamma_{0})=1$ such that the stable and unstable dimensions in \Cref{prpo:sudimm} are defined for every $(\omega,x)\in \Gamma_{0}$.

We define the following notation, depending on which  case of  \cref{thm:prod} we consider:  
\begin{align*}
h^u&=h=h_{\mu}(F\inv \mid \sigma(\zeta^+))= h_{\mu}(F \mid \sigma(\zeta^+))= h_{\mu}(F\mid \sigma(\zeta)).
\\  
h^s&=\hat{h}=h_{\mu}^{\refine}(F\mid \sigma(\zeta)).
\\ 
d^s &= \what{d}^{s}_{\zeta}.
\\d^u &= \begin{cases} 
      \bar{d}^u_{\zeta^{+}}, & \text{if } \Xi(\Csc)\prec\zeta \\
      \wcheck{d}_{\zeta^{+}}^{u}, & \text{if } \eqref{eq:BD}.
   \end{cases}
    \\  \wtd\xi^u&=\begin{cases}
 \bar\xi^u, & \text{if } \Xi(\Csc)\prec \zeta 
 \\ 
\xi^u, &\text{if } \eqref{eq:BD}.
\end{cases}  
\end{align*}

We also  define $$\wtd B(x,e^{-n}):= \pi_2 \inv (B(x, e^{-n}))= \Omega\times B(x,e^{-n}).$$ 
With this notation,  the proof of \cref{thm:prod} will follow directly from 
\cref{prop:UB,prop:LB} below. 

\subsection{Choice of finite entropy partition}
Recall we fixed a positive integer $a\in\Zbb_+$ such that  $a(\lambda_{\min}-3\epsilon_0)>1$.

 The following summarizes \cref{prpo:sudimm},  \cref{LY12.4.1},  \cref{lem:inball}, \cref{lem:chart-separation}, and \cref{prop:entcal}.  
\begin{proposition}\label{lem:approxa}
 For every  $\epsilon \in (0,1)$ sufficiently small there exist
\begin{itemize}  
\item a measurable set $\Gamma_\epsilon\subset \Gamma_0\subset X$ with $\mu(\Gamma_\epsilon)>1-\epsilon$, 
\item a positive integer $n_\epsilon=n_\epsilon(\epsilon)\in \Nbb$; 
\item a countable partition $\mathcal K$ of $(X,\mu)$;
\item constants $0<c_0<1<\ell_0<\infty$ and   $L= \lceil \log \ell_0\rceil$;
\item  a choice of $\delta>0$ in the construction of $\partP$ 
\end{itemize}
such that for all {$n\ge n_\epsilon$} and for all $(\omega,x)\in \Gamma_\epsilon$, the following hold:
\begin{enumerate}[label*={(\arabic*)}, ref= {(\arabic*)} ]
\item\label{approx0} \begin{enumerate}[label={(\alph*)}, ref={(\arabic{enumi}\alph*)}]
\item\label{approx0:chart-bound} The function $\ell\colon \Lambda_0\to [1,\infty)$ in \cref{Lyapcharts} is bounded from above by $\ell_0$  on $\Gamma_\epsilon$.
\item\label{approx0:slice-bound} The functions $c^s,c^u \colon \Lambda_0\to (0,1]$ in \cref{claim:localslice} are bounded below by $c_0$  on $\Gamma_\epsilon$.  
\item\label{approx0:scale} $\delta + 2e^{-n}\le c_0$.  
\item The conditional measure $\mu_{(\omega,x)}^{\zeta}$ is defined and further decomposes into families of conditional measures $\{\mu_{(\omega',y)}^{\zeta\vee \xi^s}: (\omega',y)\in \zeta(\omega,x) \}$ and $\{\mu_{(\omega',y)}^{\zeta^+\vee \td  \xi^u}: (\omega',y)\in \zeta(\omega,x) \}$
\end{enumerate}
\item \label{approx1} 
\begin{enumerate}[label={(\alph*)}, ref={(\arabic{enumi}\alph*)}]
\item\label{eq:both} $\displaystyle e^{-na(h+\hat{h}+\epsilon)}\le \mu^{\zeta}_{(\omega,x)}((\zeta\vee \Pcal)^{[-na,na]}(\omega,x)) \le e^{-na(h+\hat{h}-\epsilon)}$
\item\label{approx1:future} $\displaystyle e^{-na(h+\epsilon)}\le \mu^{\zeta}_{(\omega,x)}((\zeta\vee \Pcal)^{[0,na]}(\omega,x)) \le e^{-na(h-\epsilon)}$
\item\label{eq:pastzetaa} $\displaystyle e^{-na(\hat{h}+\epsilon)}\le \mu^{\zeta}_{(\omega,x)}((\zeta\vee \Pcal)^{[-na,0]}(\omega,x)) \le e^{-na(\hat{h}-\epsilon)}$

\item\label{eq:shbar} $\displaystyle e^{-na(\hat{h}+\epsilon)}\le \mu^{\xi^{s}\vee\zeta}_{(\omega,x)}((\zeta\vee \Pcal)^{[-na,na]}(\omega,x)) \le e^{-na(\hat{h}-\epsilon)}$

\item\label{eq:pasts} $\displaystyle e^{-na(\hat{h}+\epsilon)}\le \mu^{\xi^{s}\vee\zeta}_{(\omega,x)}((\zeta\vee \Pcal)^{[-na,0]}(\omega,x)) \le e^{-na(\hat{h}-\epsilon)}$
\item\label{eq:ubothu}  $\displaystyle e^{-na(h+\epsilon)}\le \mu^{\wtd{\xi}^{u}\vee\zeta^+}_{(\omega,x)}((\zeta\vee \Pcal)^{[-na,na]}(\omega,x)) \le e^{-na(h-\epsilon)}$
\item\label{eq:35u} $\displaystyle e^{-na(h+\epsilon)}\le \mu^{\wtd{\xi}^{u}\vee\zeta^+}_{(\omega,x)}((\zeta\vee \Pcal)^{[0,na]}(\omega,x)) \le e^{-na(h-\epsilon)}$
\end{enumerate}

\item \label{approx5}\label{unstableP+} 
$\displaystyle \mu_{(\omega,x)}^{\td{\xi}^u \vee\zeta^+ \vee{\Pcal^+}} ((\zeta\vee\Pcal)^{[0,an]}(\omega,x))\le e^{-an(h-\epsilon)}$

\item\label{approxball} \label{approx4}\begin{enumerate}[label={(\alph*)}, ref={(\arabic{enumi}\alph*)}]
\item\label{approx4:two-sided} $\displaystyle (\zeta\vee\Pcal)^{[-an,an]} (\omega,x)\subset \widetilde B(x,e^{-n})$
\item\label{approx4:stable} $\displaystyle \bigl(\xi^{s}\vee (\zeta\vee\Pcal)^{[-an,0]}\bigr)(\omega,x)\subset \widetilde B(x,e^{-n})$
\item\label{approx4:unstable} $\displaystyle \bigl(\wtd \xi^{u}\vee (\zeta\vee\Pcal)^{[0,an]}\bigr)(\omega,x)\subset \widetilde B(x,e^{-n})$
\item\label{approx4:stable-separation} For   any $K\in \partK$ and $(\omega,x), (\omega',y), (\omega'',z)\in \Gamma_\epsilon$ with 
\begin{itemize}
\item $(\omega,x), (\omega'',z)\in  K$,
\item  $z\in B(x,e^{-n})$,
\item $(\omega',y)\in \xi^s(\omega,x)$, and 
\item  $B(y,4e^{-n +L})\cap B(x,4 e^{-n  +L})= \emptyset$
\end{itemize}
then 
\[
\pi_2\bigl((\zeta\vee\Pcal)^{[-an,0]}(\omega'',z)\bigr)\cap B(y, e^{-n})=\emptyset.
\]
\item\label{approx4:unstable-separation} For   any $K\in \partK$ and $(\omega,x), (\omega',y), (\omega'',z)\in \Gamma_\epsilon$ with 
\begin{itemize}
\item $(\omega,x), (\omega'',z)\in  K$,
\item  $z\in B(x,e^{-n})$,
\item $(\omega',y)\in \wtd \xi^u(\omega,x)$, and 
\item  $B(y,4e^{-n +L})\cap B(x,4 e^{-n  +L})= \emptyset$
\end{itemize}
then 
\[
\pi_2\bigl((\zeta\vee\Pcal)^{[0,an]}(\omega'',z)\bigr)\cap B(y, e^{-n})=\emptyset.
\]

\end{enumerate}

\item \label{approx3}  
\begin{enumerate}[label={(\alph*)}, ref={(\arabic{enumi}\alph*)}]
\item\label{eq:ds} $\displaystyle e^{-n( d^{s}+ \epsilon)}\le \mu_{(\omega,x)}^{\xi^{s}\vee\zeta} (\widetilde B(x,e^{-n}))\le 
\mu_{(\omega,x)}^{\xi^{s}\vee\zeta} (\widetilde B(x,4e^{-n}))\le 
e^{-n( d^{s}-\epsilon)}$
\item\label{eq:du} $\displaystyle e^{-n({d}^{u}+\epsilon)}\le
\mu_{(\omega,x)}^{\td{\xi}^{u}\vee\zeta^{+}} (\widetilde B(x,e^{-n}))\le 
\mu_{(\omega,x)}^{\td{\xi}^{u}\vee\zeta^{+}} (\widetilde B(x,4e^{-n}))
\le e^{-n({d}^{u}-\epsilon)}$

\end{enumerate}
\end{enumerate}

\end{proposition}

 Using \Cref{prop:Lebdensity,prop:integral}, we further refine $\Gamma_\epsilon$.  

\begin{proposition}\label{prop:Lebdens}
There are  $\Gamma'\subset \what \Gamma \subset \Gamma_\epsilon$ with   $\mu(\what \Gamma)>1-2\epsilon$, $\mu(\Gamma')>1-3\epsilon$ and  $n'\ge \hat n\ge n_\epsilon$ such that

\begin{enumerate}[label={(\arabic*)}, ref={(\arabic*)}]
\item\label{Lebdens:unstable}{ $\mu^{\wtd{\xi}^{u}\vee\zeta^+}_{(\omega,x)}(\widetilde B(x,e^{-n})\cap  \Gamma_\epsilon\cap \partK(\omega,x))\ge e^{-n ({d}^u+\epsilon)}$ for all $(\omega,x) \in \what \Gamma$ and all $n\ge \hat n$,}
\item\label{Lebdens:stable}{ $\mu^{ {\xi}^{s}\vee\zeta}_{(\omega,x)}(\widetilde B(x,e^{-n})\cap  \Gamma_\epsilon\cap \partK(\omega,x))\ge e^{-n ({d}^s+\epsilon)}$ for all $(\omega,x) \in \what \Gamma$ and all $n\ge \hat n$,}

\item\label{Lebdens:unstable-future}\label{LY12.4.2}
 { $\mu^{\wtd{\xi}^{u}\vee\zeta^+\vee \partP^+}_{(\omega,x)}(\widetilde B(x,e^{-n})\cap  \Gamma_\epsilon)\ge e^{-n ({d}^u+\epsilon)}$ for all $(\omega,x) \in \what \Gamma$ and all $n\ge \hat n$,}

\item\label{Lebdens:stablefull} $\mu^{{\xi}^{s}\vee\zeta}_{(\omega,x)}(\widetilde B(x,e^{-n})\cap \what \Gamma)\ge e^{-n ({d}^s+\epsilon)}$ for all $(\omega,x)\in \Gamma'$ and all $n\ge n'$.
\end{enumerate}

\end{proposition}

\subsection{Rectangles, covering sets, and naive bounds}
Recall our fixed $\epsilon>0$ and $\Gamma_\epsilon$ and $n_\epsilon$  defined in \Cref{lem:approxa}.

\subsubsection{Rectangles and covering sets}\label{sss:recandcover}
For $n\ge  0$ define the following collections of ``rectangles:'' 
\begin{align*}\Rcal_{n,\zeta} &= \bigl\{(\zeta\vee\Pcal)^{[-an,an]}(\omega',x'): (\zeta\vee\Pcal)^{[-an,an]}(\omega',x')\cap \Gamma_\epsilon \neq \emptyset\bigr\},\\
\Fcal_{n,\zeta}^{s}&= \bigl\{(\zeta\vee\Pcal)^{[-an,0]}(\omega',x'): (\zeta\vee\Pcal)^{[-an,0]}(\omega',x')\cap \Gamma_\epsilon\neq\emptyset\bigr\},\\
\Fcal_{n,\zeta} ^u &= \bigl\{(\zeta\vee\Pcal)^{[0,an]}(\omega',x'):(\zeta\vee\Pcal)^{[0,an]}(\omega',x')\cap \Gamma_\epsilon\neq \emptyset\bigr\}.\end{align*}

Let $\eta^{s}=\zeta\vee\xi^s$, $\eta^u=\wtd{\xi}^{u}\vee \zeta^+$.  
For a subset $A\subset X$, $(\omega,x)\in X$, 
 $*\in\{s,u\}$, and for each $n$, we define the following subsets.   
\begin{align*}
N^{\zeta}_{(\omega,x)}(n,A)&=\bigl\{R\in \Rcal_{n,\zeta} : R\cap A\cap \zeta(\omega,x)\cap \Gamma_\epsilon\neq\emptyset\bigr\},\\
N^{*,\zeta}_{(\omega,x)}(n,A)&=N^{\zeta}_{(\omega,x)}(n,\eta^{*}(\omega,x)\cap A),\\
\widehat{N}^{*,\zeta}_{(\omega,x)}(n,A)&=\bigl\{R^*\in \Fcal_{n,\zeta} ^* : R^*\cap A \cap \zeta(\omega,x)\cap \Gamma_\epsilon\neq\emptyset\bigr\}.
\end{align*}

\subsubsection{Naive counting bounds}
We have the following upper bound on $\#N^{\zeta}_{(\omega,x)}(n,\cdot )$.  
\begin{lemma}\label{lem:2.9} For  $(\omega,x)\in X$, $A\subset X$, and $n\ge 0$, 
\begin{eqnarray*}\#N^{\zeta}_{(\omega,x)}(n,A) \le \#\widehat{N}^{s,\zeta}_{(\omega,x)}(n,A) 
\cdot \#\widehat{N}^{u,\zeta}_{(\omega,x)}(n,A)\end{eqnarray*}
\end{lemma}
\begin{proof}   
Consider $R\in N^{\zeta}_{(\omega,x)}(n,A)$.  Then $R$ is of the form \[R=(\zeta\vee\Pcal)^{[-an,an]}(\omega',x')=(\zeta\vee \partP)^{[0,an]}(\omega',x')\cap (\zeta\vee\Pcal)^{[-an,0]}(\omega',x')\] for some $(\omega',x')\in \zeta(\omega,x)\cap A\cap \Gamma_\epsilon$.   Let \[R^{s}=(\zeta\vee\Pcal)^{[-an,0]}(\omega',x')
 \textrm{ and  } 
R^{u}=(\zeta\vee\Pcal)^{[0,an]}(\omega',x')
.\]
Then $R\cap A\cap \Gamma_\epsilon\neq\emptyset$ implies $R^{s}\cap A\cap \Gamma_\epsilon\neq\emptyset $ and $R^{u}\cap A\cap \Gamma_\epsilon\neq\emptyset.$
In particular, $R\mapsto (R^s,R^u)$ determines an injective map \[N^{\zeta}_{(\omega,x)}(n,A)\to\widehat{N}^{s,\zeta}_{(\omega,x)}(n,A)\times \widehat{N}^{u,\zeta}_{(\omega,x)}(n,A).  \qedhere\] 
\end{proof}

Similarly, we have the following lower bound on $\#N^{\zeta}_{(\omega,x)}(n,\cdot)$.  
\begin{lemma} \label{lem:count} For any $(\omega,x)\in X$ and $A'\subset  A\subset  X$,
\[
\begin{split}
\# N^{\zeta}_{(\omega,x)}(n,A)
&\ge \#\what N^{s,\zeta}_{(\omega,x)}(n,A') \\
&\quad {}\cdot \min \bigl\{
\#\what N^{u,\zeta}_{(\omega,x)}(n,A\cap R^s):{R^s\in \what N^{s,\zeta}_{(\omega,x)}(n,A')}\bigr\}
\end{split}
\]
(with an appropriate interpretation that $0\cdot \min \emptyset =0$).
\end{lemma}
\begin{proof}As in \cref{lem:2.9}, given $R\in N^{\zeta}_{(\omega,x)}(n,A)$ there exist unique  elements $R^s\in  \what N^{s,\zeta}_{(\omega,x)}(n,A)$ and  $R^u\in \what N^{u,\zeta}_{(\omega,x)}(n,A\cap R^s)$ with $R= R^s\cap R^u$.
Then
\begin{align*}
\# N^{\zeta}_{(\omega,x)}&(n,A)\\
&= \sum_{R^s\in \what N^{s,\zeta}_{(\omega,x)}(n,A)} \#\what N^{u,\zeta}_{(\omega,x)}(n,A\cap R^s)\\
&\ge \sum_{R^s\in \what N^{s,\zeta}_{(\omega,x)}(n,A')} \#\what N^{u,\zeta}_{(\omega,x)}(n,A\cap R^s)\\
&\ge \#\what N^{s,\zeta}_{(\omega,x)}(n,A') \cdot \min \{
\#\what N^{u,\zeta}_{(\omega,x)}(n,A\cap R^s):{R^s\in \what N^{s,\zeta}_{(\omega,x)}(n,A')}\}.\qedhere
\end{align*}
\end{proof}

\subsubsection{Naive   dimension and entropy bounds} 

The following bounds are immediate from definition and \cref{lem:approxa}.
\begin{lemma}\label{lem:roughbda}
For every $(\omega,x)\in \Gamma_\epsilon$ and every $*\in\{s,u\}$, $A\subset X$, and 
$n\ge n_\epsilon$, 
\begin{equation*}
\#\what{N}^{*,\zeta}_{(\omega,x)}(n, A)\le e^{an(h^{*}+\epsilon)}
\end{equation*}
\end{lemma}
\begin{proof} 
We have  \[\zeta(\omega,x)=\bigcup_{(\omega',x')\in \zeta(\omega,x)}(\zeta\vee\Pcal)^{[-an,0]}(\omega',x')
=\bigcup_{(\omega',x')\in \zeta(\omega,x)}(\zeta\vee\Pcal)^{[0,an]}(\omega',x'). \] 
Applying 
\cref{lem:approxa}\ref{eq:pastzetaa} or \ref{approx1:future}  
to each $R\in  \hat{N}^{*,\zeta}_{(\omega,x)}(n,A)$,
we have  
\begin{align*}1\ge \mu^{\zeta}_{(\omega,x)}(\zeta(\omega,x))&\ge \sum_{R\in \hat{N}^{*,\zeta}_{(\omega,x)}(n,A)}\mu_{(\omega,x)}^{\zeta}(R)\\&\ge \#\what{N}^{*,\zeta}_{(\omega,x)}(n, A)\cdot  e^{-an({h}^* +\epsilon)}. \qedhere \end{align*} \end{proof}

\begin{lemma}\label{lem:roughbd2a} 
For every $*\in\{s,u\}$, $(\omega,x)\in \what \Gamma$, and $n\ge  \hat n$, 
$$ \exp(n(-d^*+a h^*-2a\epsilon))
 \le\#N^{*,\zeta}_{(\omega,x)}\bigl(n,\widetilde B(x,e^{-n} )\cap \partK(\omega,x) \bigr).
$$
\end{lemma}  
 \begin{proof} 
 By \cref{prop:Lebdens}\ref{Lebdens:unstable} or \ref{Lebdens:stable}, for $n\ge \what{n}$, we have
\begin{align}\exp(-(d^*+\epsilon)n) \le \mu_{(\omega,x)}^{\eta^*}\bigl(\widetilde B(x,e^{-n}) \cap \Gamma_\epsilon\cap \partK(\omega,x)\bigr)
\label{eqroughlow}
\end{align} 
Applying  \cref{lem:approxa}\ref{eq:shbar} or  \ref{eq:ubothu}  
to every  
$$R=(\zeta\vee\Pcal)^{[-an,an]}(\omega'',z)
\in N^{*,\zeta}_{(\omega,x)}\bigl(n,\widetilde B(x,e^{-n})\cap \partK(\omega,x)\bigr)$$
we have 
\begin{equation}\label{eq:25gooa}   \exp(-an(h^*+\epsilon)) \le \mu_{(\omega,x)}^{\eta^*}\bigl(R\bigr)
=\mu_{(\omega'',z)}^{\eta^*}\bigl(R\bigr)\le \exp(-an(h^*-\epsilon))\end{equation} 
Using \eqref{eq:25gooa}, we have 
\begin{align*}\mu_{(\omega,x)}^{\eta^*}&\bigl(\widetilde B(x,e^{-n})\cap \Gamma_\epsilon\cap \partK(\omega,x)\bigr)=\mu_{(\omega,x)}^{\eta^*}\bigl(\widetilde B(x,e^{-n})\cap\eta^*(\omega,x)\cap \Gamma_\epsilon\cap \partK(\omega,x)\bigr)\\
&\le \sum_{R\in N^{*,\zeta}_{(\omega,x)}(n,\widetilde B(x,e^{-n})\cap  \partK(\omega,x))}\mu_{(\omega,x)}^{\eta^*}(R)\\
&\le \sum_{R\in N^{*,\zeta}_{(\omega,x)}(n,\widetilde B(x,e^{-n})\cap \partK(\omega,x))}e^{-an(h^* -\epsilon)}\\
&= e^{-an(h^* -\epsilon)}\cdot \#N^{*,\zeta}_{(\omega,x)}(n,\widetilde B(x,e^{-n}) \cap \partK(\omega,x)).
\end{align*} 
The result follows combined with \eqref{eqroughlow}.  
\end{proof}

\begin{lemma} \label{lem:2.10}
For every $*\in\{s,u\}$, $(\omega,x)\in \Gamma_\epsilon$, and $n\ge n_\epsilon$, 
\begin{equation*}
\#N^{*,\zeta}_{(\omega,x)}\bigl(n,\widetilde B(x,e^{-n})  \cap \Gamma_\epsilon\cap \partK(\omega,x) \bigr)\le  \exp(n(-d^*+a h^*+2a\epsilon))
\end{equation*}
\end{lemma}
\begin{proof}
If $R\in N^{*,\zeta}_{(\omega,x)}\bigl(n,\widetilde B(x,e^{-n})  \cap \Gamma_\epsilon\cap \partK(\omega,x) \bigr)$ then 
 $$\mu_{(\omega,x)}^{\eta^*}(R) \ge e^{-an({h^*}+\epsilon)}$$ by \Cref{lem:approxa}\ref{eq:shbar} or  \ref{eq:ubothu}.   
 By  \Cref{lem:approxa}\ref{approx4:two-sided}, $R\cap \eta^{*}(\omega,x)\subset \wtd{B}(x,4e^{-n})$.  Thus 
\begin{align*}
\mu_{(\omega,x)}^{\eta^*}&(\wtd{B}(x,4e^{-n}))\\
&\ge \sum_{R\in N^{*,\zeta}_{(\omega,x)}(n,\widetilde B(x,e^{-n})  \cap \Gamma_\epsilon\cap \partK(\omega,x) )} \mu_{(\omega,x)}^{\eta^*}(R)\\
&\ge \#N^{*,\zeta}_{(\omega,x)}\bigl(n,\widetilde B(x,e^{-n})  \cap \Gamma_\epsilon\cap \partK(\omega,x) \bigr) \cdot e^{-an({h^*}+\epsilon)}. \qedhere \end{align*}\end{proof}

Similarly, we have the following.  
\begin{lemma}\label{lem:2.8} For each $(\omega,x)\in {\Gamma_\epsilon}$ and $n\ge n_\epsilon$,
\begin{equation*}
\begin{split}
&\mu_{(\omega,x)}^{\zeta}\bigl(\widetilde B(x,e^{-n})\cap \Gamma_\epsilon \cap \partK(\omega,x)\bigr)\\
&\qquad\le  \# N^{\zeta}_{(\omega,x)}\bigl(n, \widetilde B(x,e^{-n})\cap \partK(\omega,x)\bigr)\cdot \exp(-an(h+\hat{h}-\epsilon))
\end{split}
\end{equation*}
\end{lemma}

\begin{proof}[Proof of \Cref{lem:2.8}] 
By  \cref{lem:approxa}\ref{eq:both} \begin{align*} 
\mu_{(\omega,x)}^{\zeta}\bigl(\widetilde B(x,e^{-n})\cap &\Gamma_\epsilon\cap \partK(\omega,x)\bigr)
\\&\le \sum_{R\in N^{\zeta}_{(\omega,x)}(n, \widetilde B(x,e^{-n})\cap \partK(\omega,x))} \mu_{(\omega,x)}^{\zeta}(R)\\
&\le \# N^{\zeta}_{(\omega,x)}\bigl(n, \widetilde B(x,e^{-n})\cap \partK(\omega,x)\bigr) e^{-an(h+\hat{h}-\epsilon)}.\qedhere
\end{align*}

\end{proof}

\subsection{Upper bound on dimension}

Recall we write 
$$\wtd B(x,e^{-n}):= \pi_2 \inv (B(x, e^{-n}))= \Omega\times B(x,e^{-n})$$ 
for the  saturation of the Riemannian ball in the fibers.  
With the above notation, we are ready to prove the following.  
\begin{proposition}\label{prop:UB}
For $\mu$-a.e.\ $(\omega, x)\in X$, 
$$\limsup_{n\to \infty}\frac{\log\mu_{(\omega,x)}^\zeta(\wtd B(x,e^{-n}))}{\log(e^{-n})}\le {d}^u+d^s.$$
\end{proposition}

\begin{proof}[Proof of \cref{prop:UB}]
Fix $\epsilon>0$ sufficiently small and let $\Gamma'\subset \what \Gamma\subset \Gamma_\epsilon\subset \Gamma_0$ and $n'\ge \hat n\ge n_\epsilon$ be as in \cref{prop:Lebdens}.  
Fix any  $(\omega_0,x_0)\in\Gamma'$ and $n\ge n'$.   
We consider the intersection of the saturated ball 
$\wtd B(x_0,2e^{-n})= \pi_2 \inv (B(x_0,2 e^{-n}))$  and $\zeta(\omega_0,x_0)\cap \Gamma_\epsilon$; we follow the notation of \cref{sss:recandcover} to cover this intersection by atoms of iterated partitions.

Recall $N_1=\# N^{\zeta}_{(\omega_0,x_0)} (n, \wtd B(x_0,{ 2}e^{-n}))$ denotes the number of distinct $
(\zeta\vee\partP)^{[-an,an]}$-atoms intersecting $\widetilde B(x_0,{  2}e^{-n})\cap \Gamma_\epsilon\cap \zeta(\omega_0,x_0)$.  
For each $$R\in N^{\zeta}_{(\omega_0,x_0)} (n, \widetilde B(x_0,{ 2}e^{-n})),$$ \cref{lem:approxa}\ref{approx4}  implies 
$R\subset \wtd B(x_0, 4e^{-n})$ and \cref{lem:approxa}\ref{eq:both} implies $\mu^{ \zeta}_{(\omega_0, x_0)}(R) \ge e^{-na(h+\hat{h}+\epsilon)}$.  Thus, 
\begin{equation}\label{upperbd0} e^{-na(h+\hat{h}+\epsilon)}N_1\le \mu^{ \zeta}_{(\omega_0, x_0)}(\wtd B({ x_0},4e^{-n})).\end{equation}

Next,  $N_3'= \# \what N^{s,\zeta}_{(\omega_0,x_0)}(n, (\xi^s\vee\zeta)(\omega_0,x_0)\cap \wtd B(x_0,e^{-n})\cap \what \Gamma)$ denotes  the number of distinct $(\zeta\vee\partP)^{[-an,0]}$-atoms intersecting $(\xi^s\vee\zeta)(\omega_0,x_0)\cap \wtd B(x_0,e^{-n})\cap \what \Gamma$.  For every $R\in \what N^{s,\zeta}_{(\omega_0,x_0)}(n, (\xi^s\vee\zeta)(\omega_0,x_0)\cap \wtd B(x_0,e^{-n})\cap \what \Gamma) $, 
by \cref{lem:approxa}\ref{eq:pasts} we have $\mu^{{\xi}^{s}\vee\zeta}_{{ (\omega_0,x_0)}}(R)\le \exp(-na(\hat{h}-\epsilon))$ and so  by \cref{prop:Lebdens}\ref{Lebdens:stablefull},
\begin{equation*}
N_3'\exp(-na(\hat{h}-\epsilon))\ge \mu^{{\xi}^{s}\vee\zeta}_{{ (\omega_0,x_0)}}(\wtd B(x_0,e^{-n})\cap \what \Gamma)\ge e^{-n (d^s+\epsilon)}.
\end{equation*}
With $N_3=\# \what N^{s,\zeta}_{(\omega_0,x_0)}(n, \wtd B(x_0,e^{-n})\cap \what \Gamma)$, the number of distinct $(\zeta\vee\partP)^{[-an,0]}$-atoms intersecting $\zeta(\omega_0,x_0)\cap \widetilde B(x_0,e^{-n})\cap  \what \Gamma$, we have 
\begin{equation}\label{upperbd2} 
N_3\ge N_3'\ge 
 e^{na(\hat{h}-\epsilon)} e^{-n (d^s+\epsilon)}.
\end{equation}

Consider now  any  $(\omega',y)\in \wtd B(x_0,e^{-n})\cap \what \Gamma\cap \zeta(\omega_0,x_0)$ such that $R= (\zeta\vee\partP)^{[-an,0]}(\omega',y)\in {\what N^{s,\zeta}}_{(\omega_0,x_0)}(n, \wtd B(x_0,e^{-n})\cap \what \Gamma)$.  
Then $N_2'(\omega',y)=\# \what N^{u,\zeta}_{(\omega',y)}(n, (\wtd\xi^u\vee\zeta^+\vee \partP^+)(\omega',y)\cap  \wtd B(x_0,2 e^{-n}))$ denotes  the number of distinct $(\zeta\vee\partP)^{[0,an]}$-atoms intersecting $(\wtd\xi^u\vee\zeta^+\vee \partP^+)(\omega',y)\cap \Gamma_\epsilon\cap \wtd B(x_0,2 e^{-n})$. 
For any $R'=(\zeta\vee \Pcal)^{[0,an]}(\omega'',z)\in
\what N^{u,\zeta}_{(\omega',y)}(n, (\wtd\xi^u\vee\zeta^+\vee \partP^+)(\omega',y)\cap  \wtd B(x_0,2 e^{-n}))$,
\cref{lem:approxa}\ref{unstableP+}  gives
\[\mu^{\wtd{\xi}^{u}\vee\zeta^{+}\vee \partP^+}_{(\omega',y)}(R')=\mu^{\wtd{\xi}^{u}\vee\zeta^{+}\vee \partP^+}_{(\omega'',z)}(R')\le e^{-an(h-\epsilon)}.\] 
Since $(\omega',y)\in  \wtd B(x_0,e^{-n})\cap  \what \Gamma$,
we have 
 $$\wtd  B({ x_0}, 2e^{-n}) \cap \Gamma_\epsilon\supset \wtd B({y} , e^{-n}) \cap \Gamma_\epsilon$$ and thus combined with \cref{prop:Lebdens}\ref{LY12.4.2}
\begin{align*}
N_2'(\omega',y) e^{-an(h-\epsilon)} &\ge
 \mu^{\wtd{\xi}^{u}\vee\zeta^{+}\vee \partP^+}_{(\omega',y)} (\widetilde B( x_0 , 2 e^{-n}) \cap \Gamma_\epsilon)\\&\ge 
 \mu^{\wtd{\xi}^{u}\vee\zeta^{+}\vee \partP^+}_{(\omega',y)} (\widetilde B( y , e^{-n}) \cap \Gamma_\epsilon)\ge 
 e^{-n (d^u+\epsilon)}. 
\end{align*}
 Finally,  $N_2(\omega',y)= \#\what N^{u,\zeta}_{(\omega',y)}(n, (\zeta\vee\partP)^{[-an,0]}(\omega',y) \cap  \wtd B(x_0,2 e^{-n}))$
 denotes the number of distinct $(\zeta\vee\partP)^{[0,an]}$-atoms intersecting $(\zeta\vee\partP)^{[-an,0]}(\omega',y)\cap  \wtd B(x_0,2 e^{-n})\cap \Gamma_\epsilon$.
Since  $(\zeta\vee \partP)^+(\omega',y)\subset (\zeta\vee\partP)^{[-an,0]}(\omega',y)$,
\begin{equation}\label{upperbd1a} 
N_2(\omega',y)\ge N_2'(\omega',y)\ge 
 e^{na({h}-\epsilon)} e^{-n (d^u+\epsilon)}.
\end{equation}

Since the bound \eqref{upperbd1a} holds for all $(\omega',y)\in \wtd B(x_0,e^{-n})\cap  \what \Gamma$, applying   \cref{upperbd2} and \cref{lem:count} with $A'= \wtd B(x_0,e^{-n})\cap \what \Gamma$ and $A=\wtd B(x_0,2e^{-n})\cap   \Gamma_\epsilon$ 
we have 
\begin{align*} N_1 
&=\#N^{\zeta}_{(\omega_0,x_0)} (n,A)\\
&\ge \#\what N^{s,\zeta}_{(\omega_0,x_0)}(n, A')
\cdot \min \{ \#  \what N^{u,\zeta}_{(\omega_0,x_0)}(n, { A\cap R}): R\in \what N^{s,\zeta}_{(\omega_0,x_0)}(n, A')\}\\
& \ge   N_3  e^{na({h}-\epsilon)} e^{-n (d^u+\epsilon)}\\
& \ge  e^{na({h}+\hat{h}-2\epsilon)}e^{-n ({d}^u+d^s+2\epsilon)}.\end{align*}
Combined with \cref{upperbd0}, 
\begin{align*} 
 \log\mu_{(\omega_0,x_0)}^{\zeta}&(\widetilde B(x_0,4e^{-n}))\\&
\ge
(-na(h+\hat{h}+\epsilon)) + (na({h}+\hat{h}-2\epsilon)) +({-n ({d}^u+d^s+2\epsilon)})\\
&\ge - n ({d}^u+d^s) -(3a+2)n\epsilon
\end{align*}
and thus for every $(\omega_0,x_0)\in \Gamma'$, 
\begin{align*}
\limsup_{n\to \infty}\frac{\log\mu_{(\omega_0,x_0)}^{\zeta}(\wtd B(x_0,e^{-n}))}{\log(e^{-n})}
&= \limsup_{n\to \infty}\frac{\log\mu_{(\omega_0,x_0)}^{\zeta}(\wtd B(x_0,4e^{-n-2}) )}{\log(e^{-n})}\\ &\le {d}^u+d^s + (3a+2)\epsilon.   \end{align*} 

The proposition then follows using that $\epsilon>0$ was arbitrary and that the associated sets $\Gamma'$ exhaust a full-measure subset of $X$.  
\end{proof}

\subsection{Comparison of $N^*$ and $\what{N}^*$}We note that for every $(\omega,x)$ and every $n$, 
$${\# \what{N}^{*,\zeta}_{(\omega,x)}\bigl(n,\widetilde B(x,e^{-n})\cap K\bigr)}\ge {\#N^{*,\zeta}_{(\omega,x)}\bigl(n,\widetilde B(x,e^{-n})\cap K\bigr)}.$$
In order to establish  the lower bound on dimension, we need the following (upper) bound on the errors introduced when moving between $N^*$ and $\what{N}^*$ in \cref{sss:recandcover}.
However, because of the geometry used in the counting, we must restrict to sets $K\in \partK$.

\begin{lemma}\label{lem:compareN}  Fix $(\omega_{0},x_{0})\in  \Gamma_\epsilon$ and let $K= \partK(\omega_0,x_0)$.  
For $\mu_{(\omega_{0},x_{0})}^{\zeta}$-a.e.~$(\omega',x')\in \zeta(\omega_{0},x_{0})\cap \what{\Gamma}\cap K$ and $*\in\{s,u\}$,
\begin{equation*}
\limsup_{n\to\infty}\left[ \frac{\# \what{N}^{*,\zeta}_{(\omega_0,x_0)}\bigl(n,\widetilde B(x',e^{-n})\cap K\bigr)}{\#N^{*,\zeta}_{(\omega',x')}\bigl(n,\widetilde B(x',e^{-n})\cap K\bigr)}\cdot e^{-5an\epsilon}\right]=0.
\end{equation*}
\end{lemma}

\begin{proof}[Proof of \Cref{lem:compareN}]

For each $k\in \Nbb$, let $F_{k}\subset \zeta(\omega_{0},x_{0})\cap \what{\Gamma}\cap K$ denote the set of points $(\omega',x')$ for which 
\[ \limsup_{n\to\infty} \frac{\# \what{N}^{*,\zeta}_{(\omega_0,x_0)}\bigl(n,\widetilde B(x',e^{-n})\cap K\bigr)}{\#N^{*,\zeta}_{(\omega',x')}\bigl(n,\widetilde B(x',e^{-n})\cap K\bigr)}\cdot e^{-5an\epsilon}\ge \frac{1}{k}.\] 
It is enough to show that $\mu_{(\omega_{0},x_{0})}^{\zeta}(F_{k})=0$ for every $k$. 

\def\bad{\mathcal {B}}
Consider a fixed $k$.  
 For every $(\omega',x')\in F_{k}$, there is an 
 infinite set of bad times $\bad(\omega',x')\subset \Nbb$ such that for every $m\in \bad(\omega',x')$,
  \begin{equation*}
 \begin{split}\# \what{N}^{*,\zeta}_{(\omega_0,x_0)}&\bigl(m ,\widetilde B(x', e^{-m })\cap K\bigr)
 \\&\quad \quad \ge \frac{1}{2k}\#N^{*,\zeta}_{(\omega',x')}\bigl(m ,\widetilde B(x',e^{-m })\cap K\bigr)\cdot e^{5am \epsilon}.
 \end{split}
 \end{equation*}
Combined with  \cref{lem:roughbd2a}, for $(\omega',x')\in F_{k}$ and for all $m \in \bad(\omega',x')$, 
 \begin{align}\# \what{N}^{*,\zeta}_{(\omega_0,x_0)}&\bigl(m ,\widetilde B(x',e^{-m })\cap K\bigr)
 \ge \frac{1}{2k} \exp(-(d^*-a {h}^*-3a\epsilon)m).\label{eq:16kCgooa}
 \end{align}

Let \[F_{k}'=\left\{(\omega',x')\in F_{k}: \lim_{r\to0} \frac{\log(\mu_{(\omega',x')}^{\eta^*}(\widetilde B(x',r)\cap \Gamma_\epsilon\cap K))}{\log r}=d^*\right\}.\] 
Since $F_k'\subset F_k\subset \what \Gamma\subset \Gamma_\epsilon$, \cref{prop:Lebdensity} ensures
we have
 $\mu^{\zeta}_{(\omega_{0},x_{0})}(F_{k}')=\mu^{\zeta}_{(\omega_{0},x_{0})}(F_{k})$.      
If $\mu^{\zeta}_{(\omega_{0},x_{0})}(F_{k})>0$, there would exist $(\omega_{0}',x_{0}')\in F_{k}$ such that \[\mu^{\eta^*}_{(\omega_{0}',x_{0}')}(F_{k}')=\mu^{\eta^*}_{(\omega_{0}',x_{0}')}(F_{k})=\mu^{\eta^*}_{(\omega_{0}',x_{0}')}(F_{k}'\cap \eta^*(\omega_0',x_0'))>0.\]
For any such $(\omega_{0}',x_{0}')$, we have \[(\pi_2)_* \mu^{\eta^*}_{(\omega_0',x_0')} \bigl(\pi_{2}(F_{k}'\cap \eta^{*}(\omega_0',x_0'))\bigr)\ge 
\mu^{\eta^*}_{(\omega_{0}',x_{0}')}(F_{k}'\cap \eta^*(\omega_0',x_0'))>0.\]

Fix $m_0\ge n_\epsilon $ to be chosen later.  
For each $x''\in  \pi_{2}(F_{k}'\cap\eta^{*}(\omega_{0}',x_{0}'))$, there exists $\omega''\in \Omega$ so that $(\omega'',x'')\in F_{k}'\cap \eta^{*}(\omega_{0}',x_{0}')$.
For each $x''\in  \pi_{2}(F_{k}'\cap\eta^{*}(\omega_{0}',x_{0}'))$, select some choice of $\omega''$ as above and select $m=m(\omega'',x'')$ such that 
\begin{enumerate}
\item\label{prop1a} $m(\omega'',x'') \in \bad(\omega'',x'')$; 
\item $m(\omega'',x'')\ge m_0$;
\item\label{prop3a} $\mu_{(\omega'',x'')}^{\eta^*}\bigl(\widetilde B(x'', 4  e^{-m(\omega'',x''){ +L}})\cap \Gamma_\epsilon\cap K\bigr)\le \exp(-m(\omega'',x'')(d^*-\epsilon)).$ 
\end{enumerate}

The image $$ \pi_{2}(F_{k}'\cap\eta^{*}(\omega_{0}',x_{0}'))= \pi_{2}(F_{k}'\cap\eta^{*}(\omega_{0}',x_{0}')\cap \Gamma_\epsilon)$$ is a subset of a Riemannian manifold $M$.  Cover the image
$ \pi_{2}(F_{k}'\cap\eta^{*}(\omega_{0}',x_{0}'))$ by Riemannian 
 balls centered at $x'' \in \pi_{2}(F_{k}'\cap\eta^{*}(\omega_{0}',x_{0}'))$ of radius $4\exp(-m(\omega'',x'')+L)$ where $m(\omega'',x'')$ satisfies \eqref{prop1a}--\eqref{prop3a} above.

By the  Besicovitch covering lemma, there is {$0<q$ depending only on the dimension of $M$ and} 
a countable sub-cover $\mathcal D'=\{B(i)=B(x''_i, { 4e^{-m_i+L}})\} $ where $m_i= m(\omega''_i,x_i'')$
such that every point in $M$ is contained in at most $q$ balls from $\mathcal D'$.  
 By the standard volume distortion arguments there is $C_d$ (say $C_d= 4^{\dim M}$) such that each $B(i)$ with $m_i=n$ intersects at most $C_d q$ distinct balls $B(j)$ with $m_j =n$.  
Then by  \cref{lem:approxa}\ref{approx4:stable-separation} or \ref{approx4:unstable-separation}, 
 \begin{equation}\label{besicov}
\sum_{m_i=n} \# {\what N}^{*,\zeta}_{(\omega_0',x_0')}\bigl(n,  \widetilde B(x_i'',  e^{-n})\cap K\bigr)
\le C_dq \# \what {N}^{*,\zeta}_{(\omega_0,x_0)}(n, 
K
 ). 
\end{equation}

Assembling, we have  
\begin{align*} 
\quad 0<&\mu^{\eta^*}_{(\omega_{0}',x_{0}')}(F_{k}) \\&{ \le \sum_i \mu_{(\omega_i'',x_i'')}^{\eta^*}\bigl(\widetilde B(x_i'',4e^{-m_i+L})\cap\Gamma_\epsilon\cap K\bigr)}
\\&\leq \sum_i e^{-m_i (d^*-\epsilon)} \le  \sum_i e^{-m_i (d^*-a\epsilon)}  
\\
&\le \sum_{n\ge m_0} \sum_{m_i=n} 2k\#\what{N}^{*,\zeta}_{(\omega_0,x_0)}\bigl(m_i ,\widetilde B(x_i'',e^{-m_i })\cap K\bigr)\exp(-a m_{i}({h}^*+2\epsilon))& \eqnote{   \cref{eq:16kCgooa}}\\ 
&\le C_dq\sum_{n\ge m_0}   2k\#\what{N}^{*,\zeta}_{(\omega_0,x_0)}(n , K) \exp(-a n({h}^*+2\epsilon))  & \eqnote{   \cref{besicov}}\\
&\le C_dq\sum_{n\ge m_0}  2k e^{a n({h}^*+\epsilon)}
\exp(-a n({h}^*+2\epsilon)) &\hfill \eqnote{Lem.~\ref{lem:roughbda}}\\
&= C_dq\sum_{n\ge m_0}  2k  e^{- a n\epsilon}.
\end{align*}
Taking our choice of  $m_0>0$ sufficiently large, the final expression can be made arbitrarily small, contradicting the assumption   $\mu^{\eta^*}_{(\omega_{0}',x_{0}')}(F_{k})>0$.
\end{proof}

\subsection{Lower bound on dimension}
With the previous comparisons, we show the lower bound in \cref{thm:prod}.

\begin{proposition}\label{prop:LB}
For $\mu$-a.e.\ $(\omega, x)\in X$, 
$$\liminf_{n\to \infty}\frac{\log\mu_{(\omega,x)}^\zeta(\wtd B(x,e^{-n}))}{\log(e^{-n})}\ge {d}^u+d^s.$$
\end{proposition}

\begin{proof}[Proof of \Cref{prop:LB}]
Fix $(\omega_0,x_0)\in \what \Gamma$ and $K=\partK(\omega_0,x_0)$. 
For
 $\mu_{(\omega_{0},x_{0})}^{\zeta}$-a.e.~$(\omega,x)\in \what \Gamma\cap K$ and all  $n\ge \hat n$ sufficiently large, write $A= \wtd B(x,e^{-n}) \cap\Gamma_\epsilon\cap K$.  We have 
\begin{align*}
\mu_{(\omega,x)}^{\zeta}&\bigl(\widetilde{B}(x,e^{-n} ) \cap\Gamma_\epsilon \cap K\bigr)
\\
& \le   \# N^{\zeta}_{(\omega,x)}\bigl(n, \widetilde B(x,e^{-n}) \cap K \bigr)\cdot e^{-an(h+\hat{h}-\epsilon)}
&\eqnote{\cref{lem:2.8}}\\
& \le   \#\widehat{N}^{s,\zeta}_{(\omega,x)}(n,A) 
\cdot \#\widehat{N}^{u,\zeta}_{(\omega,x)}(n,A) \cdot e^{-an(h+\hat{h}-\epsilon)}
&\eqnote{\cref{lem:2.9}}\\
& \le   \# {N}^{s,\zeta}_{(\omega,x)}(n,A) 
\cdot \# {N}^{u,\zeta}_{(\omega,x)}(n,A) \cdot e^{-an(h+\hat{h}-\epsilon)}\cdot {e^{10an \epsilon}}
&\eqnote{\cref{lem:compareN}}\\
& \le   
e^{n(-d^s+a h^s+2a\epsilon)}
\cdot e^{n(-d^u+a h^u+2a\epsilon)}
 \cdot e^{-an(h+\hat{h}-\epsilon)}e^{10an \epsilon}
&\eqnote{\cref{lem:2.10}}\\
&=e^{n(-d^s -d^u +15a\epsilon)}&\eqnote{$h^s=\hat h$, $h^u=h$}
\end{align*}
and so (applying \cref{prop:Lebdensity}) 
$$\liminf_{n\to \infty}\frac{\log\mu_{(\omega,x)}^\zeta(\wtd B(x,e^{-n}))}{\log(e^{-n})}\ge {d}^u+d^s-15a\epsilon.$$
The proposition now follows using that $\epsilon>0$ was arbitrary and that the associated sets $\what \Gamma$ exhaust a full-measure subset of $X$.  
\end{proof}

\end{document}